%% file: main.tex
\documentclass[11pt, letterpaper, hidelinks]{article}
\usepackage{amsmath,amsthm, amsfonts, amssymb}		
\renewcommand{\proofname}{\textit{\textbf{Proof}.}}

\makeatletter
\renewenvironment{proof}[1][\proofname]{\par
  \pushQED{\qed}%
  \normalfont \topsep6\p@\@plus6\p@\relax
  \trivlist
  \item[\hskip\labelsep
        #1\@addpunct{}]\ignorespaces %
}{%
  \popQED\endtrivlist\@endpefalse
}
\makeatother

\usepackage{amsmath,amsthm, amsfonts, amssymb}		
\usepackage{url}
\usepackage[margin=1in]{geometry}
\usepackage{mathtools}
\mathtoolsset{showonlyrefs}
\usepackage[dvipsnames]{xcolor}
\usepackage{hyperref}
\newcommand\myshade{85}
\colorlet{mylinkcolor}{violet}
\colorlet{mycitecolor}{YellowOrange}
\colorlet{myurlcolor}{Aquamarine}
\hypersetup{
  linkcolor  = mylinkcolor!\myshade!black,
  citecolor  = mycitecolor!\myshade!black,
  urlcolor   = myurlcolor!\myshade!black,
  colorlinks = true,
}
\usepackage{enumitem, comment, xifthen}
\usepackage{graphicx}
\usepackage{framed}
\usepackage{etoolbox}
\usepackage{tikz}
\usepackage{pgfplots}
\pgfplotsset{compat=1.18} 
\usetikzlibrary{math}
\usepackage{wrapfig}
\usepackage{authblk}
\usepackage{mathabx} %
\usepackage{subcaption}
\usepackage{xfrac}
\usepackage{inconsolata}
\usepackage{mathrsfs} %
\usepackage[mathcal]{euscript} %

\DeclareMathAlphabet{\mathfrak}{U}{jkpmia}{m}{it}
\graphicspath{{figs/}}
\DeclareMathAlphabet{\mathsf}{OT1}{qhv}{m}{n}
\input{src_tex/ams_defs}

\input{src_tex/reese_defs}
\input{paper_defs}
\newcommand{\LowerP}{\underline{p}}
\title{\Large \bfseries 
Revisiting mean estimation over $\ell_p$ balls: Is the MLE optimal?}
\author[1]{Liviu Aolaritei\footnote{Alphabetical author order.}\textsuperscript{,}}
\newcommand\CoAuthorMark{\footnotemark[\arabic{footnote}]\textsuperscript{,}} %
\author[1,2,4]{Michael I.\ Jordan\protect\CoAuthorMark}
\author[1,3]{Reese Pathak\protect\CoAuthorMark}
\author[2]{Annie Ulichney\protect\CoAuthorMark}
\affil[1]{Department of Electrical Engineering and Computer Sciences (EECS), University of California, Berkeley}
\affil[2]{Department of Statistics, University of California, Berkeley}
\affil[3]{School of Operations Research and Information Engineering (ORIE), Cornell University}
\affil[4]{Inria, Ecole Normale Sup\'erieure, PSL Research University}
\date{}

\begin{document}

\maketitle

\vspace{-5mm}
\begin{abstract}
We revisit the problem of mean estimation in the Gaussian sequence model with $\ell_p$ constraints for
$p \in [0, \infty]$.
We demonstrate two phenomena for the behavior of the maximum likelihood estimator (MLE), which depend on the noise level, the radius of the (quasi)norm constraint, the dimension, and the norm index $p$. First, if $p$ lies between $0$ and $1 + \Theta(\tfrac{1}{\log \dimension})$, inclusive, or if it is greater than or equal to $2$, the MLE is minimax rate-optimal for all noise levels and all constraint radii. 
On the other hand, for the remaining norm indices---namely, if $p$ lies between $1 + \Theta(\tfrac{1}{\log \dimension})$ and $2$---there is a more striking behavior: the MLE is minimax rate-\emph{sub}optimal, despite its nonlinearity in the observations, for essentially all noise levels and constraint radii for which nonlinear estimates are necessary for minimax-optimal estimation. 
Our results imply that when given $n$ independent and identically distributed Gaussian samples, the MLE can be suboptimal by a  polynomial factor in the sample size. Our lower bounds are constructive: whenever the MLE is rate-suboptimal, we provide explicit instances on which the MLE provably incurs suboptimal risk. Finally, in the non-convex case---namely when $p < 1$---we develop sharp local Gaussian width bounds, which may be of independent interest.
\end{abstract}
\vspace{5mm}

\section{Introduction}

This paper revisits the classical problem of mean estimation in the Gaussian sequence model in $\R^\dimension$, subject to  $\ell_p$ (quasi)norm constraints. 
Specifically, suppose we are given a Gaussian observation, 
\begin{equation}\label{eqn:lp-model}
Y = \thetastar + \sigma\xi, \quad \mbox{where} \quad \|\thetastar\|_p \leq \radius, \quad \mbox{and} \quad 
\xi \sim \Normal{0}{I_\dimension}.
\end{equation}
We aim to estimate $\thetastar$ to optimal mean squared error, depending on the dimension $\dimension$, the norm index $p \in [0, \infty]$, noise level $\sigma > 0$, and radius $\radius > 0$. 
This problem is most easily motivated by viewing the $\ell_p$ constraint as a form of prior knowledge on regularity of the (unknown) parameter $\thetastar$ that we seek to estimate. For instance, depending on the choice of $p$, these constraints correspond to the assumption that $\thetastar$ is sparse or that $\thetastar$ does not have very large entries. 

\subsection{The ``nonlinearity phenomenon''}
\label{sec:DJ-theory}
The seminal work of Donoho and Johnstone~\cite{DonJoh94} established a surprising \emph{nonlinearity phenomenon}.\footnote{This terminology is taken from the monograph~\cite[p.~378]{Joh19}} Specifically, they demonstrated a phase transition in the problem described above, with respect to the optimality of linear estimation procedures over the set $\Ball^\dimension_p$ of vectors in $\R^\dimension$ obeying the $\ell_p$ constraint $\|\thetastar\|_p \leq 1$.\footnote{We assume that the radius satisfies $\radius = 1$. Results for general radii are easy to obtain via rescaling.} 

\paragraph{Nonlinear and linear minimax risks:} 
To describe their result, recall the 
(nonlinear) minimax mean squared error (MSE) for the model~\eqref{eqn:lp-model}. It is given by
\begin{equation}\label{defn:minimax-risk}
\MinimaxRisk(\Ball^\dimension_p, \sigma) 
\defn 
\inf_{\hat \theta} 
\sup_{\thetastar \in  \Ball^\dimension_p} 
\E_{Y \sim \Normal{\thetastar}{\sigma^2 I_\dimension}} 
\Big[\|\hat \theta(Y) - \thetastar\|_2^2\Big]. 
\end{equation}
The infimum above is taken over all (measurable) functions $\hat \theta \colon \R^\dimension \to \R^\dimension$. In contrast, the \emph{linear} minimax MSE further restricts the infimum to mappings $Y \mapsto \hat \theta(Y)$ which are linear, 
\begin{equation}\label{defn:linear-minimax-risk}
\LinearMinimaxRisk(\Ball^\dimension_p, \sigma) 
\defn 
\inf_{\hat \theta~\text{linear}} \;
\sup_{\thetastar \in  \Ball^\dimension_p} 
\E_{Y \sim \Normal{\thetastar}{\sigma^2 I_\dimension}} 
\Big[\|\hat \theta(Y) - \thetastar\|_2^2\Big].
\end{equation}
Of course, we have $\MinimaxRisk \leq \LinearMinimaxRisk$---the supremum risk of an optimal linear estimator bounds the minimax MSE from above. The \emph{competitive ratio} for linear estimates, defined below, quantifies by how much the risk degrades when we restrict to linear estimators,
\begin{equation}\label{defn:linear-comp-ratio}
\LinCompRatio(\Ball^\dimension_p, \sigma) 
\defn \frac{\LinearMinimaxRisk(\Ball^\dimension_p, \sigma)}{\MinimaxRisk(\Ball^\dimension_p, \sigma)}.
\end{equation}
Indeed, when $\LinCompRatio < \infty$, it provides the smallest constant $C \geq 1$ 
such that $\LinearMinimaxRisk \leq C \; \MinimaxRisk$. That is, the competitive ratio quantifies the multiplicative gap in the error of the best linear estimate to that of a minimax optimal estimator, which is allowed to be (and often must be) nonlinear.

\subsubsection{The necessity of nonlinear estimates} 
\label{sec:necessity-of-nonlinear}
Donoho and Johnstone demonstrated~\cite[Theorem 1]{DonJoh94} that the competitive ratio for linear estimators satisfies
\begin{subequations}
\label{eqn:DJresults}
\begin{equation}\label{eqn:DJLimit-subopt}
\lim_{d \to \infty} 
\LinCompRatio(\Ball^\dimension_p, \sigma)  = \begin{cases} 
1, &  \mbox{if}~2 \leq p < \infty \\ 
\infty, & \mbox{if}~0 < p < 2
\end{cases}, 
\end{equation}
where in the limit above, the noise level is taken as $\sigma \equiv \sigma(d) \defn  C/\sqrt{d}$ with $C > 0$ being any fixed constant. Furthermore, by Theorem 5 in the paper~\cite{DonJoh94}, for $p \geq 2$, the competitive ratio satisfies, for sufficiently large $d$, 
\begin{equation}\label{eqn:DJ-optimal}
\LinCompRatio(\Ball^\dimension_p, \sigma) \asymp 1 
\quad \mbox{for all}~\sigma > 0. 
\end{equation}
\end{subequations}
Simply put, the results above demonstrate the following dichotomy: 
if $p \geq 2$, then by relation~\eqref{eqn:DJ-optimal}, linear estimates are---disregarding  constant factors---minimax rate-optimal. On the other hand, if $p < 2$, then the limit relation~\eqref{eqn:DJLimit-subopt} provides a sequence of noise levels for which linear estimates are suboptimal by a multiplicative factor that diverges as the dimension grows. 

In fact, the suboptimality indicated by the limit relation~\eqref{eqn:DJLimit-subopt}---in the case the norm index satisfies $p\in (0, 2)$---is even more striking. Indeed, it holds for a wide variety of scalings of the noise level. 
To be more precise, Theorem 5 in the paper~\cite{DonJoh94} and its strengthening via the results in the monograph~\cite{Joh19}, refine the relations~\eqref{eqn:DJresults} as follows.
First, in sufficiently low noise or high noise, linear estimates are order-optimal. More precisely, 
the relation~\eqref{eqn:DJ-optimal} extends to $p \in (0, 2)$ in that
\begin{subequations}
\label{eqn:limit-relations-versus-linear}
\begin{equation}\label{eqn:limit-vs-linear-when-order-optimal}
\limsup_{d \to \infty} 
\LinCompRatio(\Ball^\dimension_p, \sigma) \asymp 1, 
\quad \mbox{if}~\sigma \ll \frac{1}{d^{1/p}}~\mbox{or}~\sigma \gg \frac{1}{\sqrt{\log d}}. 
\end{equation}
In contrast, for a large variety of ``moderate'' noise levels, linear estimates are increasingly suboptimal as the dimension grows. More precisely,
the suboptimality result~\eqref{eqn:DJLimit-subopt} extends to other noise levels in the case $p \in (0, 2)$ as follows,
\begin{equation}\label{eqn:limit-vs-linear-when-sub-optimal}
\lim_{d \to \infty} 
\LinCompRatio(\Ball^\dimension_p, \sigma)  = \infty, \quad \mbox{if}~\sigma \gg \frac{1}{d^{1/p}}~\mbox{and}~\sigma \ll \frac{1}{\sqrt{\log d}}. 
\end{equation}
\end{subequations}
To establish these results, Donoho and Johnstone show that the soft thresholding (ST) estimator, given below with its carefully chosen  regularization level, achieves the minimax MSE \cite[Theorem 3]{DonJoh94}
\begin{equation}\label{eqn:DJ-soft-thresholding}
\hat \theta^{\,\mathsf{ST}}(Y) \defn 
\argmin_{\vartheta \in \R^\dimension} 
\Big\{ \, \|\vartheta - Y\|_2^2 + 2 \lambda \|\vartheta\|_1\, \Big\} \quad \mbox{with}\quad \lambda \defn \sqrt{2 \sigma^2 \log(\dimension \sigma^{p})}.
\end{equation}
Notably, the cases $\sigma \ll \dimension^{-1/p}$ or $\sigma \gg (\log \dimension)^{-1/2}$ are uninteresting for the $\ell_p$ mean estimation problem in that 
attaining the minimax rate is ``easy.'' More rigorously, in the former regime, $\hat \theta(Y) = Y$ is order-optimal, and in the latter regime, $\hat \theta(Y) = 0$ is order-optimal. To emphasize, in these cases, knowledge of the $\ell_p$ constraint does \emph{not} meaningfully impact the rate of estimation. 
In summary, the essential takeaway from these results of Donoho and Johnstone is: 
\begin{center} 
\emph{If $p \in (0, 2)$, then nonlinear estimates are required to achieve the minimax rate \\ whenever knowledge of the $\ell_p$ constraint can be leveraged to improve the worst-case MSE.}
\end{center}

\subsection{Overview of our results}

This paper takes as a starting point the Donoho-Johnstone nonlinearity phenomenon, as described in the previous section. We study arguably the most natural nonlinear estimator for the $\ell_p$ observational model~\eqref{eqn:lp-model}: the maximum likelihood estimator (MLE).%
\footnote{Here we assume $\radius = 1$; extensions to general radii are given in Section~\ref{sec:gen-radii}.} Concretely, the MLE is given by the Euclidean projection onto $\Ball^\dimension_p$:\footnote{For any $p$, a projection exists as $\Ball^d_p$ is closed and $\vartheta \mapsto \|\vartheta - Y\|_2^2$ is continuous and coercive. For $p \geq 1$, the projection is unique by the convexity of $\Ball^d_p$. For $p \in [0, 1)$, it is almost surely unique by a theorem of Erdös~\cite{Erd45}.}
\begin{equation}\label{eqn:single-sample-mle}
\Pi_{\Ball^\dimension_p}(Y) 
\defn \argmin_{\vartheta \in \Ball^\dimension_p} 
\Big\{\, \|\vartheta - Y\|_2^2 \,\Big\}.
\end{equation}
This estimator, unlike the soft-thresholding estimator, is ``intrinsic''---it obviously depends on the geometry of the underlying constraint set. Moreover, it does not require explicit specification of the noise level $\sigma > 0$. We emphasize that since projection is a nonlinear map, the Donoho-Johnstone theory---as described in Section~\ref{sec:DJ-theory}---\emph{does not} rule out the minimaxity, nor does it establish the rate-optimality of the MLE.
We are, therefore, left with the following question: 
\begin{center}
    \emph{For which $(\sigma, p, \dimension)$ is the maximum likelihood estimator (MLE) minimax optimal?}
\end{center}
There are many additional motivations for studying the MLE. For instance, classical asymptotic theory demonstrates that the MLE is asymptotically optimal for sufficiently regular parametric models; see Chapter 7 in~\cite{van00} as well as the discussion in Section~\ref{sec:prior-work}  below. Unfortunately, these conditions do not apply to our model~\eqref{eqn:lp-model}.

To address the question regarding the optimality of the MLE for the $\ell_p$ observational model, we take a nonasymptotic approach. Indeed, in a slight departure from the asymptotic methods of Donoho and Johnstone, we do not immediately take the dimension ${\dimension \to \infty}$. We characterize the worst-case mean squared error of the MLE, allowing arbitrary dependence in the parameters $(\sigma, p, \dimension)$. Recall that the worst-case MSE of the MLE is given by
\begin{equation}\label{eqn:worst-case-MLE-risk}
\MLERisk(\Ball^\dimension_p, \sigma) 
\defn 
\sup_{\thetastar \in \Ball^\dimension_p} 
\E_{Y \sim \Normal{\thetastar}{\sigma^2 I_\dimension}}
\Big[\|\Pi_{\Ball^\dimension_p}(Y) - \thetastar\|_2^2 \Big].
\end{equation}

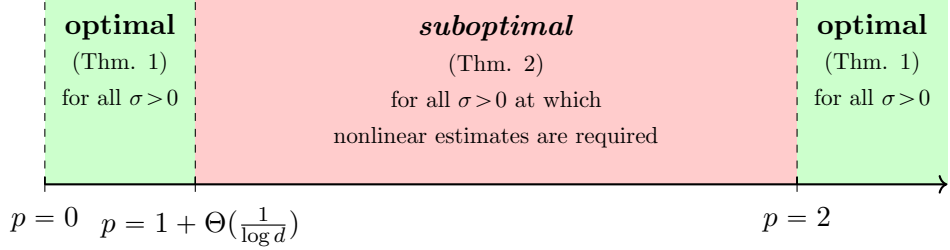
\begin{figure}
\centering
\input{figs/f-plot-p-range}
\caption{The choices of $p \in [0, \infty]$ for which MLE is order-optimal or suboptimal. If ${p \in [0, {1 + \Theta(\tfrac{1}{\log d})}]}$ or $p \geq 2$, then Theorem~\ref{thm:cases-where-mle-is-optimal} demonstrates that the MLE is optimal for every noise level $\sigma >0$. On the other hand, for intermediate $p$, Theorem~\ref{thm:cases-where-mle-is-suboptimal} demonstrates that the MLE is minimax suboptimal for essentially every noise level at which knowledge of the $\ell_p$ constraint, $\|\thetastar\|_p \leq 1$, meaningfully improves the rate of estimation.}
\label{fig:plot-of-p-range}
\end{figure}

\noindent Analogously to the case of linear estimates, we also define the \emph{competitive ratio for the MLE} as
\begin{equation}\label{defn:comp-ratio-MLE}
\MLECompRatio(\Ball^\dimension_p, \sigma) 
\defn \frac{\MLERisk(\Ball^\dimension_p, \sigma)}{\MinimaxRisk(\Ball^\dimension_p, \sigma)}.
\end{equation}
Indeed, when $\MLECompRatio < \infty$, the competitive ratio provides the smallest constant $C \geq 1$ such that the sandwich relation $\MinimaxRisk \leq \MLERisk \leq C \, \MinimaxRisk$ holds between the supremum risk of the MLE~\eqref{eqn:worst-case-MLE-risk} and the minimax MSE as defined in display~\eqref{defn:minimax-risk}.

Our main results, Theorems~\ref{thm:cases-where-mle-is-optimal} and~\ref{thm:cases-where-mle-is-suboptimal}, characterize when the MLE is order-optimal (\ie when the competitive ratio satisfies $\MLECompRatio(\Ball^\dimension_p, \sigma) \asymp 1$) and when the MLE is rate-suboptimal (\ie when $\MLECompRatio(\Ball^\dimension_p, \sigma) \gg 1$), respectively.
Figure~\ref{fig:plot-of-p-range}, shown above, illustrates the dependency of this characterization on the underlying norm index $p$. 
Our results indicate that in the same setting in which Donoho and Johnstone demonstrate the rate-suboptimality of linear estimates, the MLE may either be  minimax rate-suboptimal---if $p$ is strictly between $1 + \Theta(\tfrac{1}{\log \dimension})$ and $2$---or minimax rate-optimal if $p \in [0, 1 + \Theta(\tfrac{1}{\log \dimension})]$. 

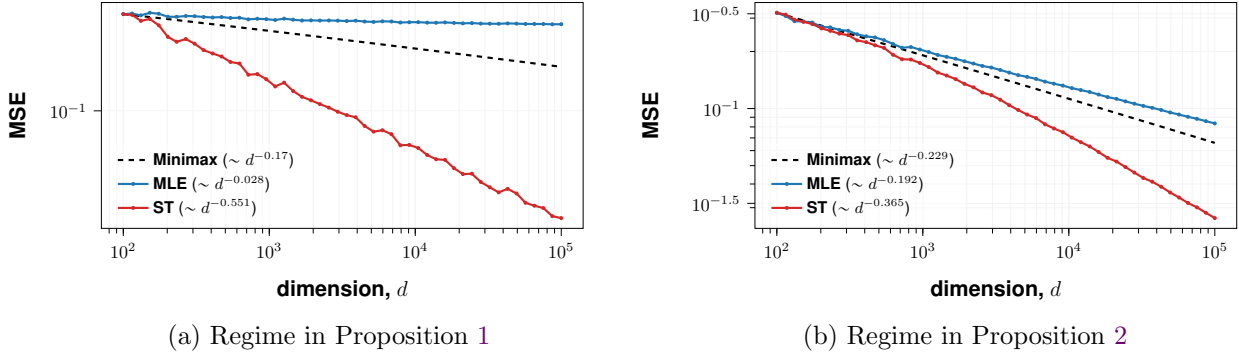
\begin{figure}
  \centering
  \begin{subfigure}[b]{0.48\textwidth}
    \centering
    \hspace{-10mm}
    \scalebox{0.68}{\input{figs/regime1}}         
    \caption{Regime in Proposition~\ref{prop:suboptimality-of-erm-easy}}
    \label{fig:regime1}
  \end{subfigure}%
  \hfill
  \begin{subfigure}[b]{0.48\textwidth}
    \centering
    \scalebox{0.68}{\hspace{-10mm}\input{figs/regime2}}
    \caption{Regime in Proposition~\ref{prop:suboptimality-of-erm-hard}}
    \label{fig:regime2}
  \end{subfigure}
\bigskip
  \caption{%
    Simulated mean squared error (MSE) of the MLE and soft‐thresholding (ST)
    estimators compared to the minimax rate, for $p=1.5$. Panel~(a) depicts a regime with the ground truth taken as
    $\thetastar=e_1$, and corresponds to the result in Proposition~\ref{prop:suboptimality-of-erm-easy} given below. Panel~(b) depicts a regime with the ground truth taken as $\thetastar=k^{-1/p}(\1_k,0_{d-k})$, and corresponds to the result in Proposition~\ref{prop:suboptimality-of-erm-hard} given below. Both panels indicate situations where the MLE is minimax rate-suboptimal. Further details on the simulations are provided in
    Appendix~\ref{sec:experiment-details}.%
  }
  \label{fig:simulations}
\end{figure}

Interestingly, in the case that $p \in (1 + \Theta(\tfrac{1}{\log \dimension}), 2)$, the MLE is rate suboptimal for essentially the same noise levels for which Donoho and Johnstone demonstrate the rate-suboptimality of linear estimators: when $\sigma >0$ satisfies $\sigma^2 \gg d^{-2/p}$ and 
$\sigma^2 \ll \log^{-1} \dimension$. 
However, our results show that there is a second phenomenon: depending on how small the noise level is, the MLE can potentially even have non-vanishing error as the dimension grows. 
To give a sense of this dichotomy, we present, in Figure~\ref{fig:simulations}, the two behaviors of the MLE when it is suboptimal. These behaviors depend on the configuration of the triple $(\sigma, p, d)$. The first of the two behaviors is that either the error can be constant (even as $d\to \infty$): this is depicted in Figure~\ref{fig:regime1}, and corresponds to the result presented below as Proposition~\ref{prop:suboptimality-of-erm-easy}. The second behavior is that for certain $(\sigma, p, d)$ the error of the MLE can be vanishing, but still minimax rate-suboptimal: this is depicted in Figure~\ref{fig:regime2} and corresponds to the result presented below as Proposition~\ref{prop:suboptimality-of-erm-hard}.

\subsection{Organization and notation}
\paragraph{Organization:} 
Formal statements of our main results are presented in Section~\ref{sec:main-results}. We then provide some extensions of these main results and a comparison to prior work in Section~\ref{sec:discussion}. Most of the proofs of the main results are presented in Section~\ref{sec:proofs}. We defer the proofs of some auxiliary results to the appendices.  

\paragraph{Notation:}
We collect the (mostly standard) notation used in the rest of the paper. 
The symbols $\lesssim$, $\gtrsim$, and $\asymp$ denote standard asymptotic notation: for $f(d), g(d)$, we write $f(d) \lesssim g(d)$ (resp. $f(d) \gtrsim g(d)$)  if there exists an absolute constant $C > 0$ such that $|f(d)| \leq C |g(d)|$ (resp. $|f(d)| \geq C |g(d)|$), and we write $f(d) \asymp g(d)$ if both $f(d) \lesssim g(d)$ and $f(d) \gtrsim g(d)$. For a set $\Theta \subset \R^d$, $\Pi_\Theta(y) = \argmin_{\vartheta \in \Theta} \| y - \vartheta \|_2^2 $ denotes the Euclidean projection operator.  The $\ell_p$-norm of $x \in \R^d$ is $\|x\|_p = (\sum_{i=1}^d |x_i|^p)^{1/p}$ for $p \in [1, \infty)$, with $\|x\|_{\infty} = \max_{1 \leq i \leq d} |x_i|$. For $p \in (0, 1)$ the $\ell_p$-quasinorm is defined by the same formula, $\|x\|_p = (\sum_{i=1}^d |x_i|^p)^{1/p}$, and we denote $\|x\|_0$ to be the number of nonzero entries in $x$. For $p \in (1, \infty)$, we denote by $q$ its Hölder conjugate, satisfying $q \defn \frac{p}{p-1}$, \ie such that $1/p + 1/q = 1$. We also consider the $(p, q)$ pairs $(1, \infty), (\infty, 1)$ to be conjugate. We write $X \sim \Normal{\mu}{\Sigma}$ for a Gaussian random vector with mean $\mu$ and covariance $\Sigma$, and $\Phi(t) \defn \P\{Z \leq t\}$ for $Z \sim \Normal{0}{1}$. The $i$th standard basis vector is $e_i$.
For two sets $A, B \subset \R^\dimension$ we define the Minkowski sum $A + B = \{a + b : a \in A, b \in B\}$. For $\tau \in \R$ and $C \subset \R^\dimension$ we also define the dilate $\tau C = \{\tau x: x \in C\}$. For a centrally symmetric, compact, convex set with nonempty interior $K \subset \R^\dimension$, we denote by $K^\circ$ its associated polar body. We occasionally use the shorthand notations $\twomin{a}{b} = \min\{a, b\}$ and $\twomax{a}{b} = \max\{a, b\}$ for the minimum and maximum, respectively, of $a, b \in \R$. We denote the nonnegative reals by $\R_+$ and the (strictly) positive reals by $\R_{++}$.

\section{Main results}
\label{sec:main-results}

In this section, we present our main results. 
In Section~\ref{sec:main-results-risk-characterization}, we characterize the cases where the MLE has 
order-optimal worst-case risk (Theorem~\ref{thm:cases-where-mle-is-optimal}) and those where the MLE has suboptimal worst-case risk (Theorem~\ref{thm:cases-where-mle-is-suboptimal}).
In Section~\ref{sec:closer-look-failure-modes}, we provide a more detailed look at the failure modes of the MLE.
In Section~\ref{sec:computing-the-MLE}, we briefly discuss the computational aspects surrounding the MLE.

\subsection{Characterizing the risk of the maximum likelihood estimator}
\label{sec:main-results-risk-characterization}

To begin with, we recall the minimax rates of estimation over $\ell_p$ balls. 
\paragraph{The minimax rate for the $\ell_p$ ball $\Ball_p^\dimension$:} 
From prior work~\cite[Theorem 11.7]{Joh19}, it is known that the minimax rate for the unit ball in $\ell_p^\dimension$ satisfies the following
relations:
\begin{subequations}\label{eqn:minimax-risk-p}
\begin{enumerate}[label=(\roman*)]
\item for $p \in [2, \infty]$,  
\begin{equation}\label{eqn:minimax-risk-p-geq-2}
\MinimaxRisk(\Ball^\dimension_p, \sigma) \asymp \twomin{d^{1 - 2/p}}{(\sigma^2 d)}, 
\end{equation}
\item for $p \in (0, 2)$,  
\begin{equation}\label{eqn:minimax-risk-p-1-2}
\MinimaxRisk(\Ball^\dimension_p, \sigma) 
  \asymp 
  \begin{cases}
  1 & \sigma^2 \geq \tfrac{1}{1 + \log \dimension} \\ 
  (\sigma^2 \log(\e \dimension \sigma^p))^{1-p/2} &  \sigma^2 \in \left[ \tfrac{1}{d^{2/p}}, \tfrac{1}{1 + \log \dimension} \right] \\ 
  \sigma^2 \dimension & \sigma^2 \leq \tfrac{1}{\dimension^{2/p}}
  \end{cases}, \quad \mbox{and,} 
\end{equation}
\item for $p = 0$ and any $s \in [d]$, 
\begin{equation}
\label{eqn:minimax-risk-p-0}
\MinimaxRisk\big(\Ball^\dimension_0(s), \sigma\big) 
  \asymp 
\sigma^2 s \log \Big(\frac{\e d}{s}\Big). 
\end{equation}
\end{enumerate}
\end{subequations}
For completeness, we provide a proof of this result, primarily to justify that the implicit constant is independent of $p$, in Appendix~\ref{sec:minimax-rate-lp}.\footnote{Throughout this paper, we have provided explicit constants wherever possible. The constants are either given in the statements of the results or at the beginning of the associated proofs. The primary purpose is to show that these constants are universal; they do not depend on $p, d, \sigma, \radius$. We have not made any effort to optimize these constants.}

Our first main result characterizes the cases where the MLE has order-optimal worst-case risk. Specifically, 
we determine the situations when the worst-case risk of the MLE coincides with the minimax rate, given above. 

\btheo [Situations when the MLE has order-optimal worst-case risk] \label{thm:cases-where-mle-is-optimal}
The competitive ratio for the risk of the maximum likelihood estimator (MLE), as defined in display~\eqref{defn:comp-ratio-MLE}, satisfies 
\begin{equation}\label{eqn:comp-ratio-limit-optimal}
\MLECompRatio(\Ball^\dimension_p, \sigma) \asymp 1, 
\end{equation}
provided that the parameters $(\sigma, p, d)$ satisfy:
\begin{enumerate}[label=(\roman*)]
\item \label{case:p-geq-2}
$p \geq 2$, for arbitrary $\sigma > 0, d \geq 1$; or  
\item \label{case:p-near-1}
$p \in (0, 1 + \tfrac{1}{1 + \log \dimension}]$, for some $d \geq 1$, and for arbitrary $\sigma > 0$; or 
\item \label{case:p-1-2-high-or-low-noise}
$p \in (1 + \tfrac{1}{1 + \log \dimension}, 2)$, for  $(\sigma, d)$ which  satisfy $\sigma \leq \tfrac{1}{d^{1/p}}$ or $\sigma \geq \tfrac{1}{\sqrt{1+\log d}}$; or 
\item \label{case:p-0}
$p = 0$, for any $\sigma > 0$, and taking in display~\eqref{eqn:comp-ratio-limit-optimal} $\Ball^\dimension_0 = \Ball^\dimension_0(s)$ for any $s \in [d]$.
\end{enumerate}
Indeed, in any of these cases, the MLE has minimax order-optimal worst-case risk.
\etheo 
\noindent See Section~\ref{sec:proof-main-optimal-result} for a proof of this result.

Our second main result essentially demonstrates that when the cases covered by Theorem~\ref{thm:cases-where-mle-is-optimal} do not hold, the MLE is minimax rate-suboptimal.
More concretely, define the threshold norm index 
\begin{equation}\label{eqn:lower-limit-on-p}
\LowerP(d) \defn 1 + \frac{1}{\log d + 1}.
\end{equation}
Then, for $p \geq 1$, the missing case in Theorem~\ref{thm:cases-where-mle-is-optimal} is when the following conditions hold:
\begin{equation}\label{eqn:not-covered-cases}
p \in \big(\LowerP(d), 2\big), \quad \mbox{and} 
\quad 
\sigma \in \Big(\frac{1}{d^{1/p}}, \frac{1}{\sqrt{\log d + 1}}\Big).
\end{equation}
Indeed, when these conditions hold, we will show that the MLE is minimax suboptimal. For a succinct statement, here we consider the case when $d \to \infty$, although we emphasize that the results underlying Theorem~\ref{thm:cases-where-mle-is-suboptimal} hold nonasymptotically; see for instance Propositions~\ref{prop:suboptimality-of-erm-easy} and~\ref{prop:suboptimality-of-erm-hard} below. 
In these asymptotics, we consider two scenarios. First, we take $p \in (1, 2)$ as fixed and independent of the dimension. Second, we allow $p$ to depend on dimension so that $p = p(d)$ with 
\begin{equation}\label{eqn:asymptotic-conditions-for-p}
p(d) \gg \LowerP(d)   
\quad \mbox{and} \quad 
\limsup_{d \to \infty}
p(d) < 2.
\end{equation}
These conditions simply separate $p = p(d)$ from the endpoints for the norm index, as given in display~\eqref{eqn:not-covered-cases}.

\btheo [Situations when the MLE is suboptimal] \label{thm:cases-where-mle-is-suboptimal}
The competitive ratio for the risk of the MLE, as defined in display~\eqref{defn:comp-ratio-MLE}, satisfies
\[
\lim_{d \to \infty} \MLECompRatio(\Ball^\dimension_p, \sigma) = \infty,
\] 
i.e., the MLE is rate-suboptimal by arbitrarily large factors as the dimension grows, in the following situations:
\begin{enumerate}[label=(\roman*)]
    \item \label{thm:suboptimal-p-independent} if $p \in (1, 2)$, independent of the dimension $d$,
    for any sequence of noise levels such that $\sigma^p d \gg 1$ and $\sigma^2 \log d \ll 1$; and, 
    \item \label{thm:suboptimal-p-dependent}
    if $p = p(d)$ and the conditions~\eqref{eqn:asymptotic-conditions-for-p} hold,   
    for any sequence of noise levels such that 
    $\sigma^p d \gg d^{2-p} (1 - \tfrac{1}{p})^{p/2}$ and $\sigma^2 \log d \ll 1$. 
\end{enumerate}
\etheo  
\noindent See Section~\ref{sec:proof-main-suboptimal-result} for a proof of this result.

We now unpack Theorem~\ref{thm:cases-where-mle-is-suboptimal} slightly. 
Define the interval
\begin{equation}\label{eqn:suboptimality-interval}
\Sigma(p, d) 
\defn 
\Big(\frac{1}{d^{1/p}}, \frac{1}{\sqrt{1 + \log \dimension}}\Big).
\end{equation}
For a fixed $p \in (1, 2)$ and sufficiently large dimension $d$, we saw in Section~\ref{sec:DJ-theory} that, for $\sigma \in \Sigma(p, d)$, nonlinear estimators are required to attain the minimax rate. Unfortunately, despite the nonlinearity of the MLE in the observation, Theorem~\ref{thm:cases-where-mle-is-suboptimal} illustrates that the MLE is minimax rate-suboptimal for such $\sigma$. Indeed, Theorem~2 says for essentially\footnote{Specifically, for all $\sigma \in \Sigma(p, d)$ except near the endpoints, for which Theorem~\ref{thm:cases-where-mle-is-optimal} guarantees order-optimality.} any such $\sigma$, the factor by which the worst-case error of the MLE exceeds the minimax risk can be made arbitrarily large, assuming the dimension is sufficiently large. 
We emphasize that this result is perhaps counterintuitive: the MLE is nonlinear,  shrinks in a way which depends on the underlying convex constraint set, and 
yet suffers from essentially the same rate-suboptimality that Donoho and Johnstone illustrated for \emph{linear} procedures.

\subsection{A closer look at the failure modes of the MLE}
\label{sec:closer-look-failure-modes}

In this section, we take a closer look at the situations in which the MLE is minimax rate-suboptimal, as described in Theorem~\ref{thm:cases-where-mle-is-suboptimal}. 
Here, we return to a nonasymptotic setting where we allow for arbitrary dependence among the parameters $(\sigma, p, d)$. 

At a high level, the lower bounds stated in the next two propositions show that the MLE suffers when estimating sufficiently sparse signals. 
It turns out that the interval of suboptimal noise levels~\eqref{eqn:suboptimality-interval} can be further subdivided---disregarding constant factors---into essentially two subintervals 
\begin{equation}\label{eqn:two-intervals}
\Sigma^{(1)}(p, d) = \Big(\frac{1}{\sqrt{q} d^{1/q}}, \frac{1}{\sqrt{1 + \log \dimension}}\Big) 
\quad \mbox{and} \quad 
\Sigma^{(2)}(p, d) =
\Big(\frac{1}{d^{1/p}}, \frac{1}{d^{1/q}}\Big).
\end{equation}
Our results show that the suboptimality behavior is different in the two cases: the risk is different in each case, and the instance on which we demonstrate the suboptimality of the MLE also differs in each case. 
First, as shown next in 
Proposition~\ref{prop:suboptimality-of-erm-easy}, in the case $\sigma \in \Sigma^{(1)}(p, d)$, the MLE suffers on a parameter with a single nonzero entry. 

\bpr 
\label{prop:suboptimality-of-erm-easy}
Suppose that the triple $(\sigma, d, p)$ satisfies
\begin{equation}\label{ineq:sigma-d-p-condition}
p \in \Big[1 + \frac{1}{1 + \log d}, 2\Big], \quad 
\sigma \geq \frac{10}{\sqrt{q} d^{1/q}}, \quad \mbox{and} \quad d \geq 45.
\end{equation}
Then, the worst-case risk of the maximum likelihood estimator satisfies the lower bound 
\[
  \MLERisk(\Ball^\dimension_p, \sigma) \gtrsim 1.
\]
Moreover, in this case, the lower bound is attained at $\thetastar = e_1$.
\epr 
\noindent See Section~\ref{sec:proof-prop-easy} for a proof. 

In the case that $\sigma \in \Sigma^{(2)}(p, d)$ as defined in the display~\eqref{eqn:two-intervals}, the MLE suffers on a $k = k(\sigma, d, p)$-sparse signal, as shown in the next result. The proof illustrates that for $\sigma \in \Sigma^{(2)}(p, d)$, we must take $k > 1$.

\bpr \label{prop:suboptimality-of-erm-hard}
Fix $p \in (1 + \tfrac{1}{\log \dimension + 1}, 2)$. Suppose that the triple $(\sigma, p, d)$
satisfies 
\[
  \sigma \geq \frac{1}{d^{1/p}} \quad \mbox{and} \quad d \geq 4.
\]
Then, the worst-case risk of the maximum likelihood estimator satisfies the lower bound
\[
\MLERisk(\Ball^\dimension_p, \sigma) \gtrsim \twomin{1}{\sigma d^{1/q}}.
\]
\epr 
\noindent See Section~\ref{sec:proof-prop-hard} for a proof. We note that the lower bound in this case is attained with $\thetastar$ which has $k = k(\sigma, d, p)$-nonzero entries equal to $k^{-1/p}$, and thus $\thetastar \in \Ball^\dimension_p$. The choice of number of nonzero entries $k(\sigma, d, p)$ is given explicitly in the proof. 

\subsubsection{Overview of the proofs of Propositions~\ref{prop:suboptimality-of-erm-easy} and~\ref{prop:suboptimality-of-erm-hard}}

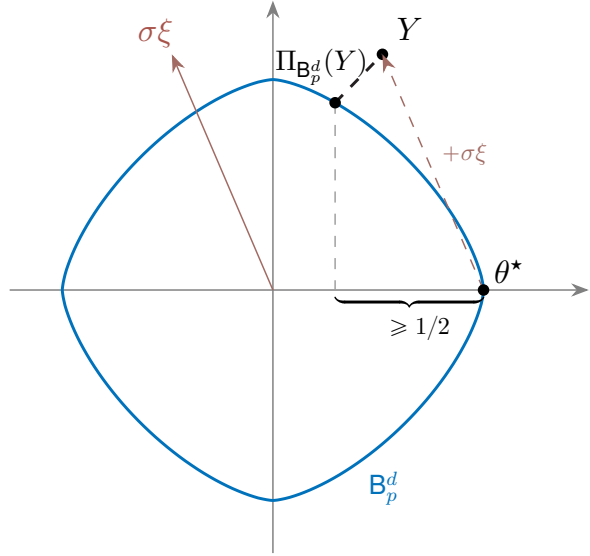
\begin{figure}[t]
    \centering
    \scalebox{1.4}{\input{figs/f-prop1-fig}}
    \caption{Depiction of the main idea underlying Proposition~\ref{prop:suboptimality-of-erm-easy}, illustrated for $d = 2$ and $p = 1.5$.}
    \label{fig:proof-schematic}
\end{figure}

We show that the MLE $\Pi_{\Ball^\dimension_p}(Y)$ incurs suboptimal risk when estimating sparse signals on the boundary of the $\ell_p$ ball. In such cases, the MLE shrinks towards denser vectors to fit the noise, which, due to the $\ell_p$ constraint, forces it to also shrink coordinates corresponding to the true signal. This tradeoff leads to estimation error exceeding the minimax rate. See Figure~\ref{fig:proof-schematic} for a schematic illustration of this phenomenon.

In order to quantify this level of shrinkage, we show that the maximum likelihood estimator is a form of data-dependent shrinkage estimator.
To be more precise, fix $\lambda, t > 0$. Let $\Psi_{p, \lambda}(t)$ denote the unique positive solution, in $\psi$, to the equation, 
\[
\psi + \lambda \psi^{p-1} = t. 
\]
Observe that $\Psi_{p,\lambda}(t)$ is precisely the proximal operator~\cite{Roc70} of the function $x \mapsto \tfrac{\lambda}{p} x^p$ evaluated at $t \geq 0$:
\[
\Psi_{p,\lambda}(t)
= \operatorname{\bf prox}_{\frac{\lambda}{p}(\cdot)^p}(t)
= \argmin_{x} \left\{ \tfrac{1}{2}(x - t)^2 + \tfrac{\lambda}{p} x^p \right\}.
\]
Viewing $\lambda$ as fixed, $\Psi_{p, \lambda}$ can be interpreted as a nonlinear shrinkage operator applied to $t$. As we show next, the MLE can then be interpreted, coordinatewise, in terms of this operator evaluated at $\lambda = \lambda^\star$, where $\lambda^\star$ is data-dependent. The following lemma will be instrumental in proving the suboptimality of the MLE.

\ble [Variational characterization of the MLE, $p \in (1, \infty)$]
\label{lem:var-characterization}
Suppose that $p \in (1, \infty)$, and take the parameter space $\Theta$ to be the $\ell_p$ ball $\Ball^\dimension_p$. 
Then, the maximum likelihood estimator satisfies the coordinatewise relations
\begin{equation}\label{eqn:coordinate-characterization}
\Pi_{\Theta}(Y)_i = \sign(Y_i) \, \Psi_{p, \lambdastar}\big(|Y_i|\big), \quad \mbox{for}~i=1,\ldots,\dimension.
\end{equation}
Here, $\lambdastar$ is defined by 
\begin{equation}\label{eqn:minimal-lambda}
\lambdastar = \inf\left\{\, \lambda \geq 0 \mid \sum_{i=1}^\dimension \Psi_{p, \lambda}\big(|Y_i|\big)^p \leq 1 \,\right\}.
\end{equation}
\ele 

We establish Lemma~\ref{lem:var-characterization} via convex duality. We now provide an overview of the ideas underlying the proofs of Propositions~\ref{prop:suboptimality-of-erm-easy} and \ref{prop:suboptimality-of-erm-hard}.

\paragraph{General proof idea:} Proposition~\ref{prop:suboptimality-of-erm-easy} and Proposition~\ref{prop:suboptimality-of-erm-hard} both follow from a two-step argument outlined below. In each case, we construct a particular sparse signal of the form $\thetastar = \frac{1}{k^{1/p}}(\1_k, 0_{d-k})$, where $k < \frac{d}{2}$ depends on $(\sigma, p, d)$. The key idea in both proofs is to leverage the coordinatewise condition,
\begin{equation} \label{eq:lagrange-fixedpoint}
|Y_i| = |\Pi_{\Ball^\dimension_p}(Y)_i| + \lambda^\star |\Pi_{\Ball^\dimension_p}(Y)_i|^{p-1}, \quad \mbox{for all}~i=1,\ldots, d,
\end{equation}
derived from Lagrange duality (see the proof of Lemma~\ref{lem:var-characterization}). These conditions allow us to relate the magnitude of the ``noise-only'' coordinates of $Y$ to a lower bound on the Lagrange multiplier, $\lambda^\star$, which controls the amount of shrinkage. In turn, this determines the estimation error on the ``signal-plus-noise'' coordinates. Our argument proceeds in two steps:
\begin{itemize}
\item[(i)] We establish that $\lambdastar$ is large, using the ``noise-only'' coordinates (\ie the last $(d - k)$ entries of $Y$); and,
\item[(ii)] By (i), as $\lambdastar$ is large, we show that the MLE must incur large estimation error in the ``signal-plus-noise'' coordinates (\ie the first $k$ entries of $Y$).
\end{itemize}
We now provide further details on the specific construction and reasoning used in each proof.

\paragraph{Proof sketch for Proposition~\ref{prop:suboptimality-of-erm-easy}:} The relatively large noise level (i.e., $\sigma \geq \frac{c}{\sqrt{q} d^{1/q}}$, for a constant $c$) allows the lower bound to be attained with the simple choice $k = 1$, resulting in $\thetastar = e_1$. If $\lambda^\star \geq 1$ and $Y_1 \leq 1$, we automatically have $\Pi_{\Ball^\dimension_p}(Y)_1 \leq \frac{1}{2}$, which implies from the fixed point equation~\eqref{eq:lagrange-fixedpoint} that $(\thetastar_1 - \Pi_{\Ball^\dimension_p}(Y)_1)^2 \geq \frac{1}{4}$. This is precisely what we will prove. To see this, we first apply \eqref{eq:lagrange-fixedpoint} and some algebraic manipulations to show that $\lambda^\star = \|Y - \Pi_{\Ball^\dimension_p}(Y)\|_q$. Since the noise is sufficiently large, we see that with large probability $Y \not \in \Theta(1) \Ball^\dimension_q$, but by definition,  $\Pi_{\Ball^\dimension_p}(Y) \in \Ball^\dimension_p \subset \Ball^\dimension_q$,
which yields $\lambdastar \gtrsim 1$ by the triangle inequality. We refer to Figure~\ref{fig:proof-schematic} for an illustration of the technical steps mentioned above.
\bigskip

An important observation is that with large probability, we have $Y \not \in \Ball^\dimension_p$ and hence $\Pi_{\Ball^\dimension_p}(Y)$ lies on the boundary of the ball $\Ball^\dimension_p$. Nonetheless, in order to fit the noise coordinates, the MLE returns a vector which is denser than the true $\thetastar$, thereby inducing suboptimal error on the sparse signal components. We now turn to Proposition~\ref{prop:suboptimality-of-erm-hard}.

\paragraph{Proof sketch for Proposition~\ref{prop:suboptimality-of-erm-hard}:} We follow the same general strategy as in Proposition~\ref{prop:suboptimality-of-erm-easy}, but now we must take $k > 1$ in order to balance the signal sparsity with the noise magnitude, which are the two crucial components of this argument. In this regime, it becomes necessary to select the parameter $k$ as a function of $(\sigma, p, d)$, specifically by 
\[
k = k(\sigma, p, d) \asymp \twomax{\Big(\frac{1}{\sigma^q d}\Big)^{\tfrac{p-1}{2-p}}}{1}.
\]
Recalling that $k$ is the sparsity of $\thetastar$, note that as $\sigma \to d^{-1/q}$ this formula shows we are selecting sparser $\thetastar$, since $k \to \Theta(1)$. On the other hand, as $\sigma \to d^{-1/p}$, we are selecting denser $\thetastar$, since $k \to \Theta(d)$. 
As a result, the hardest instance for MLE lies on the boundary of $\Ball^\dimension_p$ but is denser as the noise level shrinks. Having selected this choice of $k$, the argument now proceeds as follows. Though the noise level can be smaller in this setting, we still have $\|Y\|_p > 1$ with probability $\Omega(1)$, so the estimator lies on the boundary of the constraint set. The core idea is to show---via a pigeonhole-type argument---that there exists a sufficiently large subset of the ``noise-only'' coordinates in which the noise magnitude forces $\lambda^\star \gtrsim \sigma d^{1/q}$. Then, using a fixed-point argument and leveraging the identity \eqref{eq:lagrange-fixedpoint}, we conclude that the total estimation error in the ``signal-plus-noise'' coordinates is of order $\sigma d^{1/q}$ for small $\sigma$ and of constant order when $\sigma$ is larger.

\section{Discussion}
\label{sec:discussion}

In this section, we discuss some extensions of our main results. 
In Section~\ref{sec:gen-radii} we present the consequences of Theorem~\ref{thm:cases-where-mle-is-optimal} and Theorem~\ref{thm:cases-where-mle-is-suboptimal} when the parameter $\thetastar$ satisfies $\|\thetastar\|_p \leq \radius$ for general $\radius > 0$. In Section~\ref{sec:n-observations} we record the consequences of our main results when, in place of the single-shot observational model~\eqref{eqn:lp-model}, we observe $n \geq 1$ samples. Finally, we discuss connections to prior work in Section~\ref{sec:prior-work}.
We briefly comment on the computational aspects of the MLE over $\Ball^\dimension_p$ in Section~\ref{sec:computing-the-MLE}.

\subsection{Consequences for general radii}
\label{sec:gen-radii}

In this section, we describe the consequences of our results when the parameter space is taken to be a ball of radius $\radius > 0$ in the 
$p$ norm; that is, the case when $\Theta = \radius \Ball^\dimension_p$. Note that we only need to consider the case $p > 0$, as $p = 0$ with all sparsity and noise levels was already presented in Theorem~\ref{thm:cases-where-mle-is-optimal}.

\paragraph{The minimax rate for the scaled $\ell_p$ balls:}
By rescaling arguments\footnote{Specifically, for any parameter space $\Theta \subset \R^\dimension$ the relation $\MinimaxRisk(\radius \Theta, \sigma) = \radius^2 \MinimaxRisk( \Theta, \sigma/\radius)$ holds for any $\radius > 0$.}, the known minimax rates for the unit $\ell_p$ ball $\Ball^\dimension_p$, as described in displays~\eqref{eqn:minimax-risk-p}, extend to the minimax rate for a general radius. Specifically, we have the following relations:
\begin{subequations}\label{eqn:minimax-risk-p-general-radii}
\begin{enumerate}[label=(\roman*)]
\item for $p \in [2, \infty]$,  
\begin{equation}\label{eqn:minimax-risk-p-geq-2-general-radii}
\MinimaxRisk(\radius \Ball^\dimension_p, \sigma) \asymp \twomin{\sigma^2 d}{\radius^2 d^{1 - 2/p}}, \qquad \mbox{and,}
\end{equation}
\item for $p \in (0, 2)$,  
\begin{equation}\label{eqn:minimax-risk-p-1-2-general-radii}
  \MinimaxRisk(\radius \Ball^\dimension_p, \sigma) 
  \asymp 
  \begin{cases}
  \radius^2 & \sigma^2 \geq \tfrac{\radius^2}{1 + \log \dimension} \\ 
  \radius^p \sigma^{2-p} \left(\log\left(\e \dimension (\sigma/\radius)^p\right)\right)^{1-p/2} &  \sigma^2 \in (\tfrac{\radius^2}{d^{2/p}}, \tfrac{\radius^2}{1 + \log \dimension}) \\ 
  \sigma^2 \dimension & \sigma^2 \leq \tfrac{\radius^2}{\dimension^{2/p}}
  \end{cases}.
\end{equation}
\end{enumerate}
\end{subequations}

Using these known minimax rates, we can extend Theorem~\ref{thm:cases-where-mle-is-optimal}, which described the situations when the MLE is rate-optimal for a unit radius, to a general radius. To state our result, we recall the competitive ratio for the risk of the MLE, now taken over the scaled $\ell_p$ ball,
\begin{equation}\label{defn:comp-ratio-scaled-lp}
\MLECompRatio(\radius \Ball^\dimension_p, \sigma) 
\defn 
\frac{\MLERisk(\radius \Ball^\dimension_p, \sigma)}{\MinimaxRisk(\radius \Ball^\dimension_p, \sigma)}.
\end{equation}
\bcor 
\label{cor:general-radii-optimal}
Fix a radius $\radius > 0$. The competitive ratio for the risk of the maximum likelihood estimator (MLE), as defined in display~\eqref{defn:comp-ratio-scaled-lp}, satisfies 
\[
\MLECompRatio(\radius \Ball^\dimension_p, \sigma) \asymp 1, 
\]
provided that the parameters $(\sigma,\radius, p, d)$ satisfy:
\begin{enumerate}[label=(\roman*)]
\item 
$p \geq 2$, for arbitrary $\sigma, \radius > 0, d \geq 1$; or  
\item 
$p \in (0, 1 + \tfrac{1}{1 + \log \dimension}]$, for some $d \geq 1$, and for arbitrary $\sigma, \radius > 0$; or 
\item 
$p \in (1 + \tfrac{1}{1 + \log \dimension}, 2)$, for  $\sigma, \radius > 0$, $d \geq 1$  which satisfy $\sigma \leq \tfrac{\radius}{d^{1/p}}$ or $\sigma \geq \tfrac{\radius}{\sqrt{1+\log d}}$.
\end{enumerate}
Indeed, in any of these cases, the MLE has minimax order-optimal worst-case risk.
\ecor

We also can characterize the settings in which the MLE is suboptimal, following 
from Theorem~\ref{thm:cases-where-mle-is-suboptimal}. Although our results are nonasymptotic (\ie the lower bounds are based on Propositions~\ref{prop:suboptimality-of-erm-easy} and~\ref{prop:suboptimality-of-erm-hard}), for a succinct statement we again consider a growing dimension. Indeed, as $d \to \infty$, the next result demonstrates that for choices of $(\sigma, \radius, p, \dimension)$ which violate the conditions in Corollary~\ref{cor:general-radii-optimal}, the MLE is minimax rate-suboptimal by an arbitrarily large factor.

\bcor  [Situations when the MLE is suboptimal, general radii] 
\label{cor:general-radii-suboptimal}
In the following situations, the competitive ratio for the risk of the MLE, as defined in display~\eqref{defn:comp-ratio-scaled-lp}, satisfies
\[
\lim_{d \to \infty} \MLECompRatio(\radius \Ball^\dimension_p, \sigma) = \infty.
\] 
Hence, in the following situations, the MLE is rate-suboptimal by arbitrarily large factors as the dimension grows: 
\begin{enumerate}[label=(\roman*)]
    \item if $p \in (1, 2)$, independent of the dimension $d$,
    for any sequence of noise levels and radii $\sigma, \radius > 0$ which satisfy $\sigma d^{1/p} \gg \radius$ and $\sigma \sqrt{\log d} \ll \radius$; and, 
    \item
  if $p = p(d)$ and the conditions~\eqref{eqn:asymptotic-conditions-for-p} hold,   
    for any sequence of noise levels such that 
    $\sigma d^{1/p} \gg \radius d^{2/p-1} (1 - \tfrac{1}{p})^{1/2}$ and $\sigma \sqrt{\log \dimension} \ll \radius$. 
\end{enumerate}
\ecor 

\subsection{Consequences for mean estimation with \texorpdfstring{$n$}{n} Gaussian observations}
\label{sec:n-observations}

In this section, we record the consequences of our main results, Theorems~\ref{thm:cases-where-mle-is-optimal} and~\ref{thm:cases-where-mle-is-suboptimal}, in a setting when 
we have $n$ observations. 

\subsubsection{Setup for \texorpdfstring{$n$}{n} observations}
\label{sec:setup-n-obs}
Specifically, and analogously to the single-observation model~\eqref{eqn:lp-model}, suppose that we have 
$n$ Gaussian observations,
\begin{equation}\label{eqn:n-observation-model}
Y_i = \mustar + \xi_i 
\quad \mbox{where} \quad 
\|\mustar\|_p \leq \radius, 
\quad \mbox{and} \quad 
\xi_i \simiid \Normal{0}{\tau^2 I_\dimension}, 
\quad \mbox{for}~i=1,\ldots, \numobs.
\end{equation}
Equivalently, $Y_i \simiid \Normal{\mustar}{\tau^2 I_\dimension}$. 
Above, the parameter $\tau > 0$ denotes the per-sample noise level, and the radius parameter $\radius > 0$ is the (known) bound on the $\ell_p$ norm of the mean $\mu^\star$. 

\paragraph{Worst-case risk of MLE, minimax risk, and competitive ratio:} 
The maximum likelihood estimator (MLE) for $n$ observations in the model~\eqref{eqn:n-observation-model} takes the form 
\begin{equation}\label{eqn:mle-n-observations}
\hat \mu^{\sf MLE}(Y_1, \dots, Y_\numobs) = 
\argmin_{\mu \in \radius \Ball^\dimension_p} 
\Big\{\, \sum_{i=1}^\numobs \|Y_i - \mu\|_2^2 \, \Big\},
\end{equation}
where we recall that $\radius \Ball^\dimension_p$ denotes the set of vectors $\mu \in \R^\dimension$ with $\|\mu\|_p \leq \radius$.
The worst-case mean squared error of the MLE, under the $n$-observation model~\eqref{eqn:n-observation-model} is then given by 
\[
\MLERisk_\numobs(\radius \Ball^\dimension_p, \tau) \defn \sup_{\mustar \in \radius \Ball^\dimension_p} 
\E_{\{Y_i\}_{i=1}^\numobs \simiid \Normal{\mustar}{\tau^2 I_\dimension}}
\Big[ \, \|\hat \mu^{\sf MLE}(Y_1, \dots, Y_n)- \mustar\|_2^2 \,\Big].
\]
The minimax MSE for this setting is then given by 
\begin{equation}\label{defn:n-sample-minimax-risk}
\MinimaxRisk_\numobs(\radius \Ball^\dimension_p, \tau) 
\defn 
\inf_{\hat \mu} 
\sup_{\mustar \in \radius \Ball^\dimension_p} 
\E_{\{Y_i\}_{i=1}^\numobs \simiid \Normal{\mustar}{\tau^2 I_\dimension}}
\Big[ \, \|\hat \mu(Y_1, \dots, Y_\numobs)- \mustar\|_2^2 \,\Big].
\end{equation}
Finally, the \emph{competitive ratio for the MLE with $n$ observations} is given by 
\begin{equation}\label{def:comp-ratio-n-samples}
\MLECompRatio_\numobs(\radius \Ball^\dimension_p, \tau) 
\defn 
\frac{\MLERisk_\numobs(\radius \Ball^\dimension_p, \tau)}{\MinimaxRisk_\numobs(\radius \Ball^\dimension_p, \tau)}.
\end{equation}
When $\MLECompRatio_\numobs < \infty$, the ratio can be interpreted as the smallest $C \geq 1$ such that the sandwich relation 
$\MinimaxRisk_\numobs \leq \MLERisk_\numobs \leq C \, \MinimaxRisk_\numobs$ holds. That is, the competitive ratio quantifies how much worse the supremum risk of the MLE is, multiplicatively, as compared to a minimax optimal estimator.

\paragraph{Reductions from \texorpdfstring{$n$}{n} samples to one sample:}
Importantly, the $n$-sample model~\eqref{eqn:n-observation-model} can be formally reduced to the single-observation model~\eqref{eqn:lp-model}, both in terms of the minimax risk, and in terms of the worst-case risk of the MLE. Beginning with the minimax risk, note that the sample mean 
\begin{equation*}
\overline{Y}_n \defn 
\frac{1}{\numobs}\sum_{i=1}^\numobs Y_i
\end{equation*}
is a sufficient statistic for 
the $n$ observations~\eqref{eqn:n-observation-model}. Hence, as a consequence of the Rao-Blackwell theorem (\eg Theorem 3.28 in~\cite{Kee10}) and the fact that $\overline{Y}_n \sim \Normal{\mu^\star}{\tfrac{\tau^2}{\numobs} I_\dimension}$, we immediately deduce the following equality of minimax risks,
\begin{equation}\label{eqn:minimax-reduction-to-single-obs}
\MinimaxRisk_\numobs(\radius \Ball^\dimension_p, \tau) = 
\MinimaxRisk\bigg(\radius \Ball^\dimension_p, \sqrt{\frac{\tau^2}{\numobs}}\bigg),
\end{equation}
for every $n \geq 1$ and any $\radius, \tau > 0$.
Turning now to the worst-case risk of the maximum likelihood estimate, note first that for any $\mu \in \R^\dimension$ we have
\[
\frac{1}{\numobs} \sum_{i=1}^\numobs 
\|Y_i - \mu\|_2^2 = 
\|\overline{Y}_n - \mu\|_2^2 
+ \frac{1}{\numobs} \sum_{i=1}^\numobs \|Y_i - \overline{Y}_n\|_2^2. 
\]
This implies that the $n$-sample MLE, as given in display~\eqref{eqn:mle-n-observations}, can be equated to the single-sample MLE, as given in display~\eqref{eqn:single-sample-mle}, but applied to the sufficient statistic $\overline{Y}_n$,
\[
\hat \mu^{\mathsf{MLE}}(Y_1, \dots, Y_n) = \Pi_{\radius \Ball^\dimension_p}(\overline{Y}_n).
\]
As a consequence, we have the following reduction between the worst-case risk of the MLE for $n$ observations and a single observation,
\begin{equation}\label{eqn:mle-risk-reduction-to-single-obs}
\MLERisk_\numobs(\radius \Ball^\dimension_p,\tau) = 
\MLERisk\Big(\radius \Ball^\dimension_p, \sqrt{\frac{\tau^2}{\numobs}}\Big),
\end{equation}
for every $n \geq 1$ and any $\radius, \tau > 0$.
\paragraph{The minimax rate for $n$ observations:} From relation~\eqref{eqn:minimax-reduction-to-single-obs}, we are able to reduce the $n$-sample minimax rate~\eqref{defn:n-sample-minimax-risk} to the single-sample minimax rate given in   display~\eqref{eqn:minimax-risk-p-general-radii}. Consequently, we have,
\begin{subequations}
\begin{enumerate}[label=(\roman*)]
\item for $p \in [2, \infty]$, the $n$-sample minimax rate satisfies
\begin{equation}
\MinimaxRisk_\numobs(\radius \Ball^\dimension_p, \tau) \asymp \twomin{\frac{\tau^2 \dimension}{\numobs}}{\radius^2 d^{1 - 2/p}}, \qquad \mbox{and,}
\end{equation}
\item for $p \in (0, 2)$, the $n$-sample minimax rate satisfies
\begin{equation}
  \MinimaxRisk_\numobs(\radius \Ball^\dimension_p, \tau) 
  \asymp 
  \begin{cases}
  \radius^2 & \numobs \leq \tfrac{\tau^2}{\radius^2} (1+ \log d)  \\ 
  \radius^2 \left(\frac{\tau^2}{\numobs \radius^2} \log\left(\e \dimension (\tfrac{\tau^2}{\numobs\radius^2})^{p/2}\right)\right)^{1-p/2} &  n 
  \in \tfrac{\tau^2}{\radius^2} (1 + \log d, d^{2/p})\\ 
  \tfrac{\tau^2 \dimension}{\numobs} & n \geq \frac{\tau^2}{\radius^2} d^{2/p} 
  \end{cases}, 
  \quad \mbox{and,}
\end{equation}
\item for $p = 0$ and any $s \in [\dimension]$, the $n$-sample minimax rate satisfies
\[
\MinimaxRisk_\numobs(\Ball^\dimension_0(s), \tau) 
\asymp \frac{\tau^2}{\numobs}
s \log \Big(\frac{\e d}{s}\Big).
\]
\end{enumerate}
\end{subequations}

\subsubsection{Results for \texorpdfstring{$n$}{n} observations}

Theorem~\ref{thm:cases-where-mle-is-optimal} then has the following consequence in this slightly more general setting. 
\bcor [Situations when the MLE has order-optimal worst-case risk, $n$ samples]
The competitive ratio for the risk of the maximum likelihood estimator (MLE), as defined in display~\eqref{def:comp-ratio-n-samples} satisfies
\begin{equation}\label{eqn:compratio-n-limit}
\MLECompRatio_\numobs(\radius \Ball^\dimension_p, \tau) \asymp 1, 
\end{equation}
provided that the parameters $(n, \tau, \radius, p, d)$ satisfy:
\begin{enumerate}[label=(\roman*)]
\item
$p \geq 2$, for arbitrary $(\numobs, \tau, \radius, \dimension)$; or  
\item 
$p \in (0, 1 + \tfrac{1}{1 + \log \dimension}]$, for some $d \geq 1$, and for arbitrary $(\tau, \radius, \dimension)$; or 
\item
$p \in (1 + \tfrac{1}{1 + \log \dimension}, 2)$, for $(\numobs, \tau, \radius, \dimension)$ which satisfy $n \geq \tfrac{\tau^2}{\radius^2} d^{2/p}$ or $n \leq \tfrac{\tau^2}{\radius^2} (1 + \log \dimension)$; or, 
\item 
$p = 0$, for any $\sigma > 0$, replacing $r \Ball^\dimension_p$ with $\Ball^\dimension_0(r)$ for any $r \in [d]$ in display~\eqref{eqn:compratio-n-limit}. 
\end{enumerate}
Indeed, in any of these cases, the MLE has minimax order-optimal worst-case risk.
\ecor

We also can characterize the settings in which the MLE is suboptimal, following 
from Theorem~\ref{thm:cases-where-mle-is-suboptimal}.

\bcor  [Situations when the MLE is suboptimal, $n$ samples] 
\label{cor:n-samples-suboptimal}
In the following situations, the competitive ratio for the risk of the MLE, as defined in display~\eqref{def:comp-ratio-n-samples}, satisfies
\[
\lim_{d \to \infty} \MLECompRatio_\numobs(\radius \Ball^\dimension_p, \tau) = \infty.
\] 
Hence, in the following situations, the MLE is rate-suboptimal by arbitrarily large factors as the dimension grows: 
\begin{enumerate}[label=(\roman*)]
    \item \label{item:n-sample-subopt-part-1} 
    if $p \in (1, 2)$, independent of the dimension $d$,
    for $(n, \tau, \radius, \dimension, p)$ which satisfy $\numobs \gg \tfrac{\tau^2}{\radius^2} \log \dimension$ and $\numobs \ll \tfrac{\tau^2}{\radius^2} d^{2/p}$; or,
    \item
  if $p = p(d)$ and the conditions~\eqref{eqn:asymptotic-conditions-for-p} hold,   for $(n, \tau, \radius, \dimension, p)$ which satisfy $\numobs \gg \tfrac{\tau^2}{\radius^2} \log \dimension$ and 
    $\numobs \ll \tfrac{\tau^2}{\radius^2} d^{2 - 2/p} (1 - 1/p)^{-1}$.
\end{enumerate}
\ecor 

We now give two concrete illustrations of the suboptimality, as it arises in Corollary~\ref{cor:n-samples-suboptimal}. In these examples, we fix a scaling of the parameters $(n, \dimension, \tau, \radius, p)$.  Note that in the limits in Corollary~\ref{cor:n-samples-suboptimal}, 
the dimension and sample size both grow as $d \to \infty$, and hence we may parameterize the suboptimal instances by the sample size (rather than the dimension).

Our first example shows that when $d$ is slightly larger than $\sqrt{n}$, the competitive ratio is large. Indeed, it diverges at a rate which is polylogarithmic in the sample size.

\bex [Logarithmic suboptimality]
Fix $\delta \in (0, 1)$ and set the parameters $(\radius, \tau, p, \dimension)$ as follows
\begin{subequations}
\label{eqn:risk-bounds-n-sample-log}
\begin{equation}
\radius^2 = \tau^2 = 1, \quad p = 1 + \delta, \quad \mbox{and} \quad
\dimension = d(n) = n^{p/2} \log n = n^{(1+\delta)/2} \log n.
\end{equation}
The minimax risk and the supremum risk of the MLE satisfy, respectively, 
\begin{equation}
\MinimaxRisk_\numobs \asymp 
\sqrt{\Big(\frac{\log\log n}{n}\Big)^{1-\delta}}
\quad \mbox{and} \quad 
\MLERisk_\numobs \asymp \frac{(\log \numobs)^{\delta/(1+\delta)}}{n^{(1-\delta)/2}}.
\end{equation}
\end{subequations}
Thus, the factor by which the MLE is suboptimal is given by 
\[
\MLECompRatio_\numobs 
= \frac{\MLERisk_\numobs}{\MinimaxRisk_\numobs}
\asymp 
\frac{(\log n)^{\delta/(1+\delta)}}{(\sqrt{\log \log n})^{1-\delta}} \to \infty, 
\quad \mbox{as}~n \to \infty,
\]
as Corollary~\ref{cor:n-samples-suboptimal}\ref{item:n-sample-subopt-part-1} indicates. 
\eex 

The next example illustrates that the competitive ratio can be polynomial in the sample size $n$, for very high-dimensional problems.

\bex [Polynomial suboptimality]
For a more extreme scenario, consider again $\delta \in (0, 1)$, and set the parameters $(\radius, \tau, p, \dimension)$ as follows 
\begin{subequations}
\label{eqn:risk-bounds-n-sample-poly}
\begin{equation}
\radius^2 = \tau^2 = 1, \quad p = 1 + \delta, \quad \mbox{and} \quad
\dimension = \e^{\sqrt{n}}.
\end{equation}
The minimax risk and the supremum risk of the MLE satisfy, respectively, 
\begin{equation}
\MinimaxRisk_\numobs \asymp 
\frac{1}{n^{(1-\delta)/4}}
\quad \mbox{and} \quad 
\MLERisk_\numobs \asymp 
1.
\end{equation}
\end{subequations}
Thus, the factor by which the MLE is suboptimal is given by 
\[
\MLECompRatio_\numobs 
= \frac{\MLERisk_\numobs}{\MinimaxRisk_\numobs}
\asymp 
n^{(1-\delta)/4} \to \infty, 
\quad \mbox{as}~n \to \infty,
\]
again as Corollary~\ref{cor:n-samples-suboptimal}\ref{item:n-sample-subopt-part-1} predicts.
\eex 

Note that the relations~\eqref{eqn:risk-bounds-n-sample-log} and~\eqref{eqn:risk-bounds-n-sample-poly} were obtained by leveraging the analysis in Propositions~\ref{prop:suboptimality-of-erm-easy} and Propositions~\ref{prop:suboptimality-of-erm-hard} as well as the reductions from $n$ samples to a single sample, as discussed in Section~\ref{sec:setup-n-obs} above. We emphasize here the value of these nonasymptotic bounds; while Corollary~\ref{cor:n-samples-suboptimal} illustrates that the competitive ratio diverges upward as the dimension grows, by using the nonasymptotic bounds we gain further information: the rate at which this divergence occurs.

\subsection{Connections with prior work}
\label{sec:prior-work}

We discuss some connections between our work and the prior work on estimation over $\ell_p$ balls and related work on the maximum likelihood estimator.

\paragraph{Estimation over $\ell_p$ balls:} 
As mentioned in the introduction, the work of Donoho and Johnstone~\cite{DonJoh94} revealed a 
phase transition in which nonlinear estimates are necessary for minimax optimal estimation for $p < 2$; see Section~\ref{sec:DJ-theory} for further details. Their work built on a literature studying the rate-optimality~\cite{Pin80, IbrKha85,DonLiuMac90,DonLiu91} and rate-suboptimality~\cite{SacStr82} of linear procedures in various statistical estimation problems.

The survey paper~\cite{Zha12} also gives an account of some extensions of these results to other sparse sets of parameters. The monograph~\cite{Joh19} additionally gave extensions of the results in~\cite{DonJoh94}, including some nonasymptotic results. They also provide a detailed analysis of hard- and soft-thresholded estimates in the Gaussian sequence model.

The very recent paper~\cite[Section 3.2.6]{PraNey25} studies the suboptimality of MLE, and demonstrates that for $p \in (1, 2)$ which is fixed, and independent of the dimension $\dimension$, that MLE over $\ell_p$ balls is minimax rate-suboptimal for the particular choice of $\sigma = d^{1-1/p}$, and for sufficiently large $\dimension$. Unlike this work, their result does not characterize the full range of behaviors of rate-optimality or rate-suboptimality of the MLE, with regard to the parameters $(\sigma, \radius, p, d)$, nor does it provide nonasymptotic control on the risk of the MLE and/or the competitive ratio between the MLE and the minimax risk. From their work, for instance, it is unclear for which instances $(\thetastar, \sigma, p, d)$ the MLE would exhibit its rate-suboptimality. That work also does not study non-convex case, when $p < 1$, and also does not establish the rate-optimality of the MLE for $p \in [1, 1 + \Theta(\tfrac{1}{\log d})]$. 

In a slightly different  but related direction, the paper~\cite{PatMa24} studied a fixed-design linear regression setup with $\ell_p$ parametric constraints. They demonstrated the suboptimality of the Lasso for certain design matrices via analysis of the regularization path and its dependence on the design matrix. 

\paragraph{Suboptimality of the MLE:}

The literature on the suboptimality of the MLE is vast. Indeed, many authors have shown that the MLE can be minimax rate-suboptimal in a variety of situations; see the overview papers by Le Cam~\cite{Le90} and Stigler~\cite{Sti07} for a historical account as well as numerous examples. 
It appears that the first example of the suboptimality of the MLE was pointed out by Birgé and Massart. They demonstrated~\cite[Section 4.B]{BirMas93} the rate-suboptimality of the maximum likelihood estimate for $\alpha$-Hölder density estimation on the unit interval when $\alpha < 1/2$.
Another rate-suboptimal example for density estimation is presented in Section 6.4 of the book~\cite{DevLug01}.
In a different direction, the  paper~\cite{KurRakGun20} establishes the suboptimality of the MLE for estimating the support function of a convex body from noisy measurements. The paper~\cite{KurGaoGunSen24} establishes the rate-suboptimality of the MLE for convex regression in high dimensions.

In the special setting when the observational model consists of a signal, known to lie in a convex compact set, plus additive Gaussian noise, there are still numerous examples of the rate-suboptimality of the MLE. For instance, the elliptical example in~\cite[Lemma 2]{Zha13} and the set constructed by Chatterjee in~\cite[Proposition 1.5]{Cha14}. These examples provided parameter spaces which shrink as the dimension grows for which the MLE exhibits rate-suboptimality. Additional examples were provided in the paper~\cite{PraNey25}, including for pyramids, ellipsoids, solids of revolution, and more. 

\paragraph{Optimality of the MLE:} 
The maximum likelihood estimator (or equivalently, in Gaussian noise, the least squares estimate (LSE)) is known to be minimax rate-optimal in a variety of situations. The classical optimality theory of the MLE was based on regularity conditions such as differentiability in quadratic mean in the parametric case~\cite{Le70, van98} or, in the nonparametric case, on entropy conditions (\eg the convergence of Dudley's (bracketing) entropy integral near the origin)~\cite{van00, VanWel23}.

We now mention some concrete examples in which the optimality of the MLE was established. For instance, from the upper bounds in~\cite[Theorem 2.3]{Zha02} or~\cite[Theorem 3.2]{Bel18}
and the lower bound in~\cite[Corollary 5]{BelTsy15}, the MLE was shown to be rate-optimal in isotonic regression.  According to~\cite[pp. 3007]{WeiFanWai20}, the MLE is minimax rate-optimal for sufficiently regular elliptical parameter sets. The paper~\cite[Section 3]{Han21} gave numerous additional examples for image and edge estimation, binary classification, multiple isotonic regression, and concave density estimation.

Chatterjee gave a general, sufficient condition for the optimality of the MLE~\cite[Proposition 1.6]{Cha14}; the author, however, gives the caveat that this condition ``may be difficult to verify in examples.'' This was expanded on in the recent work~\cite{PraNey25}. Their Corollary 2.19 demonstrates that the LSE is rate-optimal on a convex body if and only if a localized Gaussian width mapping has a controlled Lipschitz seminorm. This condition also appears quite difficult to verify in particular examples.

\paragraph{Variance and bias of the MLE:} Finally, we mention the paper~\cite{KurPutRak23}, whose results provide additional context and interpretation for our main lower bound results. The following result is essentially implicit in their work; see Appendix~\ref{sec:proof-of-variance-convex-MLE} for a complete and self-contained proof for our formulation. 
\bpr [Implicit in~\cite{KurPutRak23}]
\label{prop:variance-of-mle}
There exists a constant $C \geq 1$ such that for any $\sigma > 0$, and any closed, convex set $\Theta \subset \R^\dimension$, we have 
\[
\E_{Y \sim \Normal{\thetastar}{\sigma^2 I_\dimension}}
\Big[\|\Pi_{\Theta}(Y) - \E[\Pi_{\Theta}(Y)]\|_2^2 \Big]
\, \leq 
C \, \MinimaxRisk(\Theta, \sigma).
\]
That is, the variance of the MLE is minimax rate-optimal.
\epr 
An immediate consequence of Proposition~\ref{prop:variance-of-mle} is that if $\Theta \subset \R^\dimension$ is additionally \emph{centrally symmetric}, \ie if $\Theta = - \Theta$ holds, then the MLE is minimax rate-optimal at the origin. Indeed, if $\Theta$ is centrally symmetric, then as Gaussian noise is sign-symmetric, 
it easily follows that the MLE is unbiased at the origin: for $\thetastar = 0$, we have $\E_{Y \sim \Normal{\thetastar}{\sigma^2 I_\dimension}} \Pi_{\Theta}(Y) = \thetastar$, and hence by Proposition~\ref{prop:variance-of-mle} the risk at $\thetastar = 0$ is minimax rate-optimal. 

In the context of our results, note that our constraint sets $\Ball^\dimension_p$ are centrally symmetric convex bodies, and hence we can only expect to derive lower bounds against the optimality of the MLE far from the origin. Indeed, and unsurprisingly given Proposition~\ref{prop:variance-of-mle}, our lower bound results are attained on the boundary of $\Ball^\dimension_p$; see Propositions~\ref{prop:suboptimality-of-erm-easy} and Proposition~\ref{prop:suboptimality-of-erm-hard}. Finally, by Proposition~\ref{prop:variance-of-mle}, we can conclude that the squared \emph{bias}, that is the quantity $\|\thetastar - \E_{\thetastar} \Pi_{\Ball^\dimension_p}(Y)\|_2^2$, must be larger than the minimax rate for the lower bound instances we derived. 

\subsection{Computation of the MLE}
\label{sec:computing-the-MLE}

In this section, we discuss efficient methods for computing the MLE. We will handle the cases $p = 0$, $p=\infty$, and $p > 0$ separately.

\paragraph{Computation of the MLE when $p = 0$:} 
It is straightforward to compute a projection onto the set of $s$-sparse vectors. Indeed, for any $s \in [d]$ and any $y \in \R^\dimension$, the projection is nothing more than the $s$-largest entries of $y$ by magnitude; i.e., we may take
\[
\Pi_{\Ball^\dimension_0(s)}(y) = \sum_{i \in \hat S} y_i e_i 
\quad \mbox{for any} \quad 
\hat S \in \argmax_{S \subset [d], |S| = s} \sum_{i \in S} y_i^2.
\]
This can clearly be computed in polynomial time: simply sort the entries of $y$ by magnitude, keep the largest $s$ entries, and set the remaining entries to zero; ties may be broken arbitrarily.

\paragraph{Computation of the MLE when $p = \infty$:}
It is just as straightforward to project onto an $\ell_\infty$-ball. For any $y\in\R^\dimension$, the Euclidean projection of $y$ onto $r \Ball^\dimension_\infty$ is obtained by clipping each coordinate to the interval $[-r,r]$. Specifically,
\[
  \Pi_{r \Ball^\dimension_\infty}(y) = \sum_{i=1}^\dimension \sign(y_i) \min\left\{|y_i|,r\right\} e_i.
\]
This can clearly be computed in polynomial time by coordinatewise clipping.

\paragraph{Computation of the MLE when $p \in (0, \infty)$:}
For $p \geq 1$, the 
projection is a convex program, and can be computed efficiently. In particular, see the paper~\cite{WonLanXu23}, which develops a method based on Lagrange duality and bisection that works for $p \geq 1$, in which case it can produce an ($\eps$-approximate) solution in polynomial time. Moreover, for $p \in (0, 1)$, this approach---while not guaranteed to converge---produces solutions with small duality gap.

\section{Proofs}
\label{sec:proofs}

\subsection{Proof of Theorem~\ref{thm:cases-where-mle-is-optimal}}
\label{sec:proof-main-optimal-result}

In order to characterize the risk of the MLE, we derive a general risk bound in terms of the local Gaussian width of the parameter set. Specifically, we will bound the worst-case risk of the MLE for a general parameter set, 
\begin{equation}\label{def:sup-risk-mle}
    \MLERisk(\Theta, \sigma) \defn 
    \sup_{\theta \in \Theta} 
    \E_{Y \sim \Normal{\theta}{\sigma^2 I_\dimension}}
    \Big[ \|\Pi_\Theta(Y) - \theta\|_2^2\Big].
\end{equation}
We need to introduce some additional notation for our results. By $\overline{\Theta}$, we denote the symmetral $\overline{\Theta} \defn (\Theta - \Theta)/2$. Moreover, we recall the $\eps$-local Gaussian width of $\overline{\Theta}$, given by
\begin{equation*}
    w\big(\overline{\Theta} \cap \eps \Ball^\dimension_2\big) =  \E \bigg[\sup_{\Delta \in \overline{\Theta} \cap \eps \Ball^\dimension_2} \langle \Delta, \xi \rangle \bigg].
\end{equation*}
Finally, we require the following variational quantity, 
\begin{equation}\label{defn:var-quantity-upper}
\big[\eps_\star(\Theta, \sigma)\big]^2 \defn \inf_{\eps > 0} \Big\{\,\eps^2 +   \frac{\sigma^2}{\eps^2}  w\big(\eps \Ball^\dimension_2 \cap \overline{\Theta}\big)^2 \,\Big\}.
\end{equation}
We make use of the following bound on the supremum risk of the MLE. 
\bpr [Risk bound for the MLE over convex constraint sets]
\label{prop:convex-upper}
For any closed convex $\Theta \subset \R^\dimension$, the supremum risk of the MLE~\eqref{def:sup-risk-mle} satisfies
\begin{align*}
   \MLERisk(\Theta, \sigma)
    \lesssim \threemin{\sigma^2 d}{ \rad(\overline{\Theta})^2}{\eps_\star(\Theta, \sigma)^2}.
\end{align*}
\epr 

\noindent Such risk bounds are fairly standard in the analysis of least squares estimators~\cite{van00, Cha14, Wai19} over convex sets. For completeness, and to justify our explicit constants, we include a (short) proof of Proposition~\ref{prop:convex-upper} in Appendix~\ref{sec:proof-of-convex-upper}.

The non-convex case---when $p \in [0, 1)$---is not covered by Proposition~\ref{prop:convex-upper}. 
Therefore, we separately establish the following risk bounds in the non-convex case. 
\bpr [Risk bound for the MLE when $p < 1$]
\label{prop:mle-risk-bound-sparse}
For any $d \geq 1$, and any $\sigma > 0$, we have:
\begin{enumerate}[label=(\roman*)]
\item 
for $p = 0$ and every $s \in [d]$, we have 
\[
\MLERisk\Big(\Ball^\dimension_0(s), \sigma\Big)
\lesssim 
\sigma^2 s \log\Big(\frac{\e d}{s}\Big).
\]
\label{item:hard-sparse-case}
\item 
for $p \in (0, 1)$, we have 
\[
\MLERisk\Big(\Ball^\dimension_p, \sigma\Big) 
\lesssim 
\begin{cases}
1 & \sigma \geq \tfrac{1}{\sqrt{1 + \log \dimension}} \\ 
(\sigma^2 \log(\e d \sigma^p))^{1-p/2} &
\sigma \in [\tfrac{1}{d^{1/p}}, \tfrac{1}{\sqrt{1+\log \dimension}}] \\ 
\sigma^2 \dimension &
\sigma \leq \tfrac{1}{d^{1/p}}
\end{cases}. 
\]
\label{item:weakly-sparse-case}
\end{enumerate}
\epr
\noindent We establish part~\ref{item:hard-sparse-case}, which is folklore, in Section~\ref{sec:proof-of-hard-sparse-case}. The proof of the more challenging case of weak sparsity (\ie when $p \in (0, 1)$ as in part~\ref{item:weakly-sparse-case}) is presented in Section~\ref{sec:proof-of-weakly-sparse-case}. We note that the implicit constants, given in the proof, are independent of $(\sigma, p, d)$.

\paragraph{Completing the proof of Theorem~\ref{thm:cases-where-mle-is-optimal}:}

We provide a proof for each of the assertions in the statement of the theorem.

\bigskip

\noindent\emph{Assertion~\ref{case:p-geq-2}.} 
For $p \geq 2$, if $\Theta = \Ball^\dimension_p$, then $\overline{\Theta} = \Ball^\dimension_p$. Hence, by Hölder's inequality, the radius satisfies $\rad(\Ball^\dimension_p) = \dimension^{1/2 - 1/p}$.
Thus, from the minimax rate~\eqref{eqn:minimax-risk-p-geq-2} and Proposition~\ref{prop:convex-upper}, we obtain 
\begin{align*}
    \twomin{\sigma^2 \dimension}{\dimension^{1 - 2/p}} 
    \lesssim  
    \MinimaxRisk(\Ball^\dimension_p, \sigma) 
    \leq 
    \sup_{\thetastar \in \Ball^\dimension_p}
    \E_{Y \sim \Normal{\thetastar}{\sigma^2 I_\dimension}} 
    \left[\, \|\Pi_{\Theta}(Y) - \thetastar\|_2^2 \, \right]  
    \lesssim  
    \twomin{\sigma^2 \dimension}{\dimension^{1 - 2/p}}.
\end{align*}
The claimed relation for the worst-case risk then follows directly from the above display.

\bigskip

\noindent\emph{Assertion~\ref{case:p-near-1}, for $p \geq 1$.} Using the fact that for $q = p/(p-1)$, we observe that if $p \leq 1 + {1}/({1 + \log \dimension})$, then $q \geq 2 + \log \dimension$. Therefore, for any $x \in \R^\dimension$, Hölder's inequality yields
\begin{align*}
    \|x\|_q \leq \textrm{exp}\left({\frac{\log \dimension}{2 + \log \dimension}}\right) \|x\|_\infty \leq e \|x\|_\infty,
\end{align*}
which implies that $\e^{-1} \Ball^\dimension_\infty \subset \Ball^\dimension_q$. Using this, we obtain
\begin{align*}
    \left( \eps \Ball^\dimension_2 \cap \Ball^\dimension_p \right)^\circ 
    = 
    \conv \left(\frac{1}{\eps} \Ball^\dimension_2 \cup \Ball^\dimension_{q} \right)
    \supset 
    \conv \left(\frac{1}{\eps} \Ball^\dimension_2 \cup \frac{1}{e} \Ball^\dimension_\infty \right)
    = 
    \left( \eps \Ball^\dimension_2 \cap \e \Ball^\dimension_1 \right)^\circ.
\end{align*}
This shows that $\eps \Ball^\dimension_2 \cap \Ball^\dimension_p \subset \eps \Ball^\dimension_2 \cap \e \Ball^\dimension_1$, from which we recover
\begin{align*}
    w \left(\eps \Ball^\dimension_2 \cap \Ball^\dimension_p \right) 
    \leq
    w \left(\eps \Ball^\dimension_2 \cap e \Ball^\dimension_1 \right)
    =
    \e \, w \left(\frac{\eps}{\e} \Ball^\dimension_2 \cap \Ball^\dimension_1 \right) 
    \eqcolon 
    \e \, \overline{w}(\eps/\e),
\end{align*}
where we defined $\overline{w}(\eps) \defn w(\eps \Ball^\dimension_2 \cap \Ball^\dimension_1)$. Then, by a simple rescaling argument, Proposition~\ref{prop:convex-upper} yields
\begin{align}
\label{ineq:upper-p-near-1-mid}
\begin{split}
    \sup_{\thetastar \in \Theta} \E_{Y \sim \Normal{\thetastar}{\sigma^2 I_\dimension}} \left[\, \|\Pi_{\Theta}(Y) - \thetastar\|_2^2 \, \right] 
    &\leq 
    4\e^2 \,  \inf_{\eps > 0}  \left\{ \eps^2 + \frac{\sigma^2}{\eps^2} \overline{w}(\eps)^2 \right\} \\
    &\leq 
    (2 \e c)^2 \, \inf_{\eps \in (0, 1]} \left\{\eps^2 + \frac{\sigma^2}{\eps^2} \log(\e \dimension \eps^2) \right\},
\end{split}
\end{align}
where the second inequality follows from the bound
\begin{align*}
    \overline{w}(\eps) \leq c \, \sqrt{\log(\e \dimension \eps^2)}, \quad \mbox{for all}~\eps \in (0, 1],
\end{align*}
for some $c \geq 1$, on the local Gaussian complexity (see \cite[Section 4.1]{GorLitMenPaj07} for details). Now, we set 
\begin{align*}
    {\tilde\eps}^{\,2} = (\sigma^2 \log ( \e \dimension \sigma^p))^{1 - p/2}, 
\quad \mbox{for}~\sigma^2 \in [d^{-2/p}, (1 + \log d)^{-1}].
\end{align*}
By inspection, $\tilde\eps \leq 1$. Thus, using this value in the last infimum from equation~\eqref{ineq:upper-p-near-1-mid}, we obtain
\begin{align*}
    \sup_{\thetastar \in \Theta}\E_{Y \sim \Normal{\thetastar}{\sigma^2 I_\dimension}} \left[\, \|\Pi_{\Theta}(Y) - \thetastar\|_2^2 \, \right] 
    \leq  
    (2 \e c)^2 \left( {\tilde\eps}^{\,2} + \frac{\sigma^2}{{\tilde\eps}^{\,2}} \log(\e \dimension {\tilde\eps}^{\,2}) \right) 
    \leq 
    8 \e^2 c^2 {\tilde\eps}^{\,2} \lesssim \MinimaxRisk(\Ball^\dimension_p, \sigma), 
\end{align*}
where the last inequality follows from the minimax risk relation~\eqref{eqn:minimax-risk-p-1-2}, and the penultimate inequality uses the fact that for $\phi(t) \defn \log(\e t)/t$, we have 
\begin{align*}
    \frac{\sigma^2}{{\tilde\eps}^{\,2}} \log(\eps \dimension {\tilde\eps}^{\,2} ) = 
    \sigma^2 \dimension \, \phi({\tilde\eps}^{\,2} \dimension) \leq \sigma^2 \dimension \leq {\tilde\eps}^{\,2}.
\end{align*}
The final two inequalities follow from the fact that ${\tilde\eps}^{\,2} \dimension \geq 1$, $\phi$ is decreasing on $[1, \infty)$, and 
$\sigma^2 \dimension \leq {\tilde\eps}^{\,2}$. These can be easily verified from the fact that $\sigma \mapsto {\tilde\eps}^{\,2}$ is an increasing 
map for $\sigma^2 \in [d^{-2/p}, (1 + \log d)^{-1}]$. This concludes the proof.

\bigskip

\noindent\emph{Assertion~\ref{case:p-near-1}, for $p < 1$.}
The result follows by comparing Proposition~\ref{prop:mle-risk-bound-sparse}\ref{item:weakly-sparse-case} to 
the minimax rate provided in display~\eqref{eqn:minimax-risk-p-1-2}. 

\bigskip 

\noindent\emph{Assertion~\ref{case:p-1-2-high-or-low-noise}.} The result follows by the same argument as in assertion~\ref{case:p-geq-2}; simply note that in this case, $\rad(\Ball^\dimension_p) = 1$, and use 
the minimax risk relations~\eqref{eqn:minimax-risk-p-1-2}.

\bigskip 

\noindent\emph{Assertion~\ref{case:p-0}, for $p = 0$.}
The result follows by comparing Proposition~\ref{prop:mle-risk-bound-sparse}\ref{item:hard-sparse-case} to 
the minimax rate provided in display~\eqref{eqn:minimax-risk-p-0}.

\subsection{Proof of Proposition~\ref{prop:mle-risk-bound-sparse}\ref{item:weakly-sparse-case} in the weakly-sparse case: \texorpdfstring{$p \in (0, 1)$}{p in (0, 1)}}
\label{sec:proof-of-weakly-sparse-case}

We now turn to the weakly-sparse case when $\thetastar \in \Ball^\dimension_p$, for $p \in (0, 1)$. 
We first introduce some notation. Fix $\eps > 0$. We consider the following localized set and random supremum, 
\begin{subequations}
\begin{align}
\cC_\dimension(p, \eps) &\defn 
\Big\{\, x \in \R^\dimension : 
\|x\|_p^p \leq 2, \|x\|_0 \leq \eps\, \Big\}, \quad \mbox{and,}  \\
\Phi_{p, \eps}(\xi) &\defn 
\sup_{x \in \cC_\dimension(p, \eps)} 
\langle x, \xi\rangle.
\end{align}
\end{subequations}
We also define 
\begin{equation}\label{def:fixed-point-weak-sparse}
\eps_p(\sigma)^2 
\defn 
\inf_{\eps > 0} 
\Big\{\, \eps^2  + \frac{\sigma^2}{\eps^2} \|\Phi_{p, \eps}\|_{L^2}^2 \, \Big\}.
\end{equation}
Our first result bounds the worst-case risk of the MLE, for $p \in (0, 1)$, 
using $\eps_p(\sigma)$.
\ble 
\label{lem:weak-sparse-upper-mle-bound}
Fix $p \in (0, 1)$. 
We have 
\[
\MLERisk(\Ball^\dimension_p, \sigma) 
\leq 4 \, 
\min\{\sigma^2 d, 1, \eps_p(\sigma)^2\}.
\]
\ele 
\noindent See Section~\ref{sec:proof-of-weak-sparse-upper-lem} for the proof. 

In order to control $\eps_p(\sigma)$, we require effective control of the localized Gaussian width of the set $\cC_\dimension(p, \eps)$.
The next result provides a sharp---up to constant pre-factor---bound on the localized Gaussian width.
To state the result we define the \emph{$p$th order convex hull} of a set $T \subset \R^\dimension$ as 
\[
\conv_p(T) \defn 
\bigg\{\, x \in \R^\dimension \mid
x = \sum_{j=1}^m \lambda_j t_j~\mbox{for}~t_j \in T, \lambda_j \geq 0, \sum_{j=1}^m \lambda_j^p = 1, m \in \N \,\bigg\}.
\]
\ble 
\label{lem:bound-on-loc-complexity}
Let $d \geq 1$, $p \in (0, 1)$, and fix 
$\eps$ such that $d^{-(2-p)/(2p)} \leq \eps \leq 1$. Set 
\[
s \equiv s(p, \eps) \defn \ceil{\eps^{-\tfrac{2p}{2 - p}}}.
\]
Then we have the inclusions 
\begin{subequations}
\begin{equation}\label{ineq:inclusions-of-the-sparse-set}
\eps \, \conv_p\Big(\Ball^\dimension_0(s) \cap \Ball^\dimension_2\Big) 
\subset 
\cC_\dimension(p, \eps) 
\subset
3 \eps \, \conv\Big(\Ball^\dimension_0(s) \cap \Ball^\dimension_2\Big),  
\end{equation}
and consequently 
\begin{equation}\label{ineq:random-supremum-bounds}
\eps^2 \E \Big[\max_{|S| \leq s(p, \eps), S\subset [\dimension]} \|\xi_S\|_2^2\Big] 
\leq \|\Phi_{p, \eps}\|_{L^2}^2
\leq 9 \eps^2 \E \Big[\max_{|S| \leq s(p, \eps), S\subset [\dimension]} \|\xi_S\|_2^2\Big].
\end{equation}
\end{subequations}
\ele 
\noindent See Section~\ref{sec:proof-of-loc-width-weak-sparse} for the proof.

We remark that inequality~\eqref{ineq:random-supremum-bounds} also implies that we have sharply characterized the Gaussian width of the set $\cC_d(p, \eps) = 2^{1/p} \Ball^\dimension_p \cap \eps \Ball^\dimension_2$. Indeed, we have 
\[
w(\cC_d(p, \eps)) 
\asymp 
\eps w\big(\Ball^\dimension_0(s(p,\eps)) 
\cap \Ball^\dimension_2\big)
\asymp 
\eps^{-\tfrac{p}{2-p}} 
\sqrt{\log\big(\e d \eps^{\tfrac{2p}{2-p}}\big)}.
\]

\paragraph{Completing the proof of Proposition~\ref{prop:mle-risk-bound-sparse}\ref{item:weakly-sparse-case}:}
We will establish that the following bound holds,
\begin{equation}
\label{ineq:desired-bound-weakly-sparse}
\MLERisk\Big(\Ball^\dimension_p, \sigma\Big) 
\leq 540 \, 
\begin{cases}
1 & \sigma \geq \tfrac{1}{\sqrt{1 + \log \dimension}} \\ 
(\sigma^2 \log(\e d \sigma^p))^{1-p/2} &
\sigma \in [\tfrac{1}{d^{1/p}}, \tfrac{1}{\sqrt{1+\log \dimension}}] \\ 
\sigma^2 \dimension &
\sigma \leq \tfrac{1}{d^{1/p}}
\end{cases}. 
\end{equation}
First, note that from Lemma~\ref{lem:weak-sparse-upper-mle-bound}, we immediately 
obtain inequality~\eqref{ineq:desired-bound-weakly-sparse}---indeed, with an improved constant pre-factor---when $\sigma \geq \tfrac{1}{\sqrt{1 + \log \dimension}}$ or $\sigma \leq \tfrac{1}{d^{1/p}}$. Therefore, we need to establish the desired bound~\eqref{ineq:desired-bound-weakly-sparse} in intermediate regime when $\sigma \in [\tfrac{1}{d^{1/p}}, \tfrac{1}{\sqrt{1+\log \dimension}}]$; we make this assumption throughout the remainder of the argument. 

In particular, from Lemma~\ref{lem:weak-sparse-upper-mle-bound}, to establish~\eqref{ineq:desired-bound-weakly-sparse}, it suffices to show that 
\begin{equation}\label{ineq:upper-bound-on-fixed-point-mid-p}
\eps_p(\sigma)^2 \leq 270 \, 
\Big[\sigma^2 \log(\e d \sigma^p)\Big]^{1-p/2}
\quad 
\mbox{when} \quad 
\sigma \in [\tfrac{1}{d^{1/p}}, \tfrac{1}{\sqrt{1+\log \dimension}}].
\end{equation}
From Lemma~\ref{lem:bound-on-loc-complexity} and~\ref{lem:maximal-ineq-norm} applied with $k = s(p, \eps)$, and the definition of $\eps_p(\sigma)$ as given in display~\eqref{def:fixed-point-weak-sparse}, we have 
\begin{align}
\eps_p(\sigma)^2 
&\leq 
54 \, \inf_{\eps \in [d^{-(1/p - 1/2)}, 1]}
\Big\{\eps^2 + \sigma^2 s(p,\eps) 
\log (\e d/s(p,\eps))\Big\} \nonumber \\
&\leq 
108 \, \inf_{\eps \in [d^{-(1/p - 1/2)}, 1]}
\Big\{\eps^2 + \sigma^2 \eps^{-\tfrac{2p}{2-p}} 
\log (\e d \eps^{\tfrac{2p}{p-2}})\Big\} \nonumber \\
&= 
108 \, 
\inf_{\tau \in [d^{-1/p}, 1]}
\Big\{\, \tau^{2-p} + \sigma^2 \tau^{-p} 
\log(\e d \tau^p) \,\Big\}, 
\label{ineq:rescaled-opt-problem-weak-sparse}
\end{align}
where we used the fact that $s(p, \eps) \leq 2 \eps^{-(2p)/(2-p)}$ and made the change of variables $\eps^2 = \tau^{2-p}$. 
We now make the choice of $\tau = \tilde \tau$ where
\[
\tilde \tau \defn \sigma \sqrt{\log (\e d \sigma^p)}. 
\]
Since $\sigma \in [d^{-1/p}, 1/\sqrt{\log(\e d)}]$, we clearly have
$\tilde \tau \in [d^{-1/p}, 1]$. Therefore, 
from inequality~\eqref{ineq:rescaled-opt-problem-weak-sparse}, we obtain 
\begin{align*}
\eps_p(\sigma)^2 
&\leq 
108\, \Big( [\sigma^2 \log(\e d\sigma^p)]^{1-p/2} 
+ \sigma^{2-p} \log^{-p/2}(\e d \sigma^p) 
[\log (\e d \sigma^p) + \tfrac{p}{2} \log \log (\e d \sigma^p)]\Big) \\
&\leq 
270 \, [\sigma^2 \log (\e d \sigma^p)]^{1-p/2},
\end{align*}
where the final inequality uses $\log(t) \leq t$ and $p\leq 1$. The final display above establishes inequality~\eqref{ineq:upper-bound-on-fixed-point-mid-p} and thereby establishes~\eqref{ineq:desired-bound-weakly-sparse}, completing the proof.

\subsubsection{Proof of Lemma~\ref{lem:bound-on-loc-complexity}}
\label{sec:proof-of-loc-width-weak-sparse}

We split the argument into two parts, first we show that \eqref{ineq:inclusions-of-the-sparse-set} implies inequality~\eqref{ineq:random-supremum-bounds}. Then, we establish the bounds in inequality~ \eqref{ineq:inclusions-of-the-sparse-set}.

\paragraph{Inequality~\eqref{ineq:inclusions-of-the-sparse-set} implies inequality~\eqref{ineq:random-supremum-bounds}:}

For arbitrary $T \subset \R^\dimension$, define the functional 
\[
h(T) \defn \E\Big[ \sup_{t \in T} \, \langle t, \xi \rangle^2\Big].
\]
By its definition and the convexity of the map $t \mapsto \langle t, \xi \rangle^2$, the functional $h$ satisfies the properties: 
\begin{equation}\label{ineq:functional-properties}
\begin{gathered}
h(\alpha T) \stackrel{\rm(i)}{=} \alpha^2 h(T), 
\quad h(T) \stackrel{\rm(ii)}{=} h(\conv(T)), \\
h\big(\cC_d(p,\eps)\big) \stackrel{\rm(iii)}{=} \|\Phi_{p, \eps}\|_{L^2}^2, 
\quad \mbox{and} \quad 
h(T) \leq h(S)~\quad \mbox{for}~T \subset S. 
\end{gathered}
\end{equation}
From inclusions~\eqref{ineq:inclusions-of-the-sparse-set} we obtain the inclusions 
\begin{equation}\label{ineq:inclusions-for-cvx-hull}
\eps \, \conv(\Ball^\dimension_0(s) \cap \Ball^\dimension_2) 
\stackrel{\rm(i)}{\subset} 
\conv(\cC_\dimension(p, \eps)) \stackrel{\rm(ii)}{\subset} 
3 \eps \, \conv(\Ball^\dimension_0(s) \cap \Ball^\dimension_2).
\end{equation}
Therefore, we have from the upper inclusion~\eqref{ineq:inclusions-for-cvx-hull}(ii) that
\begin{subequations}
\begin{equation}\label{ineq:phi-upper}
\|\Phi_{p, \eps}\|_2^2 
\stackrel{{\rm \eqref{ineq:functional-properties}(iii)}}{=} 
h\big(\cC_d(p,\eps)\big)
\stackrel{{\rm \eqref{ineq:functional-properties}(iv)}}\leq 
h\big( 3\eps \, \conv(\Ball^\dimension_0(s) \cap \Ball^\dimension_2)\big)
\stackrel{{\rm \eqref{ineq:functional-properties}(i, ii)}}{=}
9 \eps^2 
h\big(\Ball^\dimension_0(s) \cap \Ball^\dimension_2\big).
\end{equation}
In the reverse direction, 
we have from the lower inclusion~\eqref{ineq:inclusions-for-cvx-hull}(i) that
\begin{equation}\label{ineq:phi-lower}
\|\Phi_{p, \eps}\|_2^2 
\stackrel{{\rm \eqref{ineq:functional-properties}(iii)}}{=} 
h\big(\cC_d(p,\eps)\big)
\stackrel{{\rm \eqref{ineq:functional-properties}(iv)}}\geq 
h\big( \eps \, \conv(\Ball^\dimension_0(s) \cap \Ball^\dimension_2)\big)
\stackrel{{\rm \eqref{ineq:functional-properties}(i, ii)}}{=}
\eps^2 
h\big(\Ball^\dimension_0(s) \cap \Ball^\dimension_2\big).
\end{equation}
Finally, we note that 
\begin{equation}\label{ineq:functional-for-sparse}
h\big(\Ball^\dimension_0(s) \cap \Ball^\dimension_2\big) = 
\E \Big[ \sup_{\|x\|_0 \leq s, \|x\|_2 \leq 1} 
\langle x, \xi \rangle^2 \Big] 
= 
\E \Big[\max_{|S| \leq s(p, \eps), S\subset [\dimension]} \|\xi_S\|_2^2\Big].
\end{equation}
\end{subequations}
Combining inequalities~\eqref{ineq:phi-upper},~\eqref{ineq:phi-lower},~and~\eqref{ineq:functional-for-sparse}, we obtain the claimed inequality~\eqref{ineq:random-supremum-bounds}.

\paragraph{Proof of inequality~\eqref{ineq:inclusions-of-the-sparse-set}, upper inclusion:}

Throughout, we use the notation $s = s(p, \eps)$, and we write $m = \ceil{d/s}$. Given any vector 
$x \in \R^\dimension,$ partition its coordinates into $m$ subsets of coordinates with at most $s$ coordinates each, denoted by $I_1, I_2, \dots, I_m$. The top $s$-entries lie in $I_1$, next $s$ in $I_2$ and so on. 
We may decompose such $x$ as follows,
\begin{equation*}
x = 3\eps \cdot \sum_{j=1}^m \frac{\|x_{I_j}\|_2}{3\eps} \cdot \frac{x_{I_j}}{\|x_{I_j}\|_2}
\end{equation*}
To establish the claim, it suffices to demonstrate that 
\begin{equation}\label{ineq:desired-bound-on-sum-of-weights}
\sum_{j=1}^m 
\|x_{I_j}\|_2 \leq 3\eps 
\quad \mbox{for any}~x \in \cC_d(p, \eps).
\end{equation}
To that end, fix $x \in \cC_d(p, \eps)$. 
For any $j \geq 1$, since the sets $I_{j}$ contain progressively smaller by magnitude coordinates, we have 
\[
\|x_{I_{j+1}}\|_{2}^p 
\leq 
|I_{j+1}|^{p/2} \max_{i \in I_{j + 1}} |x_i|^p
\leq |I_{j+1}|^{p/2} \Big(\frac{1}{|I_j|} \sum_{i \in I_j} |x_i|^p\Big)
\leq 
s^{p/2-1} \|x_{I_j}\|_p^p.
\]
Note that for any $x \in \cC_d(p, \eps)$ we additionally have 
\[
\|x_{I_j}\|_2 \leq \|x\|_2^{1-p}\|x_{I_j}\|_2^p  
\leq \eps^{1-p} \|x_{I_j}\|_2^p. 
\]
Therefore, we have 
\begin{equation*}
\sum_{j} 
\|x_{I_j}\|_2 
\leq 
\eps^{1-p} \sum_j \|x_{I_j}\|_2^p 
\leq 
\eps + \eps^{1-p} \sum_{j=1}^{m-1} \|x_{I_j}\|_2^p
\leq 
\eps + \eps^{1-p} s^{p/2 - 1} 
\|x\|_p^p 
\leq 
\eps + 2 \eps^{1-p} s^{p/2 - 1}
\leq 3 \eps, 
\end{equation*}
where we have used the definition of $s = s(p, \eps)$
in the final display. 
Since $x \in \cC_d(p, \eps)$ was arbitrary, 
we have established inequality~\eqref{ineq:desired-bound-on-sum-of-weights}, as needed. 

\paragraph{Proof of inequality~\eqref{ineq:inclusions-of-the-sparse-set}, lower inclusion:}
Fix any $x \in \eps \conv_p(\Ball^\dimension_0(s) \cap \Ball^\dimension_2)$. We may write for some integer $m \geq 1$, 
\[
x = \sum_{j=1}^m \lambda_j t_j 
\quad \mbox{where} \quad t_j \in \Ball^\dimension_0(s) \cap \Ball^\dimension_2,~\lambda_j \geq 0,~\sum_{j=1}^m \lambda_j^p = 1. 
\]
By the quasi-triangle inequality for $\ell_p$, we have by Hölder's inequality that
\begin{subequations}
\begin{equation}\label{ineq:lower-incl-p}
\|x\|_p^p \leq 
\eps^p \sum_{j=1}^m \lambda_j^p \|t_j\|_p^p 
\leq \eps^p s^{(2-p)/2} 
\leq 2, 
\end{equation}
where we used the fact that $(2-p)/2 \in (0, 1)$ 
so that 
\[
s^{(2-p)/2} \leq \eps^{-p} + 1 \leq 2\eps^{-p}.
\]
Additionally, we have 
\begin{equation}\label{ineq:lower-incl-2}
\|x\|_2 \leq \eps \sum_{j=1}^m \lambda_j \leq \eps \sum_{j=1}^m \lambda_j^p \leq \eps,
\end{equation}
\end{subequations}
since $0 < p < 1$. 
Combining inequalities~\eqref{ineq:lower-incl-p} and~\eqref{ineq:lower-incl-2}, we obtain 
$x \in \cC_\dimension(p, \eps)$, as required.

\subsubsection{Proof of Lemma~\ref{lem:weak-sparse-upper-mle-bound}}
\label{sec:proof-of-weak-sparse-upper-lem}
Let $\hat \theta = \Pi_{\Ball^\dimension_p}(Y)$. Note that we have
\[
\|\hat \theta - Y\|_2^2 \leq \|\thetastar - Y\|_2^2.
\]
Expanding both sides, we obtain
\begin{equation}\label{eqn:basic-ineq}
\|\hat \theta - \thetastar\|_2^2 
\leq 2 \sigma \langle \xi, 
\hat \theta - \thetastar\rangle.
\end{equation}
Fix $\eps > 0$. There are two cases, either $\|\hat \theta - \thetastar\|_2^2 \leq \eps$. Or, if not, then we have, 
\[
\|\hat \theta - \thetastar\|_2^2 
\leq 
4 \frac{\sigma^2}{\eps^2} \, 
\Big\langle \xi, \frac{\eps}{\|\hat \theta - \thetastar\|_2} (\hat \theta - \thetastar)\Big\rangle^2 
\leq 
4 \frac{\sigma^2}{\eps^2}
\sup_{\|u\|_p^p \leq 2, \|u\|_2 \leq \eps} 
\langle \xi, u\rangle^2. 
\]
Therefore, combining both cases, we obtain
\begin{subequations}
\begin{align}
\E \Big[\|\hat \theta - \thetastar\|_2^2 \Big] 
&\leq 
\inf_{\eps > 0} 
\Big\{\eps^2 + 4 \frac{\sigma^2}{\eps^2} 
\E\Big[\sup_{\|x\|_p^p \leq 2, \|x\|_2 \leq \eps} \langle x, \xi\rangle^2\Big]\Big\} \nonumber \\ 
&\leq 
4 \eps_p(\sigma)^2 \label{ineq:bound-eps-p-fixed-point}
\end{align}
where we used the definition of $\eps_p(\sigma)$, as given in display~\eqref{def:fixed-point-weak-sparse}. Note that 
$\Ball^\dimension_p \subset \Ball^\dimension_1 \subset \Ball^\dimension_2$,
for all $p \in (0, 1)$. Hence, 
\begin{equation}\label{ineq:diam-bound-p}
\E \Big[\|\hat \theta - \thetastar\|_2^2 \Big] 
\leq \diam(\Ball^\dimension_p)^2 \leq 4.
\end{equation}
On the other hand, we also have from the basic inequality~\eqref{eqn:basic-ineq} and Cauchy-Schwarz that 
\begin{equation}\label{ineq:noise-bound-p}
\E\Big[\|\hat\theta - \thetastar\|_2^2\Big] 
\leq 4 \sigma^2 \E \|\xi\|_2^2 = 
4 \sigma^2 \dimension. 
\end{equation}
\end{subequations}
Combining inequalities~\eqref{ineq:bound-eps-p-fixed-point},~\eqref{ineq:diam-bound-p}, and~\eqref{ineq:noise-bound-p} yields the claim.

\subsection{Proof of Theorem~\ref{thm:cases-where-mle-is-suboptimal}}
\label{sec:proof-main-suboptimal-result}

We provide a proof for the two assertions in the statement of the theorem. We will make use of the following lemma.

\ble \label{lem:upper-bound-p-in-(0,1)}
For $p \in (1, 2)$, if $\sigma \in [\frac{1}{d^{1/p}},\frac{1}{\sqrt{1 + \log d}}]$, then the risk of the MLE satisfies
\begin{equation}
\MLERisk(\Ball^\dimension_p, \sigma) \lesssim \twomin{1}{\sigma \sqrt{q}\dimension^{1/q}}.
\end{equation}
\ele

\medskip

\noindent\emph{Assertion~\ref{thm:suboptimal-p-independent}.} Fix $p \in (1,2)$, independent of the dimension $d$. From Lemma~\ref{lem:upper-bound-p-in-(0,1)} and Proposition~\ref{prop:suboptimality-of-erm-hard}, we have $\MLERisk(\Ball^\dimension_p, \sigma) \gtrsim \min \{1, \sigma \dimension^{1/q}\}$. 
Additionally, the minimax rates given in display~\eqref{eqn:minimax-risk-p-1-2} gives $\MinimaxRisk(\Ball^\dimension_p, \sigma) \asymp (\sigma^2 \log(\e \dimension \sigma^p))^{1-p/2}$. Therefore, setting $\alpha \defn \tfrac{p -1}{p}, \beta \defn 1 - \tfrac{p}{2}$, we  have
\begin{align} \MLECompRatio(\Ball^\dimension_p, \sigma) 
&\gtrsim 
\twomin{\frac{\sigma \dimension^{1/q}}{(\sigma^2 \log(\e \dimension \sigma^p))^{1-p/2}}}{ 
\frac{1}{(\sigma^2 \log(\e \dimension \sigma^p))^{\beta}}}\nonumber \\ 
&=  
\Big(\twomin{\frac{(d \sigma^p)^{\alpha/\beta}}{\sigma^2 \log(\e \dimension \sigma^p)}}{ 
\frac{1}{\sigma^2 \log(\e \dimension \sigma^p)}}\Big)^\beta. \label{ineq:lower-part-1-comp-ratio-diverge}
\end{align}
Fix a constant $C > 1$. 
Note that for $d$ sufficiently large
$(d \sigma^p)^{\alpha/\beta} \geq C \log (\e d \sigma^p) \geq C \sigma^2 \log (\e d \sigma^p)$, 
where we have used the fact that $\sigma^2 \leq 1$ holds for large $d$ and that $\sigma^p d \gg 1$. Now, note that since $\sigma^2 \leq 1$, for $d$ large, we also have 
$\sigma^2 \log (\e d \sigma^p) \leq \sigma^2 \log(\e d) \leq 1/C$. 
Thus, combining these two cases, for $d$ sufficiently large, the bound in inequality~\eqref{ineq:lower-part-1-comp-ratio-diverge} provides, for sufficiently large $d$, that
\[
\MLECompRatio(\Ball^\dimension_p, \sigma) \gtrsim C^\beta. 
\]
Since $C > 1$ was arbitrary, the above display establishes the assertion in this case.

\bigskip

\noindent\emph{Assertion~\ref{thm:suboptimal-p-dependent}.}
Throughout, we use the shorthand notation $p$ for the norm index, but recall here that $p = p(d)$; that is, the norm index is a function of the dimension. 
By the condition $\limsup_{d \to \infty}
p < 2$ in \eqref{eqn:asymptotic-conditions-for-p}, for sufficiently large $d$, 
we have $1 - \tfrac{p}{2} > \eps > 0$ for some fixed $\eps > 0$.
Moreover, note that the condition $\sigma^p d \gg d^{2-p} (1 - \tfrac{1}{p})^{p/2}$ is equivalent to $\sigma \gg \frac{1}{\sqrt{q}d^{1/q}}$, where we recall that $q = p/(p-1)$; it is the conjugate exponent. In this regime, we know from Lemma~\ref{lem:upper-bound-p-in-(0,1)} and Proposition~\ref{prop:suboptimality-of-erm-easy} that the risk of the MLE is order one, \ie that $\MLERisk(\Ball^\dimension_p, \sigma) \asymp 1$. Recalling the minimax rate from \eqref{eqn:minimax-risk-p-1-2}, the competitive ratio for the risk of the MLE then satisfies
\[
    \MLECompRatio(\Ball^\dimension_p, \sigma) \asymp \frac{1}{(\sigma^2 \log(\e \dimension \sigma^p))^{1-p/2}} \geq \left(\frac{1}{ \sigma^2 \log(\e \dimension \sigma^p)}\right)^\varepsilon \geq \left(\frac{1}{\sigma^2 \log(\e \dimension)}\right)^\varepsilon.
\]
Fix $C > 1$. Since $\sigma^2 \log d \ll 1$, it follows that for $d$ sufficiently large, $\sigma^2 \log(\e d) \leq 1/C$. Hence,
\[
\MLECompRatio(\Ball^\dimension_p, \sigma)
\gtrsim C^\eps.
\]
Since $C > 1$ was arbitrary, the display above establishes the desired claim.

\subsubsection{Proof of Lemma~\ref{lem:upper-bound-p-in-(0,1)}}

From Proposition~\ref{prop:convex-upper}, we know that
\[
\MLERisk(\Ball^\dimension_p, \sigma) 
\lesssim 
\twomin{1}{\eps_\star(\Ball^\dimension_p, \sigma)^2},
\]
where we used the fact that $\rad(\Ball^\dimension_p) = 1$ when $1 \leq p \leq 2$. Since $\Ball^\dimension_p$ is a centrally symmetric set, we see by duality that 
\[
w(\overline{\Ball^\dimension_p} \cap \eps \Ball^\dimension_2) = 
w(\Ball^\dimension_p \cap \eps \Ball^\dimension_2) \leq w(\Ball^\dimension_p) = \E \|\xi\|_q, \quad \mbox{for any}~\eps > 0. 
\]
By definition of the variational quantity $\eps_\star(\Ball^\dimension_p, \sigma)$, see display~\eqref{defn:var-quantity-upper}, we have
\begin{equation}\label{eqn:upper-on-var-quantity}
    \eps_\star(\Ball^\dimension_p, \sigma)^2 
    \leq 
    \inf_{\eps > 0} \left\{\,\eps^2 +   \frac{\sigma^2}{\eps^2}  \E \left[\|\xi\|_q \right]^2 \,\right\}
    = 2 \sigma \E \left[\|\xi\|_q \right].
\end{equation}
Since $q \geq 1$, Jensen's inequality yields
\begin{equation}\label{eqn:upper-on-q-norm}
\E\|\xi\|_q  \leq \E[\|\xi\|_q^q]^{1/q} = d^{1/q} \|\xi_1\|_{L^q} \lesssim \sqrt{q} d^{1/q}.
\end{equation}
Above, we used the standard fact that $\|\Normal{0}{1}\|_{L^q} \lesssim \sqrt{q}$ for $q \geq 1$, where we emphasize the implicit constant is independent of $q$; see, \eg equation (2.15) in the book~\cite{Ver18}.
Combining inequalities~\eqref{eqn:upper-on-var-quantity} and~\eqref{eqn:upper-on-q-norm} yields the claim.

\subsection{Proof of Proposition~\ref{prop:suboptimality-of-erm-easy}}
\label{sec:proof-prop-easy}
We will establish the following lower bound,
\begin{equation}\label{eqn:desired-bound-easy}
\MLERisk(\Ball^\dimension_p, \sigma) 
\geq \frac{1}{16}.
\end{equation}
Throughout the argument, we use the shorthand $\Pi(Y) = \Pi_{\Ball^\dimension_p}(Y)$ and we work on the event
\begin{equation}\label{eqn:good-event-prop-easy}
\cE \defn \Big\{ \|Y\|_q \geq 2 \Big\} \cap \Big\{ \xi_1 \leq 0 \Big\}.
\end{equation}
Recall, by Lemma~\ref{lem:var-characterization}, that we have 
\[
|Y_i| = |\Pi(Y)_i| + \lambda^\star |\Pi(Y)_i|^{p-1} \quad \mbox{for all}~i=1,\ldots, d.
\]
Additionally, since the vectors $\Pi(Y)$ and $Y$ have the same signs coordinatewise, we can rearrange the display above to obtain 
\[
|Y_i - \Pi(Y)_i| = \lambda^\star |\Pi(Y)_i|^{p-1}, 
\quad \mbox{for all}~i=1,\ldots, d.
\]
In particular, we can rearrange the above display to conclude that
\begin{equation}\label{eqn:lambda-equation-prop-easy}
\lambdastar = \frac{\|Y - \Pi(Y)\|_q}{\|\Pi(Y)\|_p^{p/q}}.
\end{equation}
On the event $\cE$ in display~\eqref{eqn:good-event-prop-easy}, we have $\|Y\|_p \geq \|Y\|_q \geq 2 > 1$, since $q \geq p$. 
Therefore, $Y \not \in \Ball^\dimension_p$, and consequently, $\|\Pi(Y)\|_p = 1$.  Therefore, under the assumption that $\|Y\|_q \geq 2$, we have from display~\eqref{eqn:lambda-equation-prop-easy} and the triangle inequality that
\[
\lambdastar = \frac{\|Y - \Pi(Y)\|_q}{\|\Pi(Y)\|^{p/q}_p} = \|Y - \Pi(Y)\|_q 
\geq \|Y\|_q - \|\Pi(Y)\|_q \geq 
\|Y\|_q - 1 \geq 1.
\]
Above, the penultimate inequality follows from $\|\Pi(Y)\|_q \leq \|\Pi(Y)\|_p = 1$, where again we use $q \geq p$.

We now apply Lemma~\ref{lem:var-characterization} again, this time to the first coordinate. Using $\lambda^\star \geq 1$ and the fact that $p - 1 \leq 1$ and $\Pi(Y) \in \Ball^\dimension_p \subset \Ball^\dimension_\infty$, we have 
\begin{equation}\label{eqn:lower-bound-on-Y-prop-easy}
|Y_1| = |\Pi(Y)_1| + \lambdastar |\Pi(Y)_1|^{p-1} \geq 2 |\Pi(Y)_1|.
\end{equation}
Note that on the event $\cE$ given in display~\eqref{eqn:good-event-prop-easy}, we have $\xi_1 \leq 0$. There are two cases. First, if $Y_1 \leq 0$, then since $\Pi(Y)_1 \leq 0$, we clearly have that $(\thetastar_1 - \Pi(Y)_1)^2 \geq 1$. On the other hand, if $Y_1 \in [0, 1]$, then by display~\eqref{eqn:lower-bound-on-Y-prop-easy}, we have 
  $|\Pi(Y)_1| \leq 1/2$ and consequently $(\thetastar_1 - \Pi(Y)_1)^2 \geq 1/4$. 
  Therefore, on $\cE$, we have shown that $(\thetastar_1 - \Pi(Y)_1)^2 \geq 1/4$. We immediately obtain,
  \begin{equation}\label{ineq:lower-risk}
  \E \|\Pi(Y) - \thetastar\|_2^2 \geq \frac{1}{4} \P(\cE).
  \end{equation}
  We now make use of the following lemma to control the probability of event $\cE$.
  
  \ble \label{lem:small-ball}
    Let $D \geq 44$ be an integer. Suppose that 
    $\xi \sim \Normal{0}{I_D}$. Then for $r \in [2, 2 \log D]$, we have 
    \[
    \P \Big\{ \|\xi\|_r \geq \frac{1}{\sqrt{32\e}} \sqrt{r} D^{1/r}\Big\}  \geq \half. 
    \]
    \ele 
    \noindent See Section~\ref{sec:small-ball-proof} for a proof.
  
  Note that by the independence of $\xi_1, \xi_{2:\dimension}$ and the sign-symmetry of $\xi_1$, we have 
  \begin{equation}\label{ineq:event-lower}
  \P(\cE) 
 \geq \half \P\Big\{\|\xi_{2:d}\|_q \geq \frac{2}{\sigma} \Big\} \stackrel{\rm(a)}{\geq} 
  \half \P\Big\{\|\xi_{2:d}\|_q \geq 
  \frac{1}{\sqrt{32\e}} \sqrt{q} (d-1)^{1/q}\Big\} 
  \stackrel{\rm (b)}{\geq} \frac{1}{4}.
  \end{equation}
Above, to establish inequality~(a), we used our lower bound on $\sigma$: 
\[
\frac{2}{\sigma} \leq 
\underbrace{\frac{\sqrt{32\e}}{10}(1 + \tfrac{1}{d-1})^{1/q}}_{\defn (\star)} 
\cdot \frac{1}{\sqrt{32\e}}\sqrt{q} (d-1)^{1/q}
\leq \frac{1}{\sqrt{32\e}}\sqrt{q} (d-1)^{1/q}.
\]
The final inequality used $d \geq 45$ to conclude that $(\star) \leq \tfrac{\sqrt{32(1.025) \e}}{10}\leq 1$. 
Inequality~\eqref{ineq:event-lower}(b) follows from  Lemma~\ref{lem:small-ball}, taking the parameters $D = d-1$ and $r= q$, and noting that for $d \geq 45$ we have $D \geq 44$ and $q \leq 2 + \log d \leq 2 \log(d-1)$.
Combining displays~\eqref{ineq:lower-risk} and~\eqref{ineq:event-lower} yields the desired inequality~\eqref{eqn:desired-bound-easy} and hence the claim.

\subsubsection{Proof of Lemma~\ref{lem:small-ball}}
\label{sec:small-ball-proof}
Let $u(t) = \P\{|\Normal{0}{1}| \geq t\}$. Note that for any $1 \leq s \leq D/2$, we have 
\[
\P\{|\xi_s^\star| \geq t\} = 
\P \Big\{ \# \{ i : |\xi_i| \geq t \} \geq s \Big\}
\geq 
\P\Big\{\mathsf{Bin}(D, u(t)) \geq s \Big\} 
\geq \frac{1}{2}, 
\]
if $t = u^{-1}(s/D) \geq \sqrt{2 \log \tfrac{D}{s}}$. Therefore, we have $\E |\xi_s^\star| \geq \sqrt{\frac{1}{2} \log \frac{D}{s}}$. 
Consequently, for any integer $1 \leq s \leq D/2$, we have
\begin{equation*}
\E\|\xi\|_r \geq 
\frac{1}{s^{1 - 1/r}} \E \sup_{\substack{x : \|x\|_0 = s, \\ \|x\|_\infty = 1}} 
\langle x, \xi \rangle 
= s^{1/r} 
\E \frac{1}{s} \sum_{i=1}^s |\xi_i^\star| 
\geq 
 s^{1/r} 
\E |\xi_s^\star| \geq 
 s^{1/r}\sqrt{\half \log \frac{D}{s}}.
\end{equation*}
Note that $s = \floor{\frac{D}{\e^{r/2}}}$
satisfies $s \in [1, D/2]$ and hence we have 
\[
\E\|\xi\|_r \geq 
\half \sqrt{r} \floor{\frac{D}{\e^{r/2}}}^{1/r} \geq 
\frac{1}{2^{1 + 1/r} \sqrt{\e}} \sqrt{r} D^{1/r} 
\geq 
\frac{1}{\sqrt{8\e}} \sqrt{r} D^{1/r}.
\]
Since $\xi \mapsto \|\xi\|_r$ is $1$-Lipschitz with respect to the Euclidean norm (as $r \geq 2$), we have 
\[
\Med(\|\xi\|_r) \geq \E \|\xi\|_r - \Var(\|\xi\|_r) \geq \E 
\|\xi\|_r - 1 \geq 
\frac{1}{\sqrt{8\e}} \sqrt{r} D^{1/r} - 1 
\geq \frac{1}{\sqrt{32 \e}} \sqrt{r} D^{1/r},
\]
Above, we used the Gaussian Poincaré inequality~\cite[Proposition 4.1.1]{BakGenLed14} to control the variance, and the final inequality arose by observing 
\[
\frac{\sqrt{r} D^{1/r}}{\sqrt{32 \e}} \geq \frac{\sqrt{D}}{4 \sqrt{\e}} \geq 1.
\]
Above, the first inequality follows since  $r \mapsto \sqrt{r} D^{1/r}$ is increasing when $r \in [2, 2 \log D]$, and the second inequality holds as soon as $D \geq 44$. To conclude, we note 
\[
\P \Big\{ \|\xi\|_r \geq \frac{1}{ \sqrt{32\e}} \sqrt{r} D^{1/r}\Big\} 
\geq \P\Big\{ \|\xi\|_r \geq m \Big\} \geq \half,
\]
as needed.

\subsection{Proof of Proposition~\ref{prop:suboptimality-of-erm-hard}}
\label{sec:proof-prop-hard}

In order to establish the lower bound in Proposition~\ref{prop:suboptimality-of-erm-hard}, 
we will show, 
\begin{equation}\label{eqn:desired-bound-prop-hard}
\MLERisk(\Ball^\dimension_p, \sigma) 
\geq \frac{3}{409600} \, \min\{1, \sigma d^{1/q}\},
\end{equation}
which clearly establishes the claim. 
It is convenient to introduce the following notation. First, we may write for $d \geq 4$, 
\[
d = m + n, \quad \mbox{where}~m \in 4 \N \quad \mbox{and} \quad n \in \{0, 1, 2, 3\}.
\]
Let us define the random variable
\begin{equation}\label{eqn:noise-term-definition}
\NoiseTerm{q}{m} 
\defn 
\min_{\substack{I \subset \{m/2 + 1, \dots, m\} \\ 
|I| = m/4}} 
\frac{4}{m} \sum_{i \in I} |\xi_i|^q.
\end{equation}
We also define the event 
\begin{equation}\label{eqn:event-cE}
\cE(\delta) \defn \left\{\NoiseTerm{q}{m} \geq \delta \right\}.
\end{equation}
We now make use of the following lemma.

\ble \label{lem:suboptimality-of-erm-hard}
Fix $p \in (1, 2)$. Fix $\delta \in (0, 1)$. 
Suppose that the triple $(\sigma, d, p)$ satisfies 
\begin{equation}\label{ineq:sigma-d-p-condition-2}
  \sigma \geq \frac{1}{\delta^{1/q}}\frac{16}{d^{1/p}} \quad \mbox{and} \quad d \geq 4.
\end{equation}
Define the parameters
\[
\lamlower \defn \frac{\sigma m^{1/q}}{2} \left( \frac{\delta}{4} \right)^{1/q} 
\quad \mbox{and} \quad 
k = \ceil{ \lamlower^{-p/(2-p)} }
\]
We set $\thetastar = \frac{1}{k^{1/p}} \left(\1_k, 0_{d-k}\right) \in \Ball^\dimension_p$. Then the MLE 
satisfies the risk bound, 
\begin{equation}\label{ineq:risk-bound-lower-suboptimal}
\E \|\Pi_{\Ball^\dimension_p}(Y) - \thetastar\|_2^2 \geq \frac{\P(\cE(\delta)) \delta^{1/q}}{8 \sqrt{32}} \, \min\{1, \sigma d^{1/q}\}.
\end{equation}
\ele 

The following lemma will be used to bound the event $\cE(\delta)$. 

\ble \label{lem:noise-term-lower-bound}
For any $t \in (0, 1)$, $d \geq 4$, and $q \in [2, 2 + \log(d)]$, we have
\[
\P\big\{\NoiseTerm{q}{m} \geq \half \left(\frac{3 t}{10}\right)^q\Big\} \geq \frac{(1-t)^2}{80}.
\]
\ele 

If we take $t = 1/2$ in Lemma~\ref{lem:noise-term-lower-bound} 
we obtain that with $\delta = \tfrac{1}{2} (\tfrac{3}{20})^q$ that 
\[
\P\left(\cE(\delta)\right) \geq \frac{1}{320}.
\]
Combining the above display with inequality~\eqref{ineq:risk-bound-lower-suboptimal}, we find that 
\[
\E \|\Pi_{\Ball^\dimension_p}(Y) - \thetastar\|_2^2 \geq \frac{1}{2^{1/q}}\frac{3}{51200 \sqrt{32}} \, \min\{1, \sigma d^{1/q}\}
\geq \frac{3}{409600} \min\{1, \sigma d^{1/q}\},
\]
proving the desired lower bound for $\sigma \geq \frac{1}{\delta^{1/q}}\frac{16}{d^{1/p}}$, with $\tfrac{16}{\delta^{1/q}} \leq 151$. In what follows, through a risk monotonicity argument, we extend this lower bound up to $\sigma \geq \frac{1}{d^{1/p}}$.

\ble [Risk monotonicity] \label{lem:monotonicity-risk}
Fix $0 < \sigma \leq \nu$, and let $\Theta \subset \R^\dimension$ be closed and convex. Then for any $\thetastar \in \Theta$, we have 
\begin{subequations}
\begin{equation}
\label{eqn:monotone-as}
\|\Pi_\Theta(\thetastar + \sigma \xi) - \thetastar\|_2 \leq 
\|\Pi_\Theta(\thetastar + \nu \xi) - \thetastar\|_2.
\end{equation}
This implies, for any $\thetastar \in \Theta$, that
\begin{equation}
\label{eqn:monotone-pointwise}
\E_{Y \sim \Normal{\thetastar}{\sigma^2 I_\dimension}} 
\Big[\|\Pi_\Theta(Y) - \thetastar\|_2^2\Big]
\leq 
\E_{Y \sim \Normal{\thetastar}{\nu^2 I_\dimension}} 
\Big[\|\Pi_\Theta(Y) - \thetastar\|_2^2\Big].
\end{equation}
Consequently, 
\begin{equation}
\label{eqn:monotone-sup}
\cR^{\sf MLE}(\Theta, \sigma) 
\leq 
\cR^{\sf MLE}(\Theta, \nu).
\end{equation}
\end{subequations}
\ele 

\paragraph{Completing proof of Proposition~\ref{prop:suboptimality-of-erm-hard}:}
Applying \eqref{eqn:monotone-sup} from Lemma~\ref{lem:monotonicity-risk}, with $\Theta = \Ball^\dimension_p$ and any $\sigma \in [\frac{1}{d^{1/p}}, \frac{151}{d^{1/p}}]$, we obtain
\begin{multline*}
\cR^{\sf MLE}(\Ball^\dimension_p, \sigma) 
\geq 
\cR^{\sf MLE}\left(\Ball^\dimension_p, \frac{1}{d^{1/p}}\right) 
\geq 
\MinimaxRisk \left(\Ball^\dimension_p, \frac{1}{d^{1/p}} \right) \\
\stackrel{\rm(a)}{\geq} 
\frac{1}{868} \frac{d^{1/q}}{d^{1/p}}
\stackrel{\rm(b)}{\geq} 
\frac{1}{868} \frac{\sigma d^{1/q}}{151}
\geq 
\frac{1}{131068} \min\{\sigma d^{1/q}, 1\}.
\end{multline*}
Above, inequality~(a) follows from the minimax rate in \eqref{eqn:minimax-risk-p-1-2}, where the constant $\tfrac{1}{868}$ was obtained by tracing through  the proof of Theorem 11.7 in the monograph \cite{Joh19}; a formal statement with proof is provided in Appendix~\ref{sec:minimax-rate-lp}. 
Inequality~(b) follows from the assumption that $\sigma \leq \frac{151}{d^{1/p}}$. The desired inequality~\eqref{eqn:desired-bound-prop-hard}, and hence the claim, now follows by noting that $\tfrac{1}{131068} > \tfrac{3}{409600}$. 

\subsubsection{Proof of Lemma~\ref{lem:suboptimality-of-erm-hard}}
  By assumption $m \geq 4$, and moreover $\tfrac{d}{2} \leq m \leq d$. 
  We assume that $\cE(\delta)$ occurs, and use the notation $\Psi_i \defn |\Pi_{\Ball^\dimension_p}(Y)_i|$ We define the set 
  \[
  \cI(\lamlower) \defn \Big\{ i \in [m] : \Psi_i \leq \lamlower \Psi_i^{p-1} \Big\}.
  \] 
  Evidently, we have that if $i \not \in \cI(\lamlower), i \in [m]$, then $\Psi_i^p > \lamlower^{p/(2-p)}.$ Hence,
  \[
  \Big|[m] \setminus \cI(\lamlower)\Big| \, \lamlower^{p/(2-p)} \leq \|\Psi\|_p^p \leq 1. 
  \]
  Rearranging, we obtain that $|\cI(\lamlower)| \geq 3m/4$, since $\lamlower^{p/(2-p)} \geq 4/m$, indeed:
  \begin{equation}\label{ineq:lamlower-bound}
    \lamlower^{p/(2-p)} = \left(\frac{\sigma^q m}{2^q} \frac{\delta}{4} \right)^{\tfrac{p-1}{2-p}} 
    \geq 
    \left(\frac{\sigma^q d}{2^q} \frac{\delta}{8} \right)^{\tfrac{p-1}{2-p}}
    \geq \frac{8}{d} \geq \frac{4}{m},
  \end{equation}
  where the penultimate inequality follows from~\eqref{ineq:sigma-d-p-condition-2} and $p \geq 1$. 
  Therefore, by the pigeonhole principle, there is a subset $I_\star \subset \{m/2 + 1, \dots, m\}$ of size $m/4$ for which 
  $I_\star \subset \cI(\lamlower)$. 
  Moreover, for an index $i \in I_\star$, we have $\thetastar_i = 0$, since by~\eqref{ineq:lamlower-bound} we have 
  $k = \ceil{\lamlower^{-p/(2-p)}} \leq m/4 < m/2$. 
  Now, for the sake of contradiction, suppose that $\lambdastar < \lamlower$. Then, we have 
  \[
  \sigma |\xi_i| = |Y_i| = \Psi_i + \lambdastar \Psi_i^{p-1} < 2 \lamlower \Psi_i^{p-1}, 
  \quad \mbox{for all}~i \in I_\star. 
  \]
  In particular, we have $\Psi_i^p > (\tfrac{\sigma |\xi_i|}{2 \lamlower})^{q}$, and consequently
  \[
  1 \geq \sum_{i \in I_\star} \Psi_i^p > \left(\frac{\sigma}{2 \lamlower}\right)^{q} \sum_{i \in I_\star} |\xi_i|^q 
  \geq \left(\frac{\sigma}{2 \lamlower}\right)^{q} \frac{m}{4} \;  \NoiseTerm{q}{m}, 
  \]
  where we used the definition of the noise term in~\eqref{eqn:noise-term-definition}. However, this is a contradiction, 
  since on $\cE(\delta)$, we have by definition of $\lamlower$ that 
  \[
    \left(\frac{\sigma}{2 \lamlower}\right)^{q} \frac{m}{4} \;  \NoiseTerm{q}{m} 
    \geq 
    \left(\frac{1}{\lamlower}\right)^{q} \frac{\sigma^q m }{2^q} \frac{\delta}{4} = 1. 
  \]
  Consequently, $\lambdastar \geq \lamlower$. If $\xi_1 < 0, Y_1 \leq 0$, then $\Pi_{\Ball^\dimension_p}(Y)_1 \leq 0$, and 
  we certainly have 
  \[
  (\thetastar_1 - \Pi_{\Ball^\dimension_p}(Y)_1)^2 \geq (\thetastar_1)^2.
  \]
  On the other hand, if $\xi_1 < 0, Y_1 \geq 0$, then note that $\Psi_1 = \Pi_{\Ball^\dimension_p}(Y)_1$, and 
  \[
  \thetastar_1 + \sigma \xi_1 = \Psi_1 + \lambdastar \Psi_1^{p-1} = f(\Psi_1) \quad \mbox{where} \quad 
  f(t) \defn t + \lambdastar t^{p-1}.
  \]
  Note that since $\lambdastar \geq \lamlower$ we have 
  \[
  f(\thetastar_1/2) \geq \thetastar_1 \Big\{\half + \frac{1}{2^{p-1}} \lamlower (\thetastar_1)^{p-2} \Big\} 
  \geq \thetastar_1
  \]
  For the final inequality we used $(\thetastar_1)^{p-2} = k^{(2-p)/p} \geq \lamlower^{-1}$ and $p \leq 2$. 
  Since $f$ is increasing and since $f(\Psi_1) = \thetastar_1 + \sigma \xi_1 < \thetastar_1$, this implies 
  $\Psi_1 \leq \thetastar_1/2$, and hence 
  \[
  \left(\thetastar_1 - \Pi_{\Ball^\dimension_p}(Y)_1\right)^2 \geq \frac{(\thetastar_1)^2}{4}.
  \]
Observe that the first $k$ coordinates of $\Pi_{\Ball^\dimension_p}(Y) - \thetastar$ have the same distribution. Therefore, 
\begin{subequations}  
\label{ineq:lower-bound-risk-hard}
\begin{equation}
\label{ineq:lower-bound-risk-hard-1}
  \E \|\Pi_{\Ball^\dimension_p}(Y) - \thetastar\|_2^2 \geq k \E \left[\left(\thetastar_1 - \Pi_{\Ball^\dimension_p}(Y)_1\right)^2 \right] 
  \geq \frac{k^{1-2/p}}{4} \P\Big\{\cE(\delta) \land \xi_1 \leq 0\Big\} 
  \geq \frac{k^{1-2/p}}{8} \P(\cE(\delta)).
\end{equation}

  The final inequality used the independence of $\cE(\delta)$ and $\xi_1$ and the fact that $\xi_1$ is symmetric. 
  Since $k \leq \lamlower^{-r} + 1$, where $r = p/(2-p)$, we have 
\begin{equation}
\label{ineq:lower-bound-risk-hard-2}
  k^{1 - 2/p} = k^{-1/r} \geq 2^{-1/r} \min\{1, \lamlower\} \geq
  2^{-1/r} \min\Big\{1, \frac{\sigma d^{1/q}}{2} \left( \frac{\delta}{8} \right)^{1/q}\Big\} 
  \geq
  \frac{1}{\sqrt{32}} \delta^{1/q} \min\{1, \sigma d^{1/q}\}.
\end{equation}
\end{subequations}
  The proof is complete by combining the displays~\eqref{ineq:lower-bound-risk-hard}.

\subsubsection{Proof of Lemma~\ref{lem:noise-term-lower-bound}}

  Put $k = m/2$ and let $g \sim \Normal{0}{I_k}$. By $\{g^\star_i\}_{i=1}^k$, we denote the order statistics of $g$: 
  \[
  |g_1^\star| \geq |g_2^\star| \geq \cdots \geq |g_k^\star|.
  \] 
  Clearly, we have the equality in distribution:
  \[
  \NoiseTerm{q}{m} = S_k \defn \frac{2}{k} \sum_{i= k/2 + 1}^{k} |g_i^\star|^q.
  \]
  Note that if $|g_{3k/4}^\star| \geq t$, then $S_k \geq \frac{t^q}{2}$. Therefore, 
  \begin{equation}\label{eqn:lower-bound-on-noise-term-1}
  \P\{\NoiseTerm{q}{m} \geq t^q/2\} \geq \P\{|g_{3k/4}^\star| \geq t\}.
  \end{equation}
  From Example 10 in the paper~\cite{GorLitSchWer06}, we have 
  \begin{equation}\label{ineq:bounds-on-g-3k-4-star}
    \sqrt{\frac{\pi}{32}} \leq \frac{k/4 + 1}{k + 1}\leq \E |g_{3k/4}^\star| \leq \sqrt{2 \pi} \frac{k/4 + 1}{k + 1} \leq \sqrt{2\pi}
  \end{equation}
  On the other hand, we may also write 
  \[
  |g_{3k/4}^\star| = \phi(g) \defn \inf_{\substack{T \subset [k] \\ |T| = k/4 + 1}} \sup_{\substack{x \in \R^k : \|x\|_1 = 1 \\ \supp(x) = T}} \langle x, g \rangle.
  \]
  Note that if $g, g' \in \R^k$, then for any $T \subset [k], |T| = k/4 + 1$ and any $x \in \R^k, \|x\|_1 = 1, \supp(x) = T$, we have 
  \[
  \langle x, g \rangle \leq \langle x, g' \rangle + \|x\|_2 \|g - g'\|_2 \leq 
  \sup_{\substack{x \in \R^k : \|x\|_1 = 1 \\ \supp(x) = T}} \langle x, g' \rangle + \|g - g'\|_2,
  \]
  where we used $\|x\|_2 \leq \|x\|_1 = 1$. Passing to the supremum over $x$ on the lefthand side, we obtain 
  \[
  \phi(g) \leq \sup_{\substack{x \in \R^k : \|x\|_1 = 1 \\ \supp(x) = T}} \langle x, g' \rangle + \|g - g'\|_2.
  \]
  Passing to the infimum over $T$ on the righthand side, we obtain $\phi(g) - \phi(g') \leq \|g - g'\|_2$, and hence $\phi$ is $1$-Lipschitz. 
  Now, the Gaussian Poincaré inequality~\cite[Proposition 4.1.1]{BakGenLed14} yields 
  \[
  \E[|g_{3k/4}^\star|^2] = \E \phi^2(g) = \Var(\phi) + \E[\phi(g)]^2 \leq 1 + 2\pi,
  \]
  where we used the inequalities~\eqref{ineq:bounds-on-g-3k-4-star}. 
  Then for any $\theta \in(0, 1)$, we have by the Paley-Zygmund inequality that 
  \[
  \P\{|g_{3k/4}^\star| \geq \theta \E |g_{3k/4}^\star| \} \geq (1 -\theta)^2 \frac{ \E[|g_{3k/4}^\star|]^2}{ \E[|g_{3k/4}^\star|^2]} \geq (1-\theta)^2 
  \frac{\pi}{32 (1 + 2\pi)} \geq \frac{(1-\theta)^2}{80}
  \]
  In particular, for $t \in (0, 1)$, we have $\tfrac{3t}{10} \leq t \sqrt{\frac{\pi}{32}} \leq t \E |g_{3k/4}^\star|$. Therefore, combining the previous display with 
  inequality~\eqref{eqn:lower-bound-on-noise-term-1}, we obtain 
  \[
  \P\Big\{S_k \geq \frac{1}{2} \left(\frac{3t}{10}\right)^{q}\Big\} \geq \P\{|g_{3k/4}^\star| \geq t \E |g_{3k/4}^\star|\} \geq  
  \frac{(1-t)^2}{80}, 
  \quad \mbox{for all}~t \in (0, 1).\qedhere 
  \]

\subsubsection{Proof of Lemma~\ref{lem:monotonicity-risk}}

Clearly, claims~\eqref{eqn:monotone-pointwise} and~\eqref{eqn:monotone-sup} follow from inequality~\eqref{eqn:monotone-as}; we thus establish~\eqref{eqn:monotone-as}.
For $\lambda > 0$, put ${p_\lambda \defn \Pi_\Theta(\thetastar + \lambda \xi)}$.
By the variational inequality for projections onto closed, convex sets 
\begin{align}
\|p_\nu - \thetastar\|_2^2 - 
\|p_\sigma - \thetastar\|_2^2 &= 
\|p_\nu - p_\sigma\|_2^2 - 
2 \langle p_\nu - p_\sigma,  (\thetastar + \sigma \xi) - p_\sigma\rangle 
+ 2\sigma \langle p_\nu - p_\sigma, \xi\rangle \nonumber \\
&\geq 
2\sigma \langle p_\nu - p_\sigma, \xi\rangle. \label{ineq:lower-bound-on-diff}
\end{align}
Additionally, projections onto closed convex sets are firmly nonexpansive~\cite[Proposition 4.16]{BauCom17}, yielding
\begin{equation}\label{ineq:gap-lower}
\langle p_\nu - p_\sigma, \xi \rangle 
= 
\frac{\langle p_\nu - p_\sigma, (\thetastar + \nu \xi) - (\thetastar +\sigma \xi)\rangle}{\nu - \sigma}
\geq 
\frac{\|p_\nu - p_\sigma\|_2^2}{\nu - \sigma} \geq 0. 
\end{equation}
Combining inequalities~\eqref{ineq:lower-bound-on-diff} and~\eqref{ineq:gap-lower}, we 
obtain ${\|p_\nu - \thetastar\|_2^2 \geq \|p_\sigma - \thetastar\|_2^2}$, as needed.

\subsection{Proof of Lemma~\ref{lem:var-characterization}}
The existence of the MLE follows from the continuity of the mapping $\vartheta \mapsto \|Y - \vartheta\|_2^2$ 
and the compactness of $\Ball^\dimension_p$. Moreover, its uniqueness follows from the strict convexity of $\vartheta \mapsto \|Y - \vartheta\|_2^2$. Now, by the Karush-Kuhn-Tucker (KKT) conditions (\eg see Corollary 28.2.1 and Theorem 28.3 in~\cite{Roc70}), 
we know $\hat \theta \defn \Pi_{\Theta}(Y) \in \Ball^\dimension_p$ is associated with some $\lambdastar \geq 0$ that satisfies%
\footnote{Here, by a slight abuse of notation, we define $\sign(0) = [-1, 1]$.}
\begin{subequations}
\begin{align}
\label{inc:subgradient-kkt}
    Y_i &= \hat \theta_i + \lambdastar |\hat \theta_i|^{p-1} \sign(\hat \theta_i), \quad \mbox{for}~i=1,\ldots, \dimension, \mbox{ and} \\ 
\label{eqn:slackness}
    0 &=\lambdastar(\|\hat \theta\|_p^p - 1).
\end{align}
\end{subequations}
From equation~\eqref{inc:subgradient-kkt}, we see that $\sign(\hat{\theta}_i) = \sign(Y_i)$. Thus, $\hat{\theta}_i$ can be recovered from the solution $\Psi_{\lambdastar}(|Y_i|)$ to $|Y_i| = \psi + \lambdastar \psi^{p-1}$, as
\begin{align}
\label{eqn:hat-theta-i}
    \hat \theta_i = \sign(Y_i)\Psi_{\lambdastar}(|Y_i|), \quad \mbox{for}~i=1,\ldots, \dimension.
\end{align}
Note that $\sign(Y_i)$ is essential to recover the correct sign of $\hat{\theta}_i$ since $\Psi_{\lambdastar}(|Y_i|)$ is always positive. This establishes relation~\eqref{eqn:coordinate-characterization}. 

We now turn to relation~\eqref{eqn:minimal-lambda}. The case $Y = 0$ can be easily verified. First, if $Y = 0$, then equation~\eqref{inc:subgradient-kkt} yields $\hat\theta = 0$, and substituting this into equation~\eqref{eqn:slackness} gives $\lambdastar = 0$. Conversely, if $Y = 0$, then $\Psi_{p, \lambda}(|Y_i|) = 0$ for all $i = 1, \ldots, d$, so the constraint $\sum_{i=1}^d \Psi_{p, \lambda}(|Y_i|)^p \leq 1$ is trivially satisfied for any $\lambda \geq 0$. Therefore, the minimal value in \eqref{eqn:minimal-lambda} is achieved at $\lambdastar = 0$. 

We now focus on the case $Y \neq 0$. If $Y \in \Ball^\dimension_p$, then $\hat{\theta} = \Pi_{\Theta}(Y) = Y$, and from equation~\eqref{inc:subgradient-kkt} we have that $\lambdastar = 0$. However, from $Y \in \Ball^\dimension_p$, we also have that $\sum_{i=1}^d \Psi_{p, \lambda}(|Y_i|)^p \leq 1$, for all $\lambda \geq 0$, since $\Psi_{p, \lambda}(|Y_i|) \leq |Y_i|$. Therefore, the minimal value in \eqref{eqn:minimal-lambda} is achieved at $\lambdastar = 0$. Finally, if $Y \notin \Ball^\dimension_p$, then $\|\hat{\theta}\|_p = 1$, and therefore $\lambdastar > 0$ by the complementary slackness condition~\eqref{eqn:slackness}. Moreover, from the display equation~\eqref{eqn:hat-theta-i}, $\lambda^\star$ satisfies $\sum_{i=1}^d \Psi_{p, \lambdastar}(|Y_i|)^p = 1$. Now, since the function $\lambda \mapsto \sum_{i=1}^d \Psi_{p, \lambda}(|Y_i|)^p$ is continuous and strictly decreasing, we also have that the minimal value in \eqref{eqn:minimal-lambda} is achieved at the same $\lambdastar$. 

\subsection*{Acknowledgements}

We thank Matey Neykov for a helpful email exchange and Gil Kur for helpful conversations on his own prior work and for mentioning relevant references in the literature on the optimality and suboptimality of the MLE.

This work was supported in part by the Vannevar Bush Faculty Fellowship program under grant number N00014-21-1-2941. 
Liviu Aolaritei acknowledges support from the Swiss National Science Foundation through the Postdoc.Mobility Fellowship (grant agreement P500PT\_222215).
Reese Pathak was partially supported by a Berkeley ARCS Foundation Fellowship and the NSF under grant DMS-2503579. 
Annie Ulichney acknowledges support from the National Science Foundation Graduate Research Fellowship Program under Grant No. DGE 2146752.

\bibliographystyle{plain}
\bibliography{references}

\clearpage

\appendix 

\section{Deferred proofs}

\subsection{Proof of Proposition~\ref{prop:convex-upper}}
\label{sec:proof-of-convex-upper}
We begin by stating two lemmas---established in Appendix~\ref{sec:proofs-remaining-for-convex-upper}---which we need to prove the result. The first lemma provides an upper bound for the risk of the MLE, which involves the following (random) supremum 
\[
\Phi_\eps(\xi) \defn \sup_{\Delta \in \eps \Ball^\dimension_2 \cap 2\overline{\Theta}} \langle \Delta, \xi \rangle.
\]
Recall that $\xi \defn (Y - \thetastar)/\sigma$.

\ble \label{lem:pointwise-risk-bound}
For any closed convex set $\Theta \subset \R^\dimension$, we have,
\[
    \MLERisk(\Theta, \sigma) \leq \eps^2 + \frac{\sigma^2}{\eps^2} \|\Phi_\eps\|_{L^2}^2,
\]
for every $\eps > 0$.
\ele 

Moreover, the second lemma relates the first and second moments of $\Phi_\eps(\xi)$ . 

\ble [Moment comparison for $\Phi_\eps$]\label{lem:second-to-first}
For any $\Theta \subset \R^\dimension$, we have 
\[
\|\Phi_\eps\|_{L^1} \leq \|\Phi_\eps\|_{L^2} \leq \sqrt{1 + \tfrac{\pi}{2}}\,  \|\Phi_\eps\|_{L^1}.
\]
\ele 

\paragraph{Completing the proof of Proposition~\ref{prop:convex-upper}:}
Armed with Lemmas~\ref{lem:pointwise-risk-bound} and~\ref{lem:second-to-first}, we proceed to prove Proposition~\ref{prop:convex-upper}. Fix $\eps > 0$. Combining these two results with $\eps$ rescaled to $2\eps$,  we obtain
\begin{equation*}
    \MLERisk(\Theta, \sigma)
    \leq
    4\eps^2 +  (1 + \tfrac{\pi}{2}) \frac{\sigma^2}{\eps^2}  w\big(\eps\Ball^\dimension_2 \cap \overline{\Theta}\big)^2 \leq 
    4 \Big[\eps^2 +  \frac{\sigma^2}{\eps^2}  w\big(\eps\Ball^\dimension_2 \cap \overline{\Theta}\big)^2\Big].
\end{equation*}
Passing to the infimum over $\eps  > 0$, 
and comparing to the definition~\eqref{defn:var-quantity-upper}, we find
\begin{subequations}
\label{eqn:convex-upper-bounds-all}
\begin{equation}
\label{eqn:convex-upper-bound-1}
    \MLERisk(\Theta, \sigma)
    \leq 
    4\, \eps_\star(\Theta, \sigma)^2.
\end{equation}
Since projections onto closed convex sets are nonexpansive, we immediately have
\begin{equation}
\label{eqn:convex-upper-bound-2}
    \MLERisk(\Theta, \sigma)
    \leq 
    \E_{Y \sim \Normal{\thetastar}{\sigma^2 I_\dimension}} \Big[\|Y - \thetastar\|_2^2 \Big] 
    = 
    \sigma^2 \dimension.
\end{equation}
Lastly, the fact that $\Pi_{\Theta}(Y) \in \Theta$ implies that $\Pi_{\Theta}(Y) - \thetastar \in 2 \overline{\Theta}$. 
Therefore, 
\begin{equation}
\label{eqn:convex-upper-bound-3}
    \MLERisk(\Theta, \sigma) \leq 4 \, \rad(\overline{\Theta})^2.
\end{equation}
\end{subequations}
Combining the inequalities~\eqref{eqn:convex-upper-bounds-all}, we obtain 
\begin{equation*}
    \MLERisk(\Theta, \sigma)\leq 4 \, \min \Big\{\sigma^2 d ,
    \eps_\star(\Theta, \sigma)^2,  \rad(\overline{\Theta})^2 \Big\},
\end{equation*}
as announced.

\subsection{Remaining proofs for Proposition~\ref{prop:convex-upper}}
\label{sec:proofs-remaining-for-convex-upper}

\subsubsection{Proof of Lemma~\ref{lem:pointwise-risk-bound}}

Throughout this argument, we use the notation $\hat \Delta \defn \Pi_{\Theta}(Y) - \thetastar$. Now, we claim the following inequality holds, 
\begin{equation}
\label{eqn:desired-bound}
    \|\hat \Delta\|_2^2 \leq \twomax{\eps^2}{  \sigma^2 \left(\frac{\Phi_\eps(\xi)}{\eps}\right)^2}
    \quad \mbox{for any}~\eps > 0.
\end{equation}
Assuming the bound~\eqref{eqn:desired-bound} holds for the moment, the claim follows by bounding the maximum of the two terms by their sum, taking the expectation on both sides and then passing to the infimum over $\eps > 0$.

\paragraph{Proof of inequality~\eqref{eqn:desired-bound}:}

To establish bound~\eqref{eqn:desired-bound}, fix an $\eps > 0$.
Note that if $\hat \Delta \in \eps \Ball^\dimension_2$, then there is nothing to prove. When $\hat \Delta \not \in \eps \Ball^\dimension_2$, note that by the variational inequality for projections onto closed convex sets, we have
\begin{equation}
\label{ineq:consequence-of-var-ineq}
    \langle \hat \Delta - (Y - \thetastar), \hat \Delta \rangle = \|\hat \Delta \|_2^2 - \sigma \langle \hat \Delta, \xi \rangle \leq 0.
\end{equation}
Consider $\alpha \defn \eps/\|\hat \Delta\|_2$. Since $\alpha \in (0, 1)$, clearly $\alpha \hat\Delta \in \eps \Ball^\dimension_2 \cap 2\overline{\Theta}$, 
where we have used the fact that the sets $\overline{\Theta}$ and $\Ball^\dimension_2$  are star-shaped around the origin. 
Combined with inequality~\eqref{ineq:consequence-of-var-ineq}, we obtain
\begin{equation}
\label{ineq:upper-bound-on-squared-norm}
    \|\hat\Delta\|_2^2 \leq  \frac{\sigma}{\alpha} \langle \alpha \hat\Delta,\xi \rangle \leq \frac{\sigma}{\alpha} \sup_{\Delta \in \eps \Ball^\dimension_2 \cap \overline{\Theta}} \langle \Delta, \xi \rangle =  \frac{\sigma}{\alpha} \Phi_\eps(\xi).
\end{equation}
Since $\|\hat\Delta\|_2 = \varepsilon/\alpha$, inequality~\eqref{ineq:upper-bound-on-squared-norm} implies
\begin{align*}
    \| \hat\Delta\|_2^2 = 
    2 \|\hat \Delta\|_2^2 - 
    \frac{\eps^2}{\alpha^2} \leq 
    2 \frac{\sigma}{\alpha} \Phi_\eps(\xi) - \frac{\varepsilon^2}{\alpha^2} \leq \sigma^2 \left( \frac{\Phi_\eps(\xi)}{\varepsilon} \right)^2.
\end{align*}
The last inequality follows from the arithmetic-geometric mean inequality. Since 
$\|\hat \Delta\|_2^2 \geq \eps^2$, by assumption, 
in this case we have again established~\eqref{eqn:desired-bound}, completing the proof of the claim.

\subsubsection{Proof of Lemma~\ref{lem:second-to-first}}
The first inequality, $\|\Phi_\eps\|_{L^1} \leq \|\Phi_\eps\|_{L^2}$, is immediate. 
We now prove the second inequality. For ease of notation, we define $T \defn  \eps \Ball^\dimension_2 \cap 2\overline\Theta$. Since $T$ is symmetric, $\Phi_\eps(\xi)$ is equal to the supremum $\sup_{t \in T} |\langle t, \xi \rangle|$. Note that
\begin{align*}
    |\langle t, \xi \rangle | \leq
|\langle t, \xi'\rangle| + \|t\|_2 \|\xi - \xi'\|_2 
\leq \Phi_\eps(\xi') + \rad(T) \|\xi - \xi'\|_2, \quad \mbox{for every}~t \in T.
\end{align*}
Therefore $\Phi_\eps$, when viewed as a function of $\xi$, is $\rad(T)$-Lipschitz with respect to the Euclidean norm. By the Gaussian Poincaré inequality \cite[Proposition 4.1.1]{BakGenLed14}, this implies $\Var(\Phi_\eps) \leq \rad(T)^2$. We conclude that 
\begin{equation}
\label{ineq:bound-l2-using-GPI}
    \|\Phi_\eps\|_{L^2}^2 = \Var(\Phi_\eps) + \|\Phi_\eps\|_{L^1}^2 \leq \rad(T)^2 + \|\Phi_\eps\|_{L^1}^2.
\end{equation}
Moreover, using the fact that $\E[|\langle t, \xi \rangle|] = \sqrt{\tfrac{2}{\pi}} \|t\|_2$, we have 
\begin{equation}
\label{ineq:bound-l1-using-rad}
\|\Phi_\eps\|_{L^1} = \E \left[ \sup_{t \in T}\; |\langle t, \xi \rangle|\right] \geq \sup_{t \in T} \E \left[|\langle t, \xi \rangle|\right] \geq \sqrt{\frac{2}{\pi}} \sup_{t \in T} \|t\|_2 = \sqrt{\frac{2}{\pi}} \rad(T).
\end{equation}
Combining bounds~\eqref{ineq:bound-l2-using-GPI} and~\eqref{ineq:bound-l1-using-rad}, we obtain $\|\Phi_\eps\|_{L^2}^2 \leq (1 + {\pi}/{2}) \|\Phi_\eps\|_{L^1}^2$, thereby completing the proof.

\subsection{Proof of Proposition~\ref{prop:mle-risk-bound-sparse} in the hard-sparse case: \texorpdfstring{$p = 0$}{p = 0}}
\label{sec:proof-of-hard-sparse-case}

We will establish that 
\begin{equation}\label{ineq:desired-bound-hard-sparse}
\MLERisk\Big(\Ball^\dimension_0(s), \sigma\Big) 
\leq 
48 \, \sigma^2 s \log\Big( \frac{\e \dimension}{s}\Big).
\end{equation}
Denote $\hat \theta = \Pi_{\Ball^\dimension_0(s)}(Y)$. We evidently have 
\[
\|\hat \theta - Y\|_2^2 \leq \|\thetastar - Y\|_2^2. 
\]
Recall the notation $\xi \defn (Y - \thetastar)/\sigma$. Expanding the squares, and assuming $\hat \theta \neq \thetastar$, we obtain
\[
\|\hat \theta - \thetastar\|_2^2 
\leq 
2 \sigma \Big\langle \xi, \frac{\hat \theta - \thetastar}{\|\hat \theta- \thetastar \|_2} \Big\rangle \|\hat \theta-\thetastar\|_2
\]
Note $\tfrac{\hat \theta - \thetastar}{\|\hat \theta- \thetastar \|_2}$ has 
no more than $2s$ nonzeros and has unit $\ell_2$ norm. We conclude that 
\[
\|\hat \theta - \thetastar\|_2^2  
\leq 4 \sigma^2 \sup_{\|u\|_0 \leq 2s, \|u\|_2 \leq 1} 
\langle u, \xi \rangle^2 
= 4\sigma^2 \max_{|S| \leq 2s, S\subset [\dimension]} \|\xi_S\|_2^2.
\]
Then taking expectations above, we obtain the desired inequality
\begin{equation}
\label{ineq:risk-bound-pointwise-p-0}
\E \Big[\|\hat \theta - \thetastar\|_2^2\Big]
\leq 4 \sigma^2 
\E \Big[\max_{|S| \leq 2s, S\subset [\dimension]} \|\xi_S\|_2^2\Big]
\stackrel{\rm(i)}{\leq} 24 \sigma^2 \min\{2s, d\} \log \frac{\e d}{\min\{2s, d\}} 
\stackrel{\rm(ii)}{\leq} 48 \sigma^2 s \log \frac{\e d}{s}. 
\end{equation}
Inequality (i) above invoked Lemma~\ref{lem:maximal-ineq-norm}, given below, with $k = \min\{2s, d\}$. Inequality (ii) uses the following observation: if $s \geq d/2$, then 
\[
\min\{2s, d\} \log \frac{\e d}{\min\{2s, d\}} = 
d \leq 2s \leq 2 s \log \frac{\e d}{s}. 
\]
On the other hand, if $s \leq d/2$, then 
\[
\min\{2s, d\} \log \frac{\e d}{\min\{2s, d\}} = 
s \log \frac{\e d}{2s} \leq 2s \log \frac{\e d}{s},
\]
where final inequality uses that 
$\tfrac{\e d}{2s}\leq (\tfrac{\e d}{s})^2$. Passing to the supremum over $\thetastar \in \Ball^\dimension_0(s)$ in display~\eqref{ineq:risk-bound-pointwise-p-0} yields the 
desired inequality~\eqref{ineq:desired-bound-hard-sparse}.

\ble 
\label{lem:maximal-ineq-norm}
For $1 \leq k \leq d < \infty$, we have 
\[
\E \Big[\max_{|S| \leq k, S\subset [\dimension]} \|\xi_S\|_2^2\Big]  
\leq 6 k \log \frac{\e d}{k},
\]
where $\xi \sim \Normal{0}{I_\dimension}$.
\ele 
\begin{proof}
    Fix $\lambda \in (0, 1/2)$. By the concavity of the logarithm, we have
    \begin{align*}
    \E \Big[\max_{|S| \leq k, S\subset [\dimension]} \|\xi_S\|_2^2\Big] 
    &\leq 
    \frac{1}{\lambda}
    \log \E \Big[\max_{|S| = k, S \subset [d]} \exp(\lambda \|\xi_S\|_2^2)\Big]  \\
    &\leq 
    \frac{1}{\lambda} 
    \log \binom{d}{k}
    + 
    \frac{1}{\lambda}
    \log \E[\exp(\lambda \chi^2_{k})]\\
    &\leq 
    \frac{1}{\lambda} k \log \frac{\e d}{k} + \frac{1}{2\lambda} k
     \log \frac{1}{1-2\lambda}, \end{align*}
    where we have used the inequality $\binom{d}{k} \leq (\tfrac{\e d}{k})^k$ and the moment generating function of the $\chi^2_k$ random variate. Setting $\lambda = 1/4$, we obtain the announced inequality
    \[
    \E \Big[\max_{|S| \leq k, S\subset [\dimension]} \|\xi_S\|_2^2\Big]  \leq 
    4 k \log \frac{\e d}{k} 
    + 
    (2 \log 2) k 
    < 6 k \log \frac{\e d}{k}. \qedhere
    \]
\end{proof}

\section{Experiment details}
\label{sec:experiment-details}
We fix $p = 1.5$ and radius $\radius = 1$. For each of the 50 ambient dimensions
\[
\dimension \in \left\{\lfloor10^{2+\tfrac{8}{39}k}\rfloor : k=0,\dots,49\right\}\subset[10^{2},10^{4}],
\]
we generate $100$ independent observations from the Gaussian sequence model 
\[
Y = \thetastar + \sigma \xi, \quad \xi \sim 
\Normal{0}{I_\dimension},
\]
and construct Monte Carlo estimates of the squared loss risk of two estimators in each of two signal-to-noise regimes that replicate the settings of Propositions~\ref{prop:suboptimality-of-erm-easy} and~\ref{prop:suboptimality-of-erm-hard}.

\paragraph{Signal and noise construction:} In the first regime (corresponding to Proposition~\ref{prop:suboptimality-of-erm-easy}), the true signal is the first canonical basis vector: $\thetastar = e_1$ where we take the noise level $\sigma = \dimension^{(1/p) - 1}$. In the second regime (corresponding to Proposition~\ref{prop:suboptimality-of-erm-hard}), we take $\sigma = \dimension^{-1/2}$ and  set $\thetastar = k^{-1/p}( 1_k, 0_{\dimension-k})$ where $k$ is computed as defined in Proposition 2, modified to take $k \in [1, \lfloor d/2 \rfloor]$ in order to avoid the trivial all-zero or full-support signals. 

\paragraph{Maximum-likelihood estimator:} We compute the MLE $\widehat{\theta}^{\text{MLE}}$ by projecting $y \in \R^\dimension$ onto the unit $\ell_p$ ball using a nested set of Brent root-finders that numerically solve the KKT conditions. The inner solver returns for a fixed Lagrange multiplier $\lambda \geq 0$ the magnitude of the $i$-th coordinate of the projection by finding the unique solution $\eta_i(\lambda)$ of the equation $\eta_i + \lambda \eta_i^{p-1} = |y_i|$ where $0 \leq \eta_i \leq |y_i|$. The outer solver adjusts $\lambda$ so that the resulting projection $x(\lambda)$ with coordinates $x_i(\lambda) = \operatorname{sgn}(y_i)\eta_i(\lambda)$ satisfies the feasibility condition $\sum_{i = 1}^\dimension |x_i|^p = 1$. The resulting optimal multiplier $\lambdastar$ gives the projection $x(\lambdastar) = \widehat{\theta}^{\text{MLE}}$.

\paragraph{Soft-thresholding estimator:} Following  \cite[Section 1.3]{PatMa24} and \cite[Theorem 3]{DonJoh94}, we implement the soft-thresholding estimator~\eqref{eqn:DJ-soft-thresholding} with the optimal choice of regularization, $\lambda = \sqrt{2 \sigma^2 \log(\e d \sigma^p)}$. 

\paragraph{Risk evaluation:} For each $\dimension$ in each regime, we compute the empirical MSE over 100 independent draws from the Gaussian sequence model for each estimator $\widehat{\theta}^{\text{MLE}}$ and $\widehat{\theta}^{\text{ST}}$. The resulting Monte Carlo risk curves plotted against the theoretical minimax rate curves from Figures~\ref{fig:regime1} and Figure~\ref{fig:regime2}.

\section{Background on the Gaussian sequence model}
\label{sec:known-results-Gaussian-seq-model}

\subsection{Minimax rate of estimation for \texorpdfstring{$\ell_p$}{lp} balls}
\label{sec:minimax-rate-lp}

The purpose of this section is to establish 
the minimax rates provided in~\eqref{eqn:minimax-risk-p}. We note that this result is implicit in the monograph~\cite{Joh19}, albeit with constants that are difficult to track down explicitly.
The lower bound and upper bound involve the following explicit constants: 
\begin{equation}\label{eqn:constants}
\LowerConst \defn \frac{1}{868} 
\quad \mbox{and} \quad 
\UpperConst \defn 6.
\end{equation}
\paragraph{Control functions:} 
We define the so-called \emph{control functions}, which we will establish as---apart from constant prefactors---the minimax MSE. 
In the case $p \in [1, 2]$, we define 
\[
\mathfrak{m}_{d, p}(\sigma) \defn 
\begin{cases}
1 & \sigma^2 \geq \tfrac{1}{1 + \log \dimension} \\ 
(\sigma^2 \log(\e \dimension \sigma^p))^{1-p/2} &  \sigma^2 \in \big[ \tfrac{1}{d^{2/p}}, \tfrac{1}{1 + \log \dimension} \big] \\ 
\sigma^2 \dimension & \sigma^2 \leq \tfrac{1}{\dimension^{2/p}}
\end{cases}.
\]
For $p \geq 2$, we set $\mathfrak{m}_{d, p}(\sigma) = \twomin{\sigma^2 d}{d^{1-2/p}}$.
\btheo [Minimax rate of estimation for $\ell_p$ balls] 
\label{thm:known-result-GSM}
For any $d \geq 1$, $p \in [1, \infty]$, and $\sigma > 0$, the minimax mean squared error, as defined in display~\eqref{defn:minimax-risk}, satisfies the inequalities
\[
\LowerConst \, \mathfrak{m}_{p,d}(\sigma) 
\leq 
\MinimaxRisk(\Ball^\dimension_p, \sigma) 
\leq 
\UpperConst \, \mathfrak{m}_{p,d}(\sigma).
\]
Above, the constants $\LowerConst, \UpperConst$ may be taken explicitly as in display~\eqref{eqn:constants}.
\etheo 

The remainder of this section is devoted to the proof of Theorem~\ref{thm:known-result-GSM}.

\subsubsection{Proof of Theorem~\ref{thm:known-result-GSM} for \texorpdfstring{$p \geq 2$}{p >= 2}}

\paragraph{Upper bound:}
In the case $p \geq 2$, it is straightforward to see 
\[
\MinimaxRisk(\Ball^\dimension_p, \sigma) \leq \mathfrak{m}_{p, d}(\sigma), \quad \mbox{for every $d \geq 1, \sigma > 0$}.
\]
The bound follows immediately by considering either the estimate $\hat \theta(Y) = 0$, and the fact that $\rad(\Ball^\dimension_p)^2 = d^{1-2/p}$, 
or $\hat \theta(Y) = Y$, and the fact that 
$\E \|\sigma \xi\|_2^2 = \sigma^2 d$, where $\xi \sim \Normal{0}{I_\dimension}$.

\paragraph{Lower bound:}
We follow a standard reduction to hyperrectangles. 
Note that $d^{-1/p} \Ball^\dimension_\infty \subset \Ball^\dimension_p$; for $p = \infty$, note that $d^{-1/p} = 1$. Since $d^{-1/p} \Ball^\dimension_\infty$ is a $d$-fold Cartesian product of the intervals $[-d^{-1/p}, d^{-1/p}]$, and the observation $Y \sim \Normal{\theta}{\sigma^2 I_\dimension}$ has independent coordinates, it follows that 
\begin{multline*}
\MinimaxRisk(\Ball^\dimension_p, \sigma) 
\geq 
\MinimaxRisk(d^{-1/p} \Ball^\dimension_\infty, \sigma) \\= 
d \, \inf_{\hat \mu} \sup_{|\mu| \leq d^{-1/p}} 
\E_{z\sim \Normal{\mu}{\sigma^2}}[(\hat \mu(z) - \mu)^2]
\geq 
\frac{2}{5} 
d \cdot \big(\twomin{d^{-2/p}}{\sigma^2}\big) = 
\frac{2}{5} \mathfrak{m}_{p, d}(\sigma),
\end{multline*}
where the penultimate inequality uses the bound (4.40) from the monograph~\cite{Joh19}.

\subsubsection{Proof of Theorem~\ref{thm:known-result-GSM} for \texorpdfstring{$p < 2$}{p <= 2}}

\paragraph{Upper bound:}
In the case $p < 2$, the following bound holds by the same argument as in the case $p \geq 2$,
\[
\MinimaxRisk(\Ball^\dimension_p, \sigma) \leq 
\twomin{1}{\sigma^2 d} \quad \mbox{for every $d \geq 1, \sigma > 0$}.
\]
We used above $\rad(\Ball^\dimension_p) = 1$, which holds as $1 \leq p < 2$. Thus, to establish the upper bound, it suffices to show that 
\[
\MinimaxRisk(\Ball^\dimension_p, \sigma) 
\leq \UpperConst\, \Big(\sigma^2 \log (\e d \sigma^p)\Big)^{1 - p/2} \quad \mbox{when}~d \geq 1,~~\mbox{and}~~ \sigma \in \Big[ \frac{1}{d^{1/p}}, \frac{1}{\sqrt{1 + \log \dimension}} \Big].
\]
To establish this bound, one can simply apply the result in Theorem~1 of~\cite{PatMa24}, noting that all the inequalities continue to hold for $1 \leq p \leq 2$.

\paragraph{Lower bound:}
For $d = 1$, the lower bound is deduced from the case $p \geq 2$; indeed we have $\mathfrak{m}_{d,p}(\sigma) = \twomin{\sigma^2}{1}$ in this case. 
Set $c \defn \tfrac{0.006}{1 + \log 2}$. We obtain by tracing through the argument underlying inequality (11.46) in the monograph~\cite{Joh19} that
\begin{equation}\label{eqn:general-lower}
\MinimaxRisk(\Ball^\dimension_p, \sigma) 
= 
\sigma^2 \MinimaxRisk(\sigma^{-1} \Ball^\dimension_p, 1) 
\geq 
c\, \Big\{
\max_{1 \leq k \leq \dimension} \twomin{k^{1-2/p}}{ \sigma^2 k \log (\e d/k)}\Big\}.
\end{equation}
Suppose that $\sigma \geq \tfrac{1}{\sqrt{1 + \log (d)}}$.
\begin{subequations}
Then taking $k = 1$ in the general lower bound~\eqref{eqn:general-lower}, we obtain 
\begin{equation}\label{eqn:known-p-geq2-high}
\MinimaxRisk(\Ball^\dimension_p, \sigma) 
\geq 
c \, (\twomin{1}{ \sigma^2 \log(\e d)}) = 
c = 
c \, \mathfrak{m}_{d,p}(\sigma).
\end{equation}
Similarly, for $\sigma \leq \tfrac{1}{d^{1/p}}$, we take $k = d$, 
and we obtain from~\eqref{eqn:general-lower} that 
\begin{equation}\label{eqn:known-p-geq2-low}
\MinimaxRisk(\Ball^\dimension_p, \sigma) 
\geq 
c (\twomin{d^{1-2/p}}{\sigma^2 d}) = 
c\sigma^2 d = 
c \mathfrak{m}_{d, p}(\sigma).
\end{equation}
For $\sigma \in [ \tfrac{1}{d^{1/p}}, \tfrac{1}{\sqrt{1 + \log \dimension}}]$ we note that 
\begin{multline*}
\MinimaxRisk(\Ball^\dimension_p, \sigma) 
\geq 
c
\max_{1 \leq x \leq \dimension} \twomin{\ceil{x}^{1-2/p}}{ \sigma^2 \ceil{x} \log (\e d/\ceil{x})}
\\
\geq 
\frac{c}{4} 
\max_{1 \leq x \leq \dimension} \twomin{x^{1-2/p}}{\sigma^2 x \log (\e d/x)}
= 
\frac{c}{4} 
\max_{\tfrac{1}{d^{1/p}} \leq \tau \leq 1} \twomin{\tau^{2-p}}{\sigma^2 \tau^{-p} \log (\e d\tau^p)},
\end{multline*}
since the map $t \mapsto t \log(\e/t)$ is nondecreasing on $[0, 1]$ and since $(\ceil{x}/x)^{1-2/p} = (x/\ceil{x})^{2/p - 1} \geq \frac{1}{4}$.
We now take 
\[
\tau = \sigma \sqrt{\log (\e \dimension \sigma^p)}.
\]
Note that by assumption that $\sigma \in  [ \tfrac{1}{d^{1/p}}, \tfrac{1}{\sqrt{1 + \log \dimension}}]$, we have 
\[
\frac{1}{d^{1/p}} \leq \sigma \leq \tau \leq \sqrt{\sigma^2 \log(\e d)} \leq 1.
\]
Thus, evaluating the general lower bound
with the choice of $\tau$ given above, we find
\begin{equation}\label{eqn:known-p-geq2-mid}
\MinimaxRisk(\Ball^\dimension_p, \sigma) 
\geq 
\frac{c}{4} 
\, 
\tau^{2-p} = 
\frac{c}{4} \mathfrak{m}_{d, p}(\sigma).
\end{equation}
\end{subequations}
Combining inequalities~\eqref{eqn:known-p-geq2-high},~\eqref{eqn:known-p-geq2-low}, and~\eqref{eqn:known-p-geq2-mid}, we find that, since $\tfrac{c}{4} \geq \tfrac{1}{868} \eqcolon c_\ell$ that
\[
\MinimaxRisk(\Ball^\dimension_p, \sigma) 
\geq \LowerConst \, 
\mathfrak{m}_{d,p}(\sigma), 
\]
holds in the case $p \geq 2$ for all $d \geq 1$ and $\sigma > 0$, as needed. 

\subsection{Proof of Proposition~\ref{prop:variance-of-mle}}
\label{sec:proof-of-variance-convex-MLE}

We essentially follow the arguments of~\cite[Theorem 1]{KurPutRak23} with minor modifications to accommodate our formulation and to obtain explicit constants. 

To simplify notation, let us set for $\xi \sim \Normal{0}{I_\dimension}$, 
\[
\hat \mu(\xi) \defn \Pi_\Theta(\thetastar + \sigma \xi), 
\quad 
h_{\eta}(\xi) \defn 
\|\hat \mu(\xi) - \eta\|_2, 
\quad 
\mbox{and} \quad 
\mu \defn \E_{\xi \sim \Normal{0}{I_\dimension}} \hat \mu(\xi).
\]
Define 
\[
\delta \defn \Big(\E \|\hat \mu(\xi) - \mu\|_2^2 \Big)^{1/2}.
\] 
We will demonstrate that 
\begin{equation}\label{ineq:desired-bound-gil}
\delta^2 \leq 128991 \, \MinimaxRisk(\Theta, \sigma)
\end{equation}
Note that the above inequality is equivalent to our claim.

Set $\GenPackSet{\delta}$ to be a maximal local packing of $\Theta$ at $\mu$; note $\mu \in \Theta$ as $\Theta$ is closed and convex. 
Specifically, for distinct $\eta, \eta' \in \GenPackSet{\delta}$, we have $\|\eta - \eta'\|_2 > \tfrac{\delta}{\sqrt{8}}$ and $\eta, \eta' \in \Theta \cap (\mu + 2\delta \Ball^\dimension_2)$, and moreover 
\[
|\GenPackSet{\delta}| = M\Big(\frac{\delta}{\sqrt{8}}, \Theta \cap (\mu + 2\delta \Ball^\dimension_2), \|\cdot\|_2\Big).
\]
Here, $M(\tau, S, d)$ denotes the packing number: it is the largest cardinality of a subset of $S' \subset S$ such that for distinct pairs $x, y \in S'$, $d(x, y) > \tau$ where $d$ is a metric; for instance here it is $d(x, y) = \|x - y\|_2$.

Suppose that 
$\delta^2 \leq c \sigma^2$ for $c = 4(4 \log(2) + 1)$. In this case, we have by Lemma~\ref{lem:trivial-lower-bound} that with $C_1 = 242$, 
\begin{equation}\label{eqn:case-delta-small}
\delta^2 \leq c  \min\{\sigma^2, \diam(\Theta)^2\} \leq C_1\, \MinimaxRisk(\Theta, \sigma).
\end{equation}
Now suppose that $\delta^2 > c \sigma^2$. 
For the sake of contradiction, suppose that 
$\sigma^2 \log |\GenPackSet{\delta}| < \tfrac{\delta^2}{16}$. 
Note that for any $\eta \in \GenPackSet{\delta}$, the map $h_\eta$ is $\sigma$-Lipschitz. Therefore, by Gaussian Lipschitz concentration~\cite[Theorem 2.26]{Wai19}, we have 
\[
\P\Big\{ \E h_\eta(\xi) - h_\eta(\xi) >   \sigma \sqrt{2 \log(2 |\GenPackSet{\delta}|)} \Big\} \leq \frac{1}{2 |\GenPackSet{\delta}|}
\]
Additionally, by Markov's inequality we have
\[
\P\Big\{ \hat \mu(\xi) \in \Theta \cap (\mu + 2\delta \Ball^\dimension_2)\Big\} \geq 1 - \P\Big\{\|\hat \mu(\xi) - \mu\|_2^2 > 4 \delta^2\Big\} \geq \frac{3}{4}.
\]
Therefore, for some $\etastar \in \GenPackSet{\delta}$, we have 
that $\hat \mu(\xi) \in \etastar + \tfrac{\delta}{\sqrt{8}} \Ball^\dimension_2$ with probability at least $\tfrac{3}{4|\GenPackSet{\delta}|}$.
A union bound then gives us that 
\[
\P\Big\{ \E h_{\etastar}(\xi) -h_{\etastar}(\xi) \leq  \sigma \sqrt{2 \log(2 |\GenPackSet{\delta}|)} 
\land 
h_{\etastar}(\xi) 
\leq \frac{\delta}{\sqrt{8}} \Big\} 
\geq 
\frac{1}{4 |\GenPackSet{\delta}|} > 0.
\]
Hence we have 
\[
\Big(\E h_{\etastar}(\xi) \Big)^2 
\leq \Big( \sigma \sqrt{2 \log(2 |\GenPackSet{\delta}|)} + \frac{\delta}{\sqrt{8}}\Big)^2 
\leq \frac{\delta^2}{4} + 
4 \log(2) \sigma^2 
+ 4 \sigma^2 \log |\GenPackSet{\delta}|.
\] 
Above we used the scalar inequality $(x+y)^2 \leq 2(x^2 + y^2)$. We can now use the display above to derive a contradiction,
\begin{align*}
\delta^2 &= \inf_{\eta \in \R^\dimension} 
\E \big[h_{\eta}(\xi)^2\big]\\ 
& \leq
\E[h_\etastar(\xi)]^2 + \Var(h_\etastar)\\
&\leq 
\frac{\delta^2}{4} + 
\frac{1}{4} c\sigma^2 
+ 
4 \sigma^2 \log |\GenPackSet{\delta}| 
\leq \frac{3\delta^2}{4} < \delta^2.
\end{align*}
Above, we used the Gaussian Poincaré inequality~\cite[Proposition 4.1.1]{BakGenLed14} to conclude 
$\Var(h_\etastar) \leq \sigma^2$. 
Consequently, we have 
\[
\frac{\delta^2}{16 \sigma^2} \leq 
\log |\GenPackSet{\delta}| = 
\log M\Big(\frac{\delta}{\sqrt{8}}, \Theta \cap (\mu + 2\delta \Ball^\dimension_2), \|\cdot\|_2\Big) 
\leq 
\sup_{\mu \in \Theta}
\log M\Big(\frac{\delta}{\sqrt{8}}, \Theta \cap (\mu + 2\delta \Ball^\dimension_2), \|\cdot\|_2\Big).
\]
Hence, by Lemma~\ref{lem:equivalent-versions-of-the-minimax-rate} which we state below, and taking the parameters $c_1 = 16, c_2 = \sqrt{8}, c_3 = 2$, we have 
\begin{equation}\label{eqn:case-delta-large}
\delta^2 \leq C_2 \, \MinimaxRisk(\Theta, \sigma),
\end{equation}
with $C_2 = 128991$.
Combining the cases~\eqref{eqn:case-delta-small} and~\eqref{eqn:case-delta-large} completes the proof and establishes the desired inequality~\eqref{ineq:desired-bound-gil}.

\subsubsection{Remaining lemmas}

\ble
\label{lem:trivial-lower-bound}
For every $\sigma > 0$, and any convex $\Theta \subset \R^\dimension$, we have 
\[
\MinimaxRisk(\Theta, \sigma) \geq \frac{1}{16} \min\{\sigma^2, \diam(\Theta)^2\}.
\]
\ele 
\begin{proof}
First, assume that $\diam(\Theta) < \infty$. Resorting to an approximating sequence if necessary, we assume that $t, s \in \Theta$ satisfy $\|t - s\|_2 = \diam(\Theta)$. We can construct for $\alpha \defn \twomin{\sigma}{\|t-s\|_2}$ two points $t_{\alpha}, s_{\alpha}$ on the line segment connecting $t, s$ such that $\|t_\alpha - s_\alpha\|_2 = \alpha$. As $t_\alpha, s_\alpha \in \Theta$, we may apply Le Cam's method and Pinsker's inequality~\cite[Chapter 15]{Wai19} to these two points, yielding
\[
\MinimaxRisk(\Theta, \sigma) \geq \frac{\|t_\alpha - s_\alpha\|_2^2}{8}
\Big\{1 - \sqrt{\frac{1}{4 \sigma^2} \|t_\alpha - s_\alpha\|_2^2}\Big\} 
\geq \frac{1}{16} \alpha^2 = 
\frac{1}{16} \min\{\sigma^2 ,\diam(\Theta)^2\}.\qedhere
\]
If $\Theta$ is unbounded, we may repeat the argument, except we simply take $\alpha = \sigma$ and any pair $(t, s)$ with distance at least $\sigma$.
\end{proof}

\ble 
\label{lem:equivalent-versions-of-the-minimax-rate}
Let $\sigma > 0$ and fix a convex set $\Theta \subset \R^\dimension$. Suppose $\delta > 0$ is such that for some constants $c_1, c_2, c_3 > 1$ we have 
\[
\frac{\delta^2}{c_1 \sigma^2} 
\leq \sup_{\mu \in \Theta} 
\log M\Big(\frac{\delta}{c_2}, 
\Theta \cap (\mu + c_3 \delta \Ball^\dimension_2), 
\|\cdot\|_2\Big). 
\]
Then there is a constant $c_4 = c_4(c_1, c_2, c_3) > 0$ such that 
\[
\delta^2 \leq c_4 \; \MinimaxRisk(\Theta, \sigma). 
\]
\ele 
\begin{proof}
     Define 
    \[
    h(\delta) \defn \sup_{\mu \in \Theta} \log M\Big(\frac{\delta}{c_2}, \Theta \cap (\mu + c_3 \delta \Ball^\dimension_2), 
\|\cdot\|_2\Big).
    \]
   Since $\Theta$ is convex, it is straightforward $h(\delta) \geq h(\delta')$. for $0 < \delta \leq \delta' < \infty$.
    Therefore, using our assumption, we find 
    \[
    \frac{\delta^2}{c_1 \sigma^2} \leq h(\delta) 
    \leq 
    h\Big(\frac{\delta}{c_3 \sqrt{c_1}}\Big) = 
    \sup_{\mu \in \Theta}
    \log M\Big(\frac{\delta}{c_2 c_3 \sqrt{c_1}}, 
    \Theta \cap (\mu + \tfrac{\delta}{\sqrt{c_1}} \Ball^\dimension_2), 
    \|\cdot\|_2\Big) 
    \eqcolon \log M_\delta
    \]
    We may assume that $M_\delta \geq 2$, otherwise the claim is trivial.
    Suppose first that $M_\delta \geq 5$. Then, Fano's inequality (\eg see \cite[Proposition 15.12]{Wai19}) applied to the packing set implied by $M_\delta$ above, gives us 
    \begin{equation}\label{ineq:lower-bound-one}
    \MinimaxRisk(\Theta, \sigma)
    \geq 
    \frac{\delta^2}{4 c_2^2 c_3^2 c_1 }
    \Big\{1 - \frac{\log 2}{\log 5} - 
    \frac{1}{2} \frac{\delta^2/(c_1 \sigma^2)}{\log M_\delta}\Big\}
    \geq \underbrace{\frac{1}{4 c_2^2 c_3^2 c_1} 
    \Big(\half - \frac{\log 2}{\log 5}\Big)}_{\eqcolon C_1} \, \delta^2
    = C_1 \delta^2. 
    \end{equation}
    Now suppose that instead we have $2 \leq M_\delta < 5$. 
    Then we can apply Le Cam's two point lower bound (\eg see \cite[Equation 15.14]{Wai19}) to 
    any pair of distinct points (say $\eta, \eta'$) in the packing.
    Denoting by $\|P - Q\|_{\rm TV}$ the total variation distance between probability measures $P, Q$ on the same space, this will give us
    \begin{align}
    \MinimaxRisk(\Theta, \sigma)
    &\geq \frac{\delta^2}{8 c_2^2 c_3^2 c_1}
    \Big\{1 - \|\Normal{\eta}{\sigma^2 I_\dimension} - 
    \Normal{\eta'}{\sigma^2 I_\dimension}\|_{\rm TV}
    \Big\} \nonumber \\
    &\geq 
    \frac{\delta^2}{8 c_2^2 c_3^2 c_1}
    \Big\{1 - 
    \sqrt{1 - \e^{- \|\eta - \eta'\|_2^2/(2\sigma^2)}}\Big\} \nonumber \\
    &\geq 
    \frac{\delta^2}{8 c_2^2 c_3^2 c_1}
    \Big\{1 - 
    \sqrt{1 - \frac{1}{( M_\delta)^2}}\Big\} \nonumber \\ 
    &\geq\underbrace{\frac{1}{8 c_2^2 c_3^2 c_1}
    \Big\{1 - 
    \sqrt{\frac{15}{16}}\Big\}}_{\eqcolon C_2} \delta^2
    = C_2 \delta^2.\label{ineq:lower-bound-two}
    \end{align}
    Above, we applied the 
    Bregtanolle-Huber inequality to bound the variational distance above by KL~\cite[Lemma 2.1]{BreHub79}.
    By combining bounds~\eqref{ineq:lower-bound-one} and~\eqref{ineq:lower-bound-two}, we obtain
    \[
    \MinimaxRisk(\Theta, \sigma) \geq \min\{C_1, C_2\} \, \delta^2.
    \]
    Taking $c_4 = 1/\min\{C_1, C_2\}$ completes the proof.
\end{proof}
\end{document}

%% file: src_tex/ams_defs.tex
\theoremstyle{plain}

\newtheorem{theo}{Theorem}

\newtheorem{lem}{Lemma}
\newtheorem{prop}{Proposition}
\newtheorem{cor}{Corollary}

\theoremstyle{definition} 

\newtheorem{nota}{Notation}
\newtheorem{de}{Definition}
\newtheorem{as}{Assumption}
\newtheorem{alg}{Algorithm}
\newtheorem{exa}{Example}
\AtBeginEnvironment{exa}{%
  \pushQED{\qed}%
}
\AtEndEnvironment{exa}{\popQED\endexa}

\newcommand{\btheo}{\begin{theo}}
\newcommand{\bde}{\begin{de}}
\newcommand{\ble}{\begin{lem}}
\newcommand{\bpr}{\begin{prop}}
\newcommand{\bno}{\begin{nota}}
\newcommand{\bex}{\begin{exa}}
\newcommand{\bcor}{\begin{cor}}
\newcommand{\spro}{\begin{proof}}
\newcommand{\bas}{\begin{as}}
\newcommand{\balg}{\begin{alg}}

\newcommand{\etheo}{\end{theo}}
\newcommand{\ede}{\end{de}}
\newcommand{\ele}{\end{lem}}
\newcommand{\epr}{\end{prop}}
\newcommand{\eno}{\end{nota}}
\newcommand{\eex}{\end{exa}}
\newcommand{\ecor}{\end{cor}}
\newcommand{\fpro}{\end{proof}}
\newcommand{\eas}{\end{as}}
\newcommand{\ealg}{\end{alg}}

\theoremstyle{plain}

\newtheorem{theos}{Theorem}
\newtheorem{props}{Proposition}
\newtheorem{lems}{Lemma}
\newtheorem{cors}{Corollary}

\theoremstyle{definition}
\newtheorem{exas}{Example}
\newtheorem{algs}{Algorithm}
\newtheorem{asss}{Assumption}
\newtheorem{defns}{Definition}

\newcommand{\btheos}{\begin{theos}}
\newcommand{\etheos}{\end{theos}}
\newcommand{\bprops}{\begin{props}}
\newcommand{\eprops}{\end{props}}
\newcommand{\bdes}{\begin{defns}}
\newcommand{\edes}{\end{defns}}
\newcommand{\blems}{\begin{lems}}
\newcommand{\elems}{\end{lems}}
\newcommand{\bcors}{\begin{cors}}
\newcommand{\ecors}{\end{cors}}
\newcommand{\bexs}{\begin{exas}}
\newcommand{\eexs}{\end{exas}}
\newcommand{\balgs}{\begin{algs}}
\newcommand{\ealgs}{\end{algs}}
\newcommand{\bass}{\begin{asss}}
\newcommand{\eass}{\end{asss}}

%% file: src_tex/reese_defs.tex
\DeclareMathDelimiter{\vvvert} {0}{matha}{"7E}{mathx}{"17}%

\makeatletter
\def\letterdef#1#2#3{\def\letterdef@##1{\expandafter\def\csname #1\endcsname{#2}}%
  \letterdef@@#3{?\@car{}}\@nil}
\def\letterdef@@#1{\@gobble#1\letterdef@{#1}\letterdef@@}
\makeatother

\DeclarePairedDelimiterX{\klx}[2]{(}{)}{%
  #1\;\delimsize\|\;#2%
}
\DeclarePairedDelimiterX{\quantklx}[3]{(}{)}{%
  #1\;\delimsize\|\;#2\;\delimsize\vert\;#3%
}
\DeclarePairedDelimiterX{\inner}[2]{\langle}{\rangle}{%
  #1,#2%
}

\newcommand{\R}{\mathbf R} %
\newcommand{\N}{\mathbf N} %

\letterdef{c#1}{\mathcal{#1}}{ABCDEFGHIJKLMNOPQRSTUVWXYZ}
\letterdef{b#1}{\mathbb{#1}}{ABCDEFGHIJKLMNOPQRSTUVWXYZ}
\letterdef{bf#1}{\mathbf{#1}}{ABCDEFGHIJKLMNOPQRSTUVWXYZ}
\let\defn\coloneq
\newcommand{\twomax}[2]{\ensuremath{#1 \lor #2}}
\newcommand{\twomin}[2]{\ensuremath{#1 \land #2}}
\newcommand{\threemin}[3]{\ensuremath{#1 \land #2 \land #3}}

\newcommand{\ceil}[1]{\left\lceil #1 \right\rceil}
\newcommand{\floor}[1]{\left\lfloor #1 \right\rfloor}

\newcommand{\e}{\mathrm{e}}
\DeclareMathOperator{\sign}{\bf sign}
\newcommand{\1}{\mathbf 1} %
\let\ones\1
\newcommand{\half}{\frac12} 
\let\epsilon\varepsilon
\newcommand{\eps}{\varepsilon}

\let\tilde\widetilde

\renewcommand{\leq}{\leqslant}
\renewcommand{\geq}{\geqslant}

\newcommand{\argmin}{\mathop{\rm arg\,min}}
\newcommand{\argmax}{\mathop{\rm arg\,max}}

\DeclareMathOperator{\conv}{ conv}

\DeclareMathOperator{\diam}{diam}

\newcommand{\iid}{\textnormal{i.i.d.}}

\newcommand{\simiid}{\stackrel{\iid}{\sim}} %

\newcommand{\Var}{\mathrm{Var}}
\newcommand{\Med}{\mathrm{Med}}

\DeclareMathOperator{\supp}{supp}
\makeatletter
\newcommand{\E}{\operatorname*{\mathbf{E}}\ilimits@}
\makeatother
\makeatletter
\renewcommand{\P}{\operatorname*{\mathbf{P}}\ilimits@}
\makeatother

\newcommand{\ie}{\textit{i}.\textit{e}., }
\newcommand{\eg}{\textit{e}.\textit{g}., }

\newcounter{algorithmctr}
\renewcommand{\thealgorithmctr}{\arabic{algorithmctr}}
   {\refstepcounter{algorithmctr}\begin{list}{}{%
       \setlength{\rightmargin}{0\linewidth}%
       \setlength{\leftmargin}{0\linewidth}}%
       \rmfamily\small
       \item[]{\setlength{\parskip}{0ex}\hrulefill\par%
        \nopagebreak{\bfseries\textsf{Algorithm \thealgorithmctr~}}}}%
   {{\setlength{\parskip}{-1ex}\nopagebreak\par\hrulefill} \end{list}}

\makeatletter
\long\def\@makecaption#1#2{
        \vskip 0.8ex
        \setbox\@tempboxa\hbox{\small {\bf #1.} #2}
        \parindent 1.5em 
        \dimen0=\hsize
        \advance\dimen0 by -3em
        \ifdim \wd\@tempboxa >\dimen0
                \hbox to \hsize{
                        \parindent 0em
                        \hfil 
                        \parbox{\dimen0}{\def\baselinestretch{0.96}\small
                                {\bf #1.} #2
                                } 
                        \hfil}
        \else \hbox to \hsize{\hfil \box\@tempboxa \hfil}
        \fi
        }
\makeatother

%% file: paper_defs.tex
\newcommand{\Normal}[2]{\mathsf{N}\left(#1,#2\right)}
\newcommand{\numobs}{n}
\newcommand{\radius}{r}
\newcommand{\dimension}{d}
\newcommand{\thetastar}{\theta^\star}
\newcommand{\etastar}{\eta^\star}

\let\hat\widehat
\newcommand{\Ball}{\mathsf{B}}
\newcommand{\lamlower}{\underline{\lambda}}
\newcommand{\lambdastar}{\lambda^\star}
\newcommand{\NoiseTermNotation}{\Xi}
\newcommand{\NoiseTerm}[2]{\NoiseTermNotation(#1, #2)}
\DeclareMathOperator{\rad}{rad}
\newcommand{\MLERisk}{\mathcal{R}^{\mathsf{MLE}}}
\newcommand{\MinimaxRisk}{\mathcal{M}}
\newcommand{\LinearMinimaxRisk}{{\MinimaxRisk^{\mathsf{L}}}}
\newcommand{\mustar}{\mu^\star}
\newcommand{\GenCompRatio}{\mathfrak{C}}
\newcommand{\LinCompRatio}{\GenCompRatio^{\mathsf{L}}}
\newcommand{\MLECompRatio}{\GenCompRatio^{\mathsf{MLE}}}
\newcommand{\LowerConst}{c_\ell}
\newcommand{\UpperConst}{c_u}
\newcommand{\GenPackSet}[1]{\mathscr{M}_{#1}}

\usetikzlibrary{positioning, calc} %
\usetikzlibrary{shadings}
\pgfdeclarelayer{background}
\pgfdeclarelayer{main}
\pgfsetlayers{background,main}
\colorlet{lightgreen}{green!20} %
\colorlet{lightred}{red!20}     %
\colorlet{darkerred}{red!40}    %
\usetikzlibrary{arrows.meta, calc}
\usetikzlibrary{decorations.pathreplacing,calligraphy}

%% file: figs/f-plot-p-range.tex
\begin{tikzpicture}[
    region_text/.style={above=10pt, align=center, text width=3.5cm},
    marker/.style={below=2pt}
]

\coordinate (p1) at (2,0);     %
\coordinate (p0) at (4,0);   %
\coordinate (p2) at (12,0);     %
\coordinate (end) at (14,0);   %

\coordinate (top) at (0, 2.5);
\coordinate (bottom) at (0,0);

\begin{pgfonlayer}{background}
    \fill[lightgreen] (p1 |- bottom) rectangle (p0 |- top); %
    \fill[lightred]   (p0 |- bottom) rectangle (p2 |- top); %
    \fill[lightgreen] (p2 |- bottom) rectangle (end |- top); %
\end{pgfonlayer}

\begin{pgfonlayer}{main}
    \draw[->, thick] (p1) -- (end); %

    \draw[dashed] (p1 |- bottom) -- (p1 |- top);
    \draw[dashed] (p0 |- bottom) -- (p0 |- top);
    \draw[dashed] (p2 |- bottom) -- (p2 |- top);

    \draw (p1) -- +(0, 0.1) -- +(0, -0.1) node[marker] {$p=0$};
    \draw (p0) -- +(0, 0.1) -- +(0, -0.1) node[marker, xshift=2pt] {$p=1 + \Theta(\frac{1}{\log d})$};
    \draw (p2) -- +(0, 0.1) -- +(0, -0.1) node[marker] {$p = 2$};

    \node[region_text, yshift=16pt] at ($(p1)!0.5!(p0)$) {\textbf{optimal}\\\scalebox{0.8}{(Thm. 1)} \\ \scalebox{0.8}{for all $\sigma\!>\!0$}};
    \node[region_text, text width=6cm] at ($(p0)!0.5!(p2)$) {\emph{\textbf{suboptimal}} \\[0.1pt] \scalebox{0.8}{(Thm. 2)} \\ \scalebox{0.8}{for all \mbox{$\sigma\!>\!0$} at which} \\ \scalebox{0.8}{nonlinear estimates are required}};
    \node[region_text, yshift=16pt] at ($(p2)!0.5!(end) - (0,0)$) {\textbf{optimal}\\\scalebox{0.8}{(Thm. 1)} \\ \scalebox{0.8}{for all $\sigma\!>\!0$}};
\end{pgfonlayer}
\end{tikzpicture}

%% file: figs/regime1.tex
\begin{tikzpicture}

\definecolor{crimson2143940}{RGB}{214,39,40}
\definecolor{darkgray176}{RGB}{244,244,244}
\definecolor{lightgray204}{RGB}{204,204,204}
\definecolor{steelblue31119180}{RGB}{31,119,180}

\begin{axis}[
width=11cm,
height=6cm,
legend cell align={left},
legend style={
  font=\small,
  fill opacity=0.5,
  text opacity=1,
  at={(0.02,0.02)},
  anchor=south west,
  draw=white
},
log basis x={10},
log basis y={10},
tick align=outside,
tick pos=left,
tick label style={font=\small\fontfamily{phv}\selectfont},
title={},
x grid style={darkgray176},
xlabel={\textbf{dimension,} $d$},
xlabel style = {font=\large\fontfamily{phv}\selectfont, yshift=-5pt},
ylabel style = {font=\large\fontfamily{phv}\selectfont, yshift=-5pt},
xmajorgrids,
xmin=70.7945784384138, xmax=141253.754462276,
xminorgrids,
xmode=log,
xtick style={color=black},
y grid style={darkgray176},
ylabel={\textbf{MSE} },
ylabel style={yshift=1em},
ymajorgrids,
ymin=0.0114391288080682, ymax=0.73731565967566,
yminorgrids,
ymode=log,
ytick style={color=black}
]
\addplot [very thick, black, dashed]
table {%
100 0.595525713383686
115 0.584867655130201
132 0.574480349160226
152 0.563983669946898
175 0.553635562486526
202 0.54323936371374
232 0.533342706296458
268 0.523176559661918
308 0.513512087217998
355 0.503787845583252
409 0.494235692357323
471 0.484856772106722
542 0.475667938774641
625 0.46648698360037
719 0.457600112209861
828 0.448788058065795
954 0.440086468772848
1098 0.431590070955334
1264 0.423219569731457
1456 0.414950809185999
1676 0.406859601635295
1930 0.398881901838981
2222 0.391051600182442
2559 0.383337422457727
2947 0.375758108860291
3393 0.368323037241937
3906 0.361024158548853
4498 0.353837153373619
5179 0.346784228397081
5963 0.339857568232149
6866 0.333052578221546
7906 0.3263685065822
9102 0.319813123341005
10481 0.313367876768846
12067 0.307047558025294
13894 0.300840065247218
15998 0.294746040273534
18420 0.288766337558213
21209 0.28289726000878
24420 0.277138319092376
28117 0.271487599919473
32374 0.26594291801568
37275 0.260503573688913
42919 0.255166481041423
49417 0.249931428000998
56898 0.244796874739734
65512 0.239760437564759
75431 0.234820398987937
86851 0.229975934020141
100000 0.225225137309457
};
\addlegendentry{\footnotesize\fontfamily{phv}\selectfont\textbf{Minimax} ($\sim d^{-0.17}$)}
\addplot [very thick, steelblue31119180, mark=*, mark size=0.5, mark options={solid}]
table {%
100 0.595525713383686
115 0.601157140875118
132 0.582391157469649
152 0.61011871078365
175 0.599779089135816
202 0.567136851200198
232 0.566627981368853
268 0.575124465032347
308 0.572871924222602
355 0.564885298459322
409 0.559590037245891
471 0.563842077637731
542 0.554759385874529
625 0.556599710761881
719 0.540918189599257
828 0.546427017108037
954 0.542638416554523
1098 0.534074854011336
1264 0.547544019756422
1456 0.536825653972472
1676 0.529243151715457
1930 0.530646581203066
2222 0.528783848440666
2559 0.529646051806132
2947 0.527674711242927
3393 0.524242387027886
3906 0.527307163884494
4498 0.5192163833376
5179 0.515429900110244
5963 0.521221022124756
6866 0.519650340322638
7906 0.509965548904371
9102 0.512739712215524
10481 0.513504974572512
12067 0.508595406778082
13894 0.507667634193879
15998 0.511156482513823
18420 0.506016390019134
21209 0.50228239952849
24420 0.504760946630937
28117 0.499967514175025
32374 0.498795164379763
37275 0.498245540756526
42919 0.501821312749135
49417 0.499522731503283
56898 0.496542731882141
65512 0.496440411727036
75431 0.496160217854316
86851 0.49336783806337
100000 0.494386683646553
};
\addlegendentry{\footnotesize\fontfamily{phv}\selectfont\textbf{MLE} ($\sim d^{-0.028}$)}
\addplot [very thick, crimson2143940, mark=*, mark size=0.5, mark options={solid}]
table {%
100 0.595525713383686
115 0.586591672549717
132 0.525333983095409
152 0.544656126219934
175 0.484228666090014
202 0.391000599865661
232 0.357460527843869
268 0.375985970015097
308 0.344968353860879
355 0.306485376730825
409 0.288133963689376
471 0.272454739040144
542 0.246224761808094
625 0.239313911298316
719 0.194604110805993
828 0.19700300022437
954 0.17923168340736
1098 0.156838438356641
1264 0.168215799712298
1456 0.144075066418947
1676 0.129490899294791
1930 0.121492166324176
2222 0.113209974567542
2559 0.106131419315463
2947 0.098329975337188
3393 0.0924145849361459
3906 0.0888160217286146
4498 0.0756461723987349
5179 0.0681500166875806
5963 0.0698827232228697
6866 0.0649792351521213
7906 0.0529808175020764
9102 0.053309542960316
10481 0.0503775286476621
12067 0.0442718826036737
13894 0.0409391436259789
15998 0.0404350210538159
18420 0.0347070504394705
21209 0.031128519864828
24420 0.031367172751336
28117 0.0268892726352623
32374 0.0244139526790781
37275 0.0222317283467106
42919 0.0237139476127896
49417 0.0217644090680427
56898 0.0183312534122777
65512 0.0173681493808635
75431 0.0166367927926198
86851 0.0143305102106632
100000 0.0138239471338333
};
\addlegendentry{\footnotesize\fontfamily{phv}\selectfont \textbf{ST} ($\sim d^{-0.551}$)}
\end{axis}

\end{tikzpicture}

%% file: figs/regime2.tex
\begin{tikzpicture}

\definecolor{crimson2143940}{RGB}{214,39,40}
\definecolor{darkgray176}{RGB}{244,244,244}
\definecolor{lightgray204}{RGB}{204,204,204}
\definecolor{steelblue31119180}{RGB}{31,119,180}

\begin{axis}[
width=11cm,
height=6cm,
legend cell align={left},
legend style={
    font=\small,
    at={(0.02,0.02)}, 
    anchor=south west,
    fill opacity=0.5,
    draw opacity=1,
    text opacity=1,
    draw=white
},
log basis x={10},
log basis y={10},
minor ytick={0.002,0.003,0.004,0.005,0.006,0.007,0.008,0.009,0.02,0.03,0.04,0.05,0.06,0.07,0.08,0.09,0.2,0.3,0.4,0.5,0.6,0.7,0.8,0.9,2,3,4,5,6,7,8,9,20,30,40,50,60,70,80,90},
tick align=outside,
tick pos=left,
tick label style={font=\small\fontfamily{phv}\selectfont},
title={},
x grid style={darkgray176},
xlabel={\textbf{dimension,} $d$},
xlabel style = {font=\large\fontfamily{phv}\selectfont, yshift=-5pt},
ylabel style = {font=\large\fontfamily{phv}\selectfont},
xmajorgrids,
xmin=70.7945784384138, xmax=141253.754462276,
xminorgrids,
xmode=log,
xtick style={color=black},
y grid style={darkgray176},
ylabel={\textbf{MSE} },
ylabel style={yshift=1em},
ymajorgrids,
ymin=0.0233374240897334, ymax=0.362460314013147,
yminorgrids,
ymode=log,
ytick style={color=black}
]
\addplot [very thick, black, dashed]
table {%
100 0.31997398117045
115 0.31023404084628
132 0.30089766017007
152 0.291616595913937
175 0.282615429468519
202 0.273718513481921
232 0.265383100164093
268 0.256954912965964
308 0.249067145429686
355 0.24125165027302
409 0.23369142688683
471 0.226379990244888
542 0.219323097260073
625 0.212376536744107
719 0.205751042305313
828 0.199276421473659
954 0.192975331865719
1098 0.186910921744668
1264 0.181021089898518
1456 0.175285065534671
1676 0.169751028692648
1930 0.164370722949265
2222 0.159163113749485
2559 0.154103581521094
2947 0.149200873643187
3393 0.144457231723113
3906 0.139863779905631
4498 0.135401983219229
5179 0.131082504516512
5963 0.126897328927149
6866 0.122840679078604
7906 0.118909237646023
9102 0.115104686952863
10481 0.111413559202846
12067 0.107841724929721
13894 0.104379767564221
15998 0.101025646476768
18420 0.0977774642055396
21209 0.0946309208148818
24420 0.091583550236151
28117 0.0886322007923558
32374 0.0857736738254786
37275 0.0830056165754677
42919 0.0803245388170277
49417 0.0777284792965672
56898 0.075214867803141
65512 0.0727807978894204
75431 0.0704237665237195
86851 0.0681417572247619
100000 0.0659322999002173
};
\addlegendentry{\footnotesize\fontfamily{phv}\selectfont\textbf{Minimax} ($\sim d^{-0.229}$)}
\addplot [very thick, steelblue31119180, mark=*, mark size=0.5, mark options={solid}]
table {%
100 0.31997398117045
115 0.307772477696195
132 0.289049528256852
152 0.288138814772184
175 0.284221562639553
202 0.269595342933159
232 0.268206056158049
268 0.259975980343397
308 0.257160441177201
355 0.246468230379366
409 0.240362748986446
471 0.237164741611631
542 0.229627392944305
625 0.218972705987596
719 0.20932910436862
828 0.21117899662182
954 0.204774693503707
1098 0.198826043721859
1264 0.191786209015072
1456 0.187318416481389
1676 0.183122468253171
1930 0.177645394451254
2222 0.172783357708577
2559 0.167714733310545
2947 0.164579494003603
3393 0.159788878170334
3906 0.154819233287238
4498 0.150097754849221
5179 0.14692164203947
5963 0.143373740440105
6866 0.138680486524588
7906 0.135228003188236
9102 0.132140975030623
10481 0.128186962190243
12067 0.125216350463053
13894 0.122310536274567
15998 0.118561929944013
18420 0.114997073438953
21209 0.11250653059053
24420 0.109153907214767
28117 0.106063872478493
32374 0.103280136128234
37275 0.100905048852507
42919 0.0983094435521919
49417 0.0953862632759192
56898 0.0929453518957151
65512 0.0904588231975573
75431 0.088280968471367
86851 0.0858181039202139
100000 0.0834945002524283
};
\addlegendentry{\footnotesize\fontfamily{phv}\selectfont\textbf{MLE} ($\sim d^{-0.192}$)}
\addplot [very thick, crimson2143940, mark=*, mark size=0.5, mark options={solid}]
table {%
100 0.31997398117045
115 0.31281561688389
132 0.294282083697571
152 0.286692573164848
175 0.279226669370361
202 0.264787868509475
232 0.256987530597772
268 0.248268639898539
308 0.243953822320737
355 0.22885717951615
409 0.224051481894983
471 0.215658123571
542 0.208879195031478
625 0.192136848315425
719 0.182073615031274
828 0.181433143678617
954 0.173990443787658
1098 0.16559061949132
1264 0.154947396668053
1456 0.149333190216948
1676 0.143021884695658
1930 0.13491863003968
2222 0.129421985257395
2559 0.121572592639079
2947 0.11771528927141
3393 0.111333619398457
3906 0.10417239714135
4498 0.0982948323788376
5179 0.0929752731218925
5963 0.0891099058363869
6866 0.0826180760616568
7906 0.0785431152426516
9102 0.0750516333995927
10481 0.0703283076427878
12067 0.0665057829132615
13894 0.0631370970444376
15998 0.0590044145141272
18420 0.0549763081569278
21209 0.0525342829835908
24420 0.0491426768403424
28117 0.0459604805264994
32374 0.0430280660743636
37275 0.04113616777212
42919 0.0386898635280456
49417 0.0360862974545899
56898 0.0339533217447784
65512 0.0317734353541834
75431 0.0301261499668354
86851 0.028195088748658
100000 0.0264361809447146
};
\addlegendentry{\footnotesize\fontfamily{phv}\selectfont\textbf{ST} ($\sim d^{-0.365}$)}
\end{axis}
\end{tikzpicture}

%% file: figs/f-prop1-fig.tex
\begin{tikzpicture}[>=Stealth]

\def\p{1.5} 
\def\r{2.0}
\def\sig{1.6}
\def\angleP{38}  

\draw[->, gray] (-2.5, 0) -- (3.0, 0);
\draw[->, gray] (0, -2.5) -- (0, 2.75);

\draw[
    color=RoyalBlue, %
    thick,           %
    variable=\t,
    domain=0:360,
    samples=200      %
] plot ({ \r * sign(cos(\t)) * abs(cos(\t))^(2/\p) },
        { \r * sign(sin(\t)) * abs(sin(\t))^(2/\p) }) 
 node [below right, xshift=-35pt, yshift=-45pt]{\scalebox{0.65}{$\mathsf{B}^d_p$}};

\coordinate (O) at (0,0); %
\coordinate (P) at (0.59060626, 1.7798624);
\coordinate (Px) at (P |- O); %
\coordinate (T) at (\r, 0); 
\coordinate (N) at ({-\r * \sig * 0.3}, 
{\r * \sig * 0.7});
\coordinate (Y) at ($(T) + (N)$);
\node at (Y) [circle,fill,inner sep=1.1pt]{};

\draw[->, thin, color=Maroon!50!Gray] (O) -- (N);
\draw[thick, dashed, color=Black] (P) -- (Y) node[below right, xshift=12pt, yshift=-20pt]{\scalebox{0.6}{\color{Maroon!50!Gray}{$+\sigma \xi$}}};
\draw[thin, dashed, color=gray!75!white] (P) -- (Px);
\draw[->, thin, dashed, color=Maroon!50!Gray] (T) -- (Y);
\node[anchor=south west] at (Y) {\scalebox{0.8}{$Y$}};
\node[anchor=south west, xshift=-13pt, yshift=-2pt] at (N) {\scalebox{0.8}{{\color{Maroon!75!Gray}{$\sigma \xi$}}}};
\node[anchor=south west, xshift=-20pt, yshift=1pt] at (P) {\scalebox{0.7}{$\Pi_{\mathsf{B}^d_p}(Y)$}};
\node[anchor=west, xshift=-1pt, yshift=5pt] at (T) {\scalebox{0.8}{$\theta^\star$}};
\node at (P) [circle,fill,inner sep=1.1pt]{};
\node at (T) [circle,fill,inner sep=1.1pt]{};

\draw [decorate, thick, 
    decoration = {calligraphic brace, mirror, raise=2pt}] (Px) --  (T) 
    node [below left, xshift=-5pt, yshift=-3pt] {\scalebox{0.6}{$ \geq 1/2$}};
\end{tikzpicture}